\documentclass[11pt]{amsart}
\usepackage[T1]{fontenc}
\usepackage{lmodern,amsmath,amssymb,mathtools,bm,mathrsfs,booktabs,array,longtable}
\usepackage[margin=27mm]{geometry}
\usepackage{microtype,xcolor,tikz,pgfplots,pgfplotstable,graphicx}
\usetikzlibrary{arrows.meta,positioning,calc}
\pgfplotsset{compat=1.18}
\usepackage[colorlinks=true,linkcolor=blue!45!black,citecolor=blue!45!black,urlcolor=blue!45!black]{hyperref}
\usepackage{bookmark,placeins}
\hypersetup{pdftitle={Pure Spinor Factorization of Belinski--Zakharov Solitons and Regular Gravitational Instantons},pdfauthor={Steffen Aksteiner}}
\allowdisplaybreaks[2]

\newtheorem{theorem}{Theorem}[section]
\newtheorem{proposition}[theorem]{Proposition}
\newtheorem{lemma}[theorem]{Lemma}
\newtheorem{corollary}[theorem]{Corollary}
\theoremstyle{definition}

\newtheorem{assumption}[theorem]{Hypotheses}
\newtheorem{conjecture}[theorem]{Conjecture}
\theoremstyle{remark}
\newtheorem{remark}[theorem]{Remark}

\DeclareMathOperator{\Pf}{Pf}
\DeclareMathOperator{\tr}{tr}
\DeclareMathOperator{\diag}{diag}
\DeclareMathOperator{\adj}{adj}

\newcommand{\ii}{\mathrm i}
\newcommand{\one}{\mathbf 1}

\newcommand{\R}{\mathbb R}

\title[Pure spinor factorization and regular instantons]{Pure Spinor Factorization of Belinski--Zakharov Solitons and Regular Gravitational Instantons}
\author{Steffen Aksteiner}
\thanks{Email: \texttt{steffen.aksteiner@gmail.com}.}
\date{1 November 2026}
\begin{document}
\begin{abstract}
We show that the Belinski--Zakharov soliton construction factorizes into signed Pfaffian tau-functions, revealing a pure-spinor description of the metric and its regularity conditions. We derive analytic restrictions on the integral rod data, prove uniqueness for the asymptotically flat balance problem at fixed rod structure and periods, and give a conditional criterion for global existence. Combining these results with exhaustive enumeration of the admissible data and rigorous existence certificates, we construct a catalogue of regular asymptotically flat solutions representing 356 rod-word classes through twelve turning points, and 52 isolated regular asymptotically locally flat solutions through eight turning points. Analytic positivity and smooth-completion results establish regularity beyond the finite numerical solve. Based on this catalogue, we investigate two organizations of the solution space: generalized rank-two root systems and compositions of simpler endpoint geometries. The endpoint composition rules reproduce exactly all retained asymptotically flat words through twelve turning points.
\end{abstract}
\maketitle
\setcounter{tocdepth}{1}
\tableofcontents
\clearpage

\section{Introduction}

A gravitational instanton is a complete oriented Ricci-flat Riemannian four-manifold with finite curvature energy. For an effective isometric two-torus action, the boundary intervals of the orbit space record the circles that collapse. These intervals are called rods; their endpoints are turning points, or nuts. A primitive integral vector specifies each collapsing circle. The ordered vectors encode topology, whereas their normalization determines the cone angles. The central question is which rod configurations admit smooth Ricci-flat realizations and how the resulting families of gravitational instantons are organized.

The discovery of the Chen--Teo instanton \cite{ChenTeo} showed that the regular asymptotically flat (AF) solutions extend beyond the Euclidean Kerr family, including its Schwarzschild limit. It disproved the proposed Riemannian analogue of black-hole uniqueness and opened the problem of describing a larger moduli space. Li and Sun subsequently proved the existence of an infinite sequence of AF toric instantons with arbitrarily many turning points \cite{LiSun}; beyond the Kerr and Chen--Teo members, these metrics are non-Hermitian. Teo gave an explicit construction of the four-turning-point member, reducing the remaining balance condition to a quintic equation \cite{Teo}. The present work seeks the structure beyond this sequence. We derive restrictions on the possible rod data, construct a certified catalogue of regular solutions, and use it to test two descriptions of the moduli space: one through root arrangements and polygon walks, the other through compositions of endpoint geometries.

The Belinski--Zakharov (BZ) method constructs metrics from finitely many spectral poles and charge parameters \cite{BZ,Chen}. Twistor theory describes the same field equations through a rational matrix on the symmetry axis, called the patching matrix \cite{Tod,DunajskiTod}. Harmonic-map methods instead start from prescribed rods and give a canonical metric, generally with conical singularities \cite{KL,LiSun}. Two recent representation results make the finite algebraic approach particularly significant. 

Dunajski, Mason and Tod give a complete twistor representation of toric gravitational instantons with asymptotically locally flat or asymptotically locally Euclidean ends. They construct the rational patching matrices recursively and recover the associated metrics by linear algebra \cite{DMT}. Andersson, Dahl and Nilsson prove that simply connected toric instantons with AF or asymptotically locally flat (ALF) ends and a finite positive number of turning points admit a flat-seed BZ representation \cite{ADN}. An ALF end has a circle of bounded length at infinity; precise conventions are given below. These results reduce geometric classification to finite data, but the cone-balance conditions still have to be solved.

Our first contribution is a factorization of the positive BZ construction by signed Pfaffian tau-functions. The algebraic starting point is the Cauchy--Pfaffian identity
\[
 \det\left(\frac{1-a_i a_j}{1-b_i b_j}\right)_{i,j=1}^{n}
 =\frac{\tau_n(a,b)\,\tau_n(-a,b)}{\prod_{i=1}^{n}(1-b_i^2)}.
\]
Here $a=(a_i)$ and $b=(b_i)$ are auxiliary variables and $\tau_n$ is the signed Pfaffian function defined in Section~\ref{sec:three}. In the BZ specialization, $b_i$ is a rescaled spectral pole and $a_i$ is its charge multiplied by a power of that pole. We show how the two factors and their neighboring shifts reconstruct the metric, its Ernst potentials and its cone-normalized rods. The rank-one determinant formulas recover Alekseev's reconstruction and its monodromy interpretation \cite{Alekseev1981,Alekseev}; an explicit change of charge variables gives the Neugebauer--Kramer determinant representation \cite{Neugebauer1980,NeugebauerKramer}. The pure-spinor interpretation explains the Pfaffian relations; the positive moment representation makes bulk positivity and axis limits explicit. Polynomial factors of the patching matrix connect these formulas to the boundary description. The closed Hankel formulas of \cite{ADN} give a complementary representation and a smooth-completion theorem; Section~\ref{subsec:ADN-dictionary} supplies the dictionary needed to apply it here.

Our second contribution is a variational description of balance. For a marked rod word $W$, let $\ell_j$ be its positive finite rod lengths, $\theta_j$ its cone angles and $M$ its mass. At fixed asymptotic periods, define
\[
 b_j=\log\frac{2\pi}{\theta_j},\qquad
 \Phi_W(\ell)=\sum_j\ell_jb_j-2M.
\]
We prove that $\nabla_\ell\Phi_W=b$ and that $\Phi_W$ is strictly convex on the positive length orthant. Consequently a balanced metric is unique at fixed marked data and asymptotic scale. The Legendre dual of $\Phi_W$ is $2M$. Global existence also requires exclusion of collision minima. We formulate a descending-opening criterion and prove identities for zero insertions in self-intersection words; the full-interval existence statements retain the explicit collision-matching hypotheses of Section~\ref{subsec:collision-matching}.

The third contribution is a catalogue of regular metrics obtained from analytic restrictions, exhaustive enumeration and rigorous existence certificates. Positive moment identities imply integral determinant bounds and restrict the finite self-intersections to $-1,0,1$ in the nonzero reduced ranks under consideration. The low-order examples include the Chen--Teo and Teo metrics \cite{ChenTeo,Teo}; exact enumeration through twelve genuine turning points gives 359 AF candidates. For 356 of these rod-word classes we prove the existence of a smooth complete AF metric at a specified interior twist, and hence of a local family; three further candidates are excluded by a positive boundary gradient. The ALF census through eight turning points gives 52 certified isolated regular metrics on 49 words. The finite certificates establish exact solutions of the full normalized rod equations. The analytic positive representation and completion theorem establish positivity throughout the bulk, smooth rods and nuts, and completeness of the end. Thus the catalogue consists of proved existence results, rather than metrics accepted on the basis of small residuals or a sampled regularity grid. Full-interval AF existence retains the stated collision hypotheses, and the ALF census is restricted to the enumerated positive reduced sectors.

The catalogues also serve as experimental templates for a classification. The first approach retains the distinct rod directions and identifies them with generalized rank-two root systems, equivalently triangulated polygons \cite{CuntzHeckenberger,CuntzArrangements}. The ordered walk on the directions is essential: the same polygon can support both obstructed and regular words. The second approach composes marked endpoint components using their boundary directions and twist residues. Nine components generate exactly the 356 retained AF words through twelve turning points. This is an exact finite combinatorial result, independent of solving the parent balance equations. It does not yet prove that every permitted assembly opens to a balanced family. The ALF catalogue supplies components for this construction, while the AF catalogue tests their compatibility and the predicted opposite endpoints.

Together with the representation result of \cite{ADN}, these restrictions and certificates provide a route to a geometric classification: one must control all permitted sectors, prove the required collision matching, and identify equivalent metrics. The recursive patching-matrix construction of \cite{DMT} addresses the representation stage; our cone potential, admissibility tests and catalogue address the regular locus within it. Neither recursive representation nor the observed endpoint rules alone establish the missing global existence statement.

\paragraph{Conventions and organization.}
The integer $n$ counts BZ dressing slots and $N$ counts genuine turning points. A phantom slot contributes no nut. An AF end has leading metric
\[
 dt^2+dr^2+r^2(d\vartheta^2+\sin^2\vartheta\,d\phi^2),
\]
with identifications $(t,\phi)\sim(t,\phi+2\pi)\sim(t+2\pi L,\phi+2\pi\beta)$ and the corresponding derivative decay. We fix $L>0$ and use $-1<\beta<0$ for the normalized nontrivial AF words. ALF is understood in the finite-cover sense specified in Appendix~\ref{app:ALF-system}. No reflection symmetry is imposed.

Sections~\ref{sec:three}--\ref{sec:patching} develop the algebraic representations; Sections~\ref{sec:regularity}--\ref{sec:analytic} treat admissibility and balance. Sections~\ref{sec:classification-roots} and~\ref{sec:endpoint-composition} present the two classification approaches consecutively. Appendix~\ref{app:numerics} contains the AF equations, explicit low-order examples and catalogue. Appendix~\ref{app:ALF-system} develops the ALF marking and records its certified solutions and endpoint components; the final appendix supplies supplementary analytic arguments.

\section{Three levels of the soliton construction}\label{sec:three}

The inverse-scattering, determinant and tau descriptions encode the same
solution. We fix their normalizations together, including the different
flat seeds for even and odd dressing. Throughout this section, $n$ denotes
the number of BZ slots; $N$ denotes the number of genuine turning points.
These numbers agree unless some slots become phantom solitons. No
asymptotic-flatness or rod-regularity condition is imposed here. Transpose is denoted by $T$ as a superscript, $I$ denotes the identity matrix of the relevant size, $\one$ is the all-ones column, $\adj$ is the classical adjugate, and $\ii^2=-1$. Norms of real coordinate vectors are Euclidean.

\subsection{The original BZ method}

Let $G$ be the symmetric $2\times2$ Gram matrix of the commuting Killing fields and let $\nu$ be the orbit-space conformal factor. Write the metric in Weyl--Papapetrou coordinates as
\begin{equation}\label{eq:WP}
 g=G_{ab}\,dx^a dx^b+e^{2\nu}(d\rho^2+dz^2),\qquad
 (x^1,x^2)=(t,\phi),\qquad \det G=\rho^2.
\end{equation}
On the open orbit space $\rho>0$, introduce the right currents
\begin{equation}
 U=\rho G_\rho G^{-1},\qquad V=\rho G_zG^{-1}.
\end{equation}
The vacuum equations are
\begin{align}
 \partial_\rho U+\partial_z V&=0,\label{eq:vacuum}\\
 \nu_\rho&=-\frac1{2\rho}
       +\frac{\tr(U^2-V^2)}{8\rho},\qquad
 \nu_z=\frac{\tr(UV)}{4\rho}.\label{eq:nu}
\end{align}
For positive definite $G$, the first equation is the harmonic-map equation
for $\Phi=\rho^{-1}G$ with values in the positive definite symmetric
matrices of determinant one. It also ensures compatibility of the two
equations for $\nu$.

\paragraph{The spectral problem.}
Let $\lambda_{\mathrm{sp}}$ be the dimensionful BZ spectral parameter; it is distinct from the finite rod lengths introduced later. The linear system is
\begin{align}
 \mathscr D_\rho\Psi
 &=\frac{\rho U+\lambda_{\mathrm{sp}} V}{\rho^2+\lambda_{\mathrm{sp}}^2}\Psi,&
 \mathscr D_\rho
 &=\partial_\rho+
       \frac{2\rho\lambda_{\mathrm{sp}}}{\rho^2+\lambda_{\mathrm{sp}}^2}\partial_{\lambda_{\mathrm{sp}}},
       \label{eq:bz-lax-original}\\
 \mathscr D_z\Psi
 &=\frac{\rho V-\lambda_{\mathrm{sp}} U}{\rho^2+\lambda_{\mathrm{sp}}^2}\Psi,&
 \mathscr D_z
 &=\partial_z-
       \frac{2\lambda_{\mathrm{sp}}^2}{\rho^2+\lambda_{\mathrm{sp}}^2}\partial_{\lambda_{\mathrm{sp}}}.
       \notag
\end{align}
Its compatibility gives \eqref{eq:vacuum}; the generating matrix is
normalized by $\Psi(0)=G$. The quantity
\begin{equation}\label{eq:bz-spectral-cover}
 w=z+\frac{\lambda_{\mathrm{sp}}^2-\rho^2}{2\lambda_{\mathrm{sp}}}
\end{equation}
is constant under both $\mathscr D_\rho$ and $\mathscr D_z$.
Thus $\lambda_{\mathrm{sp}}$ uniformizes a two-sheeted spectral cover with sheet involution
$\lambda_{\mathrm{sp}}\mapsto-\rho^2/\lambda_{\mathrm{sp}}$.

The dimensionless coordinate used below is
\begin{equation}\label{eq:bz-zeta}
 \zeta=\frac{\ii\lambda_{\mathrm{sp}}}{\rho},\qquad
 w=z-\frac{\ii\rho}{2}(\zeta+\zeta^{-1}).
\end{equation}
The sheet involution becomes $\zeta\mapsto\zeta^{-1}$, and the linear
system takes the form
\begin{align}
 \widehat{\mathscr D}_\rho\Psi
 &=\frac{U-\ii\zeta V}{\rho(1-\zeta^2)}\Psi,&
 \widehat{\mathscr D}_\rho
 &=\partial_\rho+
   \frac{\zeta(1+\zeta^2)}{\rho(1-\zeta^2)}\partial_\zeta,
   \label{eq:bz-lax-zeta}\\
 \widehat{\mathscr D}_z\Psi
 &=\frac{V+\ii\zeta U}{\rho(1-\zeta^2)}\Psi,&
 \widehat{\mathscr D}_z
 &=\partial_z+
   \frac{2\ii\zeta^2}{\rho(1-\zeta^2)}\partial_\zeta.
   \notag
\end{align}
These operators differentiate at fixed $w$, rather than at fixed $\zeta$.

\paragraph{Dressing the seed.}
Given a seed generating matrix $\Psi_0$, the BZ construction
\cite{BZ,Chen} sets
\begin{equation}\label{eq:bz-dressing}
 \Psi_{\mathrm{un}}(\lambda_{\mathrm{sp}})=\mathscr X(\lambda_{\mathrm{sp}})\Psi_0(\lambda_{\mathrm{sp}}),\qquad
 \mathscr X(\lambda_{\mathrm{sp}})=I+\sum_{i=1}^n\frac{R_i}{\lambda_{\mathrm{sp}}-\mu_i},\qquad
 R_i=\mathbf n_i m_i.
\end{equation}
Here $m_i$ is a row and $\mathbf n_i$ is a column. The moving poles lie
over fixed spectral positions $w=z_i$. We use the branch
\begin{equation}\label{eq:mu}
 \mu_i=\sqrt{\rho^2+(z-z_i)^2}-(z-z_i),\qquad
 \mu_i^2+2(z-z_i)\mu_i-\rho^2=0.
\end{equation}
The second branch is $-\rho^2/\mu_i$. For the real slice considered
later, the $z_i$ are distinct real numbers.

Each pole also carries a constant, nonzero polarization row $m_{0,i}$,
defined up to multiplication by a nonzero scalar. Set
\begin{equation}\label{eq:bz-polarizations}
 m_i=m_{0,i}\Psi_0(\mu_i)^{-1},\qquad
 N_i=m_iG_0,\qquad
 \Gamma_{ij}=\frac{m_iG_0m_j^T}{\rho^2+\mu_i\mu_j}.
\end{equation}
The projective polarization coordinate will be denoted by $C_i$ in the
seed conventions below. We call $C_i$ a BZ charge; it is not itself a
mass, angular momentum or NUT charge. Infinite $C_i$ is a projective
limit, obtained by using the other polarization chart.

The symmetric-matrix reduction is
\begin{equation}\label{eq:bz-reduction}
 G_{\mathrm{un}}
 =\mathscr X(\lambda_{\mathrm{sp}})G_0\mathscr X^T(-\rho^2/\lambda_{\mathrm{sp}}).
\end{equation}
Cancellation of its residue at $\lambda_{\mathrm{sp}}=\mu_i$ gives the finite linear
system
\begin{equation}\label{eq:bz-residue-system}
 \sum_j\Gamma_{ij}\mathbf n_j=\frac{N_i^T}{\mu_i},\qquad
 \mathbf n_i=\sum_j(\Gamma^{-1})_{ij}\frac{N_j^T}{\mu_j}.
\end{equation}
Evaluating \eqref{eq:bz-dressing} at $\lambda_{\mathrm{sp}}=0$ therefore gives
\begin{equation}\label{eq:BZ}
 G_{\mathrm{un}}=G_0-\sum_{i,j=1}^n
 \frac{(\Gamma^{-1})_{ij}}{\mu_i\mu_j}N_i^TN_j.
\end{equation}
It solves the reduced matrix equation, but generally does not yet have
determinant $\rho^2$. In fact,
\begin{equation}\label{eq:bz-det-un}
 \det G_{\mathrm{un}}
 =(-1)^n\frac{\rho^{2n}}{\prod_i\mu_i^2}\det G_0.
\end{equation}

\paragraph{Even and odd flat-seed conventions.}
Put $Q(\lambda_{\mathrm{sp}})=\rho^2-2z\lambda_{\mathrm{sp}}-\lambda_{\mathrm{sp}}^2$. We use
\begin{equation}\label{eq:bz-seeds}
\begin{array}{c|c|c|c|c}
 &G_0&\Psi_0(\lambda_{\mathrm{sp}})&m_i&N_i\\ \hline
 n\ \mathrm{even}
 &\diag(1,\rho^2)&\diag(1,Q(\lambda_{\mathrm{sp}}))
 &(C_i,\mu_i^{-1})&(C_i,\rho^2\mu_i^{-1})\\[2pt]
 n\ \mathrm{odd}
 &\diag(-1,1)&\diag(-1,1)
 &(-C_i,1)&(C_i,1)
\end{array}
\end{equation}
For even $n$, the displayed rows follow from $m_{0,i}=(C_i,-2z_i)$,
using $Q(\mu_i)=-2z_i\mu_i$; a shift of the spectral origin avoids
$z_i=0$ in this chart, and the resulting expressions extend by limits.
For odd $n$, take $m_{0,i}=(C_i,1)$. The odd seed is an indefinite
generalized seed of the reduced matrix equation, with determinant $-1$;
it is not a positive definite Gram matrix satisfying \eqref{eq:WP}.
The odd dressing and determinant normalization produce the Riemannian
section.

Define
\begin{equation}\label{eq:lambda}
 \lambda_i=\frac{\ii\mu_i}{\rho},\qquad
 \Lambda=\prod_i\lambda_i,\qquad
 \lambda_i+\lambda_i^{-1}=-\frac{2\ii(z-z_i)}{\rho}.
\end{equation}
Thus $\lambda_i$ is the value of $\zeta$ at the $i$th dressing pole.
The determinant normalizations, with their fixed signs, are
\begin{equation}\label{eq:bz-normalization}
 G=h_nG_{\mathrm{un}},\qquad
 h_n=\begin{cases}
       \Lambda,&n\ \mathrm{even},\\
       \ii\rho\Lambda,&n\ \mathrm{odd}.
      \end{cases}
\end{equation}
They give $\det G=\rho^2$. Positivity of $G$ is a further condition on
the real branch. Multiplication by $h_n$ adds scalar matrices to the
currents, so the off-diagonal currents are unchanged.

The conformal factor obtained from \eqref{eq:nu} is
\begin{equation}\label{eq:ffBZ}
 e^{2\nu}=k_0\det\Gamma
 \prod_{i<j}(\rho^2+\mu_i\mu_j)^2
 \begin{cases}
 \displaystyle\rho^{-n^2/2}\Bigl(\prod_i\mu_i\Bigr)^{3-n},
       &n\ \mathrm{even},\\[3pt]
 \displaystyle\rho^{-(n-1)^2/2}\Bigl(\prod_i\mu_i\Bigr)^{2-n},
       &n\ \mathrm{odd}.
 \end{cases}
\end{equation}
The constant $k_0$ fixes the remaining integration freedom and enters
the normalization of the rods. Finally, differentiation at fixed $z_i$
gives
\begin{equation}\label{eq:lambdaD}
 \partial_\rho\lambda_i
 =\frac{\lambda_i(1+\lambda_i^2)}{\rho(1-\lambda_i^2)},\qquad
 \partial_z\lambda_i
 =\frac{2\ii\lambda_i^2}{\rho(1-\lambda_i^2)}.
\end{equation}

\subsection{The determinant representation}

For an integer $\alpha$, define
\begin{equation}\label{eq:M}
 M^{(\alpha)}_{ij}
 =\frac{1-C_iC_j\lambda_i^\alpha\lambda_j^\alpha}
        {1-\lambda_i\lambda_j},\qquad
 d_\alpha=\det M^{(\alpha)}.
\end{equation}
Write $a^{(\alpha)}=C\lambda^\alpha$,
$D_\lambda=\diag(\lambda_i)$ and $D_\mu=\diag(\mu_i)$; vector
products and powers are componentwise. The two elementary identities
\begin{equation}\label{eq:Mupdate}
 M^{(\alpha+1)}=M^{(\alpha)}+a^{(\alpha)}a^{(\alpha)T}
 =\one\one^T+D_\lambda M^{(\alpha)}D_\lambda
\end{equation}
relate all neighboring determinants. The original BZ matrix is
\begin{equation}\label{eq:Gamma-M}
 \Gamma=\begin{cases}
 D_\mu^{-1}M^{(1)}D_\mu^{-1},&n\ \mathrm{even},\\
 \rho^{-2}M^{(0)},&n\ \mathrm{odd}.
 \end{cases}
\end{equation}
Thus the determinant center is $1$ for even dressing and $0$ for odd
dressing.

To include the off-diagonal entry, set
\begin{equation}\label{eq:det-offdiag}
 r_\alpha
 =\one^T(M^{(\alpha)})^{-1}a^{(\alpha)}
 =\frac{\det(M^{(\alpha)}+a^{(\alpha)}\one^T)}{d_\alpha}-1.
\end{equation}
Substitution into \eqref{eq:BZ}, followed by the determinant lemma, gives
\begin{align}
 G&=\begin{pmatrix}
 \Lambda d_0/d_1&-\ii\rho r_1/\Lambda\\
 -\ii\rho r_1/\Lambda&\rho^2d_2/(\Lambda d_1)
 \end{pmatrix},&&n\ \mathrm{even},\label{eq:Gdet}\\[3pt]
 G&=\begin{pmatrix}
 -\ii\rho\Lambda d_{-1}/d_0&\ii\rho r_0/\Lambda\\
 \ii\rho r_0/\Lambda&\ii\rho d_1/(\Lambda d_0)
 \end{pmatrix},&&n\ \mathrm{odd}.\label{eq:Gdet-odd}
\end{align}
For comparison, the even off-diagonal entry can also be written as
\begin{equation}
 G_{t\phi}=-\ii\rho\Lambda
 \left\{\frac{\det(M^{(1)}+C(\lambda^{-1})^T)}{d_1}-1\right\}.
\end{equation}
The equivalence follows from \eqref{eq:Mupdate}. All expressions are
understood on nonsingular charts, with removable singularities treated
by taking limits of the complete quantities.

In these variables, the conformal factor is
\begin{equation}\label{eq:nudet}
 e^{2\nu}=k_0\prod_{i<j}(1-\lambda_i\lambda_j)^2
 \begin{cases}
 \displaystyle(-1)^{n/2}\rho^{(n^2-2n)/2}
          \Lambda^{1-n}d_1,&n\ \mathrm{even},\\[3pt]
 \displaystyle-\ii\rho^{(n^2-2n-1)/2}
          \Lambda^{2-n}d_0,&n\ \mathrm{odd}.
 \end{cases}
\end{equation}
For real charges and imaginary $\lambda_i$, the complete expressions
are real. In particular, the odd prefactor $-\ii$ is fixed by
\eqref{eq:ffBZ} and \eqref{eq:lambda}; it cannot be replaced independently
by a real sign.

The determinant formulas now reduce to signed tau products. The diagonal entries use the common denominator and its neighboring shifts; the off-diagonal entry uses the rank-one correction in \eqref{eq:det-offdiag}.

\subsection{The tau representation}

Let $a=(a_1,\ldots,a_n)$ and $b=(b_1,\ldots,b_n)$ be charge-weight and spectral-coordinate vectors. For ordered indices $[n]=\{1,\ldots,n\}$, put
\begin{equation}
 K(b_i,b_j)=\frac{b_i-b_j}{1-b_ib_j},\qquad
 \Delta_b(S)=\prod_{\substack{i<j\\i,j\in S}}K(b_i,b_j).
\end{equation}
Empty products are one. Define the signed subset sum
\begin{equation}\label{eq:tau}
 \tau_n(a,b)=\sum_{S\subset[n]}(-1)^{\sum_{i\in S}i}
 \left(\prod_{i\in S}a_i\right)\Delta_b(S)\Delta_b(S^c),\qquad
 \tau_\pm^{(\alpha)}=\tau_n(\pm C\lambda^\alpha,\lambda).
\end{equation}
The index sign is part of the definition. In particular, inserting or
removing a slot changes the ordering signs of subsequent slots. Set
\begin{equation}\label{eq:Pi}
 \Pi=\prod_i(1-\lambda_i^2).
\end{equation}

\begin{proposition}[Determinant factorization and symmetric shift]
\label{prop:factor}
For either parity of $n$,
\begin{align}
 d_\alpha&=\frac{\tau_+^{(\alpha)}\tau_-^{(\alpha)}}{\Pi},
          \label{eq:factor}\\
 \tau_+^{(\alpha+1)}\tau_-^{(\alpha-1)}
 +\tau_+^{(\alpha-1)}\tau_-^{(\alpha+1)}
 &=2\tau_+^{(\alpha)}\tau_-^{(\alpha)}.
          \label{eq:bilinear}
\end{align}
\end{proposition}
\begin{proof}
The Cauchy--Pfaffian factorization for
$((1-a_ia_j)/(1-b_ib_j))$ gives the product of the two signed subset
sums in \eqref{eq:tau}, divided by $\prod_i(1-b_i^2)$
\cite{CauchyPfaffian}. Its neighboring-weight identity gives
\eqref{eq:bilinear}. The even and bordered odd Pfaffian representatives
are described in Section~\ref{sec:geometry}. Equivalently, both identities
follow by expanding in the $a_i$ and grouping complementary subsets.
After clearing denominators, they extend across the singular charts.
\end{proof}
For example,
$\tau_1(a,b)=1-a_1$ and
$\tau_2(a,b)=(1-a_1a_2)K(b_1,b_2)+a_2-a_1$.

\begin{proposition}[Antisymmetric shift and parity]
\label{prop:antisymmetric-shift}
With the same conventions,
\begin{align}
 &\tau_+^{(\alpha+1)}\tau_-^{(\alpha-1)}
 -\tau_+^{(\alpha-1)}\tau_-^{(\alpha+1)}
 \notag\\[-2pt]
 &\hspace{12mm}
 =-2(-1)^n\frac{\Pi}{\Lambda}
       \one^T\adj(M^{(\alpha)})a^{(\alpha)}.
       \label{eq:antisymmetric-shift}
\end{align}
When $d_\alpha\ne0$, the right-hand side is
$-2(-1)^n\Pi d_\alpha r_\alpha/\Lambda$.
\end{proposition}
\begin{proof}
Work first where $M=M^{(\alpha)}$ is invertible, and write
$a=a^{(\alpha)}$, $d=d_\alpha$ and $S=M^{-1}$. The two rank-one updates give
\begin{equation}
 d_{\alpha+1}=d(1+a^TSa),\qquad
 d_{\alpha-1}=\frac{d}{\Lambda^2}(1-\one^TS\one).
\end{equation}
Since $M+aa^T-\one\one^T=D_\lambda M D_\lambda$, the rank-two
determinant lemma also gives
\begin{equation}
 (1+a^TSa)(1-\one^TS\one)+(\one^TSa)^2=\Lambda^2.
\end{equation}
Consequently,
\begin{equation}
 d_\alpha^2-d_{\alpha+1}d_{\alpha-1}
       =\frac{d_\alpha^2}{\Lambda^2}(\one^TSa)^2.
\end{equation}
By Proposition~\ref{prop:factor}, the square of the left-hand side of
\eqref{eq:antisymmetric-shift} is therefore the square of its proposed
right-hand side. The sign is fixed by setting $a_1=\cdots=a_{n-1}=0$.
Writing $b=\lambda$ and $D_* =\Delta_b([n])\Delta_b([n-1])$, the
antisymmetric tau product is then
$2(-1)^na_n(b_n-b_n^{-1})D_*$, while the Cauchy cofactor identity gives
\begin{equation}
 \one^T\adj(M)a
       =\frac{\Lambda}{\Pi}a_n(b_n^{-1}-b_n)D_*.
\end{equation}
This proves the stated sign. Rational continuation gives the adjugate
identity without the invertibility assumption.
\end{proof}
Together with the symmetric shift, the proposition also gives the
polarized determinant identity
\begin{equation}\label{eq:polarized-Cauchy}
 \tau_+^{(\alpha+1)}\tau_-^{(\alpha-1)}
 =\Pi\det\left(M^{(\alpha)}
       -\frac{(-1)^n}{\Lambda}a^{(\alpha)}\one^T\right).
\end{equation}
Reversing all charges reverses the rank-one term and exchanges the two
shifted products.

\begin{remark}[Alekseev's determinant reconstruction]\label{rem:Alekseev}
Alekseev's determinant form \cite{Alekseev1981} is the BZ residue formula with each inverse-matrix contraction replaced by a rank-one determinant quotient. To make the comparison explicit in our conventions, put $\xi_a=(N_{ia}/\mu_i)_{i=1}^n$ for each Killing-coordinate index $a$ and choose nonzero scalars $c_a$. The determinant lemma rewrites \eqref{eq:BZ} as
\begin{equation}\label{eq:Alekseev-determinant}
 G_{ab}=h_n\left\{(G_0)_{ab}+c_ac_b\left[
 1-\frac{\det\!\left(\Gamma+\xi_a\xi_b^T/(c_ac_b)\right)}{\det\Gamma}
 \right]\right\},
\end{equation}
with no summation over $a,b$. The scalars cancel from the complete expression. The diagonal rescalings \eqref{eq:Gamma-M} then give \eqref{eq:Gdet}--\eqref{eq:Gdet-odd}, including the off-diagonal correction \eqref{eq:det-offdiag}. Thus the determinant numerator and denominator reconstruct the same metric, with the normalization $h_n$ and reality conditions fixed above. The relation to monodromy data is the gauge equivalence of the BZ and monodromy spectral problems described in \cite{Alekseev}.

In tau variables, the rank-one correction has a more specific meaning. The denominator is the product $\tau_+^{(\alpha)}\tau_-^{(\alpha)}/\Pi$. At coefficient $-(-1)^n/\Lambda$, the update $a^{(\alpha)}\one^T$ replaces it by $\tau_+^{(\alpha+1)}\tau_-^{(\alpha-1)}/\Pi$, as in \eqref{eq:polarized-Cauchy}; reversing the update gives the opposite pair of shifts. Their difference reconstructs $G_{t\phi}$. This changes the two tau weights in opposite directions while retaining the same poles and charges; it is not an additional dressing step.
\end{remark}

The symmetric and antisymmetric shifts recover respectively the
diagonal and off-diagonal metric data. The complete Gram matrices are
\begin{equation}\label{eq:Gtau}
 G=\frac1{\tau_+^{(1)}\tau_-^{(1)}}
 \begin{pmatrix}
 \Lambda\tau_+^{(0)}\tau_-^{(0)}&
 \dfrac{\ii\rho}{2}
 (\tau_+^{(2)}\tau_-^{(0)}-\tau_+^{(0)}\tau_-^{(2)})\\[3pt]
 \dfrac{\ii\rho}{2}
 (\tau_+^{(2)}\tau_-^{(0)}-\tau_+^{(0)}\tau_-^{(2)})&
 \dfrac{\rho^2}{\Lambda}\tau_+^{(2)}\tau_-^{(2)}
 \end{pmatrix},\quad n\ \mathrm{even},
\end{equation}
and
\begin{equation}\label{eq:Gtau-odd}
 G=\frac{\rho}{\tau_+^{(0)}\tau_-^{(0)}}
 \begin{pmatrix}
 -\ii\Lambda\tau_+^{(-1)}\tau_-^{(-1)}&
 \dfrac{\ii}{2}
 (\tau_+^{(1)}\tau_-^{(-1)}-\tau_+^{(-1)}\tau_-^{(1)})\\[3pt]
 \dfrac{\ii}{2}
 (\tau_+^{(1)}\tau_-^{(-1)}-\tau_+^{(-1)}\tau_-^{(1)})&
 \dfrac{\ii}{\Lambda}\tau_+^{(1)}\tau_-^{(1)}
 \end{pmatrix},\quad n\ \mathrm{odd}.
\end{equation}
In both cases, \eqref{eq:bilinear} gives $\det G=\rho^2$ directly.
Replacing $d_1$ or $d_0$ in \eqref{eq:nudet} by \eqref{eq:factor}
expresses the conformal factor through the same tau data. Thus the
signed tau functions, the pole data and $k_0$ determine the full metric
without a further quadrature.

\subsection{Ernst potentials from tau ratios}

\paragraph{Metric and Killing descriptions.}
Write
\begin{equation}\label{eq:Ernst-Gram}
 G=\begin{pmatrix}A&B\\B&D\end{pmatrix},\qquad AD-B^2=\rho^2.
\end{equation}
The metric Ernst pair is defined algebraically by
\begin{equation}\label{eq:Ernst-metric}
 \mathcal E_{M,\pm}=\frac{\rho\mp\ii B}{A}.
\end{equation}
For a real Riemannian metric the two components are complex conjugates.
Conversely, setting $E_\pm=\mathcal E_{M,\pm}$ reconstructs
\begin{equation}\label{eq:Ernst-reconstruct}
 A=\frac{2\rho}{E_++E_-},\qquad
 B=\frac{\ii\rho(E_+-E_-)}{E_++E_-},\qquad
 D=\frac{2\rho E_+E_-}{E_++E_-}.
\end{equation}

For a chosen Killing field $\xi$, define instead its norm and twist by
\begin{equation}\label{eq:Ernst-Killing}
 F_\xi=g(\xi,\xi),\qquad
 d\chi_\xi=*_{g}(\xi^\flat\wedge d\xi^\flat),\qquad
 \mathcal E_{\xi,\pm}=F_\xi\pm\chi_\xi.
\end{equation}
The twist one-form is closed in vacuum, so its primitive exists locally;
its additive constant must be specified. Fix the orientation
$dt\wedge d\phi\wedge d\rho\wedge dz$. For the two coordinate Killing
fields the defining equations are
\begin{align}
 \chi_{t,\rho}&=\frac{AB_z-BA_z}{\rho},&
 \chi_{t,z}&=-\frac{AB_\rho-BA_\rho}{\rho},
       \label{eq:twist-t}\\
 \chi_{\phi,\rho}&=-\frac{DB_z-BD_z}{\rho},&
 \chi_{\phi,z}&=\frac{DB_\rho-BD_\rho}{\rho}.
       \label{eq:twist-phi}
\end{align}
The base conformal factor cancels from these equations. Both components
of $\mathcal E_{\xi,\pm}$ are real on the Riemannian slice; they are
generally distinct and are not a complex-conjugate pair.

Both descriptions satisfy the paired vacuum Ernst equation
\begin{equation}\label{eq:Ernst-paired}
 \mathcal O_\pm(E_+,E_-)
 :=(E_++E_-)\Delta E_\pm
       -2\nabla E_\pm\cdot\nabla E_\pm=0,
\end{equation}
where
\begin{equation}
 \Delta=\partial_\rho^2+\rho^{-1}\partial_\rho+\partial_z^2,
 \qquad \nabla f\cdot\nabla h=f_\rho h_\rho+f_z h_z.
\end{equation}
The dot product is bilinear, without complex conjugation. For the metric
pair, \eqref{eq:Ernst-paired} becomes
$(\Re\mathcal E_M)\Delta\mathcal E_M
=\nabla\mathcal E_M\cdot\nabla\mathcal E_M$.
For the real Killing pair, it is equivalent to
\begin{equation}\label{eq:Ernst-real}
 F_\xi\Delta F_\xi
 =\nabla F_\xi\cdot\nabla F_\xi
       +\nabla\chi_\xi\cdot\nabla\chi_\xi,
 \qquad
 F_\xi\Delta\chi_\xi=2\nabla F_\xi\cdot\nabla\chi_\xi.
\end{equation}
These are different real forms of the same paired equation. The
conversion from the metric description to the Killing description is
the norm--twist B\"acklund transformation, expressed differentially by
\eqref{eq:twist-t}--\eqref{eq:twist-phi}; see also \cite{Tod}.

\paragraph{The twist from spectral infinity.}
The first coefficient of the dressing matrix already supplies the
required primitives. At fixed $(\rho,z)$, expand
\begin{equation}\label{eq:dressing-infinity}
 \mathscr X(\lambda_{\mathrm{sp}})=I+\frac{\mathcal W}{\lambda_{\mathrm{sp}}}+O(\lambda_{\mathrm{sp}}^{-2}),\qquad
 \mathcal W=\sum_iR_i
  =\sum_{i,j}(\Gamma^{-1})_{ij}\frac{N_j^Tm_i}{\mu_j}.
\end{equation}
This is spectral infinity, not an asymptotic limit in spacetime.
Let $(U_0,V_0)$ and $(U_{\mathrm{un}},V_{\mathrm{un}})$ denote the currents of the seed and unnormalized dressed metric, respectively. Substituting into the dressing equations gives
\begin{equation}\label{eq:W-currents}
 \mathcal W_\rho=V_{\mathrm{un}}-V_0,\qquad
 \mathcal W_z=-U_{\mathrm{un}}+U_0.
\end{equation}
Indeed, the two Lax coefficients start as $V/\lambda_{\mathrm{sp}}$ and $-U/\lambda_{\mathrm{sp}}$,
respectively, while their spectral-derivative terms do not contribute
at order $\lambda_{\mathrm{sp}}^{-1}$. The seed currents are diagonal in both parity
conventions, and scalar determinant normalization leaves off-diagonal
currents unchanged. Equations \eqref{eq:twist-t}--\eqref{eq:twist-phi}
therefore identify the dressing-normalized primitives as
\begin{equation}\label{eq:W-twists}
 \chi_t=\mathcal W_{12},\qquad \chi_\phi=-\mathcal W_{21}.
\end{equation}

There is no remaining integration in the evaluation of these entries.
For example, when $n$ is even and $S=(M^{(1)})^{-1}$,
\begin{equation}
 \mathcal W_{12}=\sum_{i,j}(\mu_iS_{ij}\mu_j)
             \frac{C_j}{\mu_j}\frac1{\mu_i}
         =\one^TSC.
\end{equation}
Using \eqref{eq:bz-seeds} and \eqref{eq:Gamma-M} for both parities gives
\begin{equation}\label{eq:twists-matrix}
\begin{array}{c|c|c}
 &\chi_t&\chi_\phi\\ \hline
 n\ \mathrm{even}
 &\one^T(M^{(1)})^{-1}C
 &-\rho^2(\lambda^{-1})^T(M^{(1)})^{-1}(C\lambda)\\[3pt]
 n\ \mathrm{odd}
 &\ii\rho\one^T(M^{(0)})^{-1}(C\lambda^{-1})
 &\ii\rho C^T(M^{(0)})^{-1}\lambda^{-1}
\end{array}
\end{equation}

\begin{proposition}[Three Ernst descriptions of one BZ metric]
\label{prop:Ernst-three}
In the Killing frame \eqref{eq:WP}, with the orientation and twist
constants fixed above, the three Ernst pairs of the same flat-seed
metric are
\begin{equation}\label{eq:Ernst-tau-table}
\renewcommand{\arraystretch}{1.8}
\begin{array}{c|c|c}
 &n\ \mathrm{even}&n\ \mathrm{odd}\\ \hline
 \mathcal E_{M,\pm}
 &\displaystyle\frac{\rho}{\Lambda}
       \frac{\tau_\pm^{(2)}}{\tau_\pm^{(0)}}
 &\displaystyle\frac{\ii}{\Lambda}
       \frac{\tau_\pm^{(1)}}{\tau_\pm^{(-1)}}\\
 \mathcal E_{\partial_t,\pm}
 &\displaystyle\Lambda
       \frac{\tau_\pm^{(-1)}}{\tau_\pm^{(1)}}
 &\displaystyle-\ii\rho\Lambda
       \frac{\tau_\pm^{(-2)}}{\tau_\pm^{(0)}}\\
 \mathcal E_{\partial_\phi,\pm}
 &\displaystyle\frac{\rho^2}{\Lambda}
       \frac{\tau_\pm^{(3)}}{\tau_\pm^{(1)}}
 &\displaystyle\frac{\ii\rho}{\Lambda}
       \frac{\tau_\pm^{(2)}}{\tau_\pm^{(0)}}
\end{array}
\end{equation}
Each row satisfies \eqref{eq:Ernst-paired}. The metric pair is complex
conjugate on the real slice; the two Killing pairs are real.
\end{proposition}
\begin{proof}
The first row follows from \eqref{eq:Ernst-metric},
\eqref{eq:Gtau}--\eqref{eq:Gtau-odd} and the symmetric shift.
For the other rows, the symmetric shift identifies the half-sums with
$A$ and $D$, while the antisymmetric shift identifies the half-differences
with \eqref{eq:twists-matrix}. The rank-one updates
\eqref{eq:Mupdate} transfer the inverse matrices to the required center.

For example, in the even case,
\begin{equation}
 \frac{d_0}{d_1}\one^T(M^{(0)})^{-1}C
       =\one^T(M^{(1)})^{-1}C.
\end{equation}
Using \eqref{eq:antisymmetric-shift} at $\alpha=0$ consequently gives
\begin{equation}
 \chi_t=\frac{\Lambda}{2}
 \left(\frac{\tau_+^{(-1)}}{\tau_+^{(1)}}
       -\frac{\tau_-^{(-1)}}{\tau_-^{(1)}}\right),
\end{equation}
while the corresponding half-sum is $\Lambda d_0/d_1=A$.
For the $\phi$-potential use the shift centered at $2$ and
$M^{(2)}=\one\one^T+D_\lambda M^{(1)}D_\lambda$.
In the odd case the same argument uses centers $-1$ and $1$,
respectively, with the parity sign in
\eqref{eq:antisymmetric-shift}. The identification with the geometric
twists follows from \eqref{eq:W-currents}--\eqref{eq:W-twists}, rather
than from the Ernst equation alone.
\end{proof}

\begin{remark}[Neugebauer--Kramer determinants]\label{rem:NK-dictionary}
The determinant construction of Neugebauer and Kramer
\cite{Neugebauer1980,NeugebauerKramer} gives the same Killing Ernst pair.
To make the comparison explicit, use the monomial convention of
\cite{ErnstNeugebauer} and let $n=2m$. Set
\begin{equation}\label{eq:NK-dictionary}
 K_i=z_i,\qquad R_i=\sqrt{(z-K_i)^2+\rho^2},\qquad
 X_i=\frac{C_i+1}{C_i-1}.
\end{equation}
Define $U_{\mathrm{NK}}(X)$ as the determinant of the $n\times n$ matrix
whose successive row pairs are
$(K_i^jX_iR_i)_{i=1}^n$ and $(K_i^j)_{i=1}^n$, for
$j=0,\ldots,m-1$. Define $W_{\mathrm{NK}}(X)$ by replacing row $n-1$
(numbered from one) by $(K_i^m)_{i=1}^n$. Then
\begin{equation}\label{eq:NK-tau}
 \frac{U_{\mathrm{NK}}(X)-W_{\mathrm{NK}}(X)}
      {U_{\mathrm{NK}}(X)+W_{\mathrm{NK}}(X)}
 =\Lambda\frac{\tau_+^{(-1)}}{\tau_+^{(1)}}
 =\mathcal E_{\partial_t,+}.
\end{equation}
Replacing $X_i$ by $X_i^{-1}$ gives the minus component, since this is
precisely $C_i\mapsto-C_i$. The comparison uses the complexified
determinant formulas with the Riemannian reality conditions of this
paper; the Lorentzian conditions on $X_i$ select a different real slice.

For the algebraic identity, clear column $i$ by multiplying by $C_i-1$
and shift the polynomial row basis to $K_i-z$. Writing
$K_i-z=(\mu_i-\rho^2/\mu_i)/2$ and
$R_i=(\mu_i+\rho^2/\mu_i)/2$, the column expansion selects one of the
two sheets $\mu_i$ and $-\rho^2/\mu_i$, equivalently $\lambda_i$ and
$\lambda_i^{-1}$. The resulting Vandermonde expansion is exactly the
signed subset sum \eqref{eq:tau}. More explicitly, with
$D_C=\prod_i(C_i-1)$, a common nonzero spectral factor $\mathfrak f$,
independent of the charges, satisfies
\begin{equation}\label{eq:NK-factors}
 \begin{split}
 D_C(U_{\mathrm{NK}}+W_{\mathrm{NK}})&=\mathfrak f\,\tau_+^{(1)},\\
 D_C(U_{\mathrm{NK}}-W_{\mathrm{NK}})&=\mathfrak f\,\Lambda\tau_+^{(-1)}.
 \end{split}
\end{equation}
Thus the monomial determinants already represent individual tau
functions, up to common factors. In the BZ determinant expressions the
product factorization \eqref{eq:factor} and its polarized version
\eqref{eq:polarized-Cauchy} cancel the common tau factor to give the
same ratio. No change of the turning-point positions is required.
The dictionary is written on the generic chart $C_i\ne1$; after clearing
denominators, the phantom values $C_i=\pm1$ correspond to the projective
limits $X_i=\infty,0$. This comparison concerns the Killing potential
of $\partial_t$, not the algebraic metric potential
$\mathcal E_{M,\pm}$.
\end{remark}

\paragraph{The shifted family.}
Introduce
\begin{equation}\label{eq:Falpha}
 \mathcal F_{\alpha,\pm}
 =\frac{\rho^\alpha}{\Lambda}
       \frac{\tau_\pm^{(\alpha+1)}}{\tau_\pm^{(\alpha-1)}},\qquad
 \widehat{\mathcal F}_{\alpha,\pm}=s_n\mathcal F_{\alpha,\pm},
 \qquad
 s_n=\begin{cases}1,&n\ \mathrm{even},\\ \ii,&n\ \mathrm{odd}.
 \end{cases}
\end{equation}
If $k_n=1$ for even $n$ and $k_n=0$ for odd $n$,
Proposition~\ref{prop:Ernst-three} becomes
\begin{equation}\label{eq:Ernst-neighbors-flat}
 \mathcal E_{M,\pm}=\widehat{\mathcal F}_{k_n,\pm},\qquad
 \mathcal E_{\partial_t,\pm}
       =\widehat{\mathcal F}_{k_n-1,\pm}^{-1},\qquad
 \mathcal E_{\partial_\phi,\pm}
       =\widehat{\mathcal F}_{k_n+1,\pm}.
\end{equation}
The neighboring-shift relation is not restricted to these two centers.

\begin{proposition}[Shifted Ernst family]
\label{prop:Ernst-shifted}
For either parity of $n$, the pair
$\widehat{\mathcal F}_{\alpha,\pm}$ and its reciprocal solve
\eqref{eq:Ernst-paired}. This holds for integer $\alpha$, and locally
for arbitrary $\alpha$ after consistent choices of the powers.
For the metric reconstructed from
$\mathcal E_{M,\pm}^{[\alpha]}=
\widehat{\mathcal F}_{\alpha,\pm}$, the dressing-normalized Killing
potentials satisfy
\begin{equation}\label{eq:Ernst-neighbors-general}
 \mathcal E_{\partial_t,\pm}^{[\alpha]}
       =\widehat{\mathcal F}_{\alpha-1,\pm}^{-1},\qquad
 \mathcal E_{\partial_\phi,\pm}^{[\alpha]}
       =\widehat{\mathcal F}_{\alpha+1,\pm}.
\end{equation}
This is a statement about the complexified local construction, before
imposing a real slice and global regularity.
\end{proposition}
\begin{proof}
For even $n$, dress the diagonal power-law seed
\begin{equation}\label{eq:shifted-seeds-even}
 G_0^{[\alpha]}=\diag(\rho^{1-\alpha},\rho^{1+\alpha}),\qquad
 \Psi_0^{[\alpha]}
 =\diag(Q^{(1-\alpha)/2},Q^{(1+\alpha)/2}).
\end{equation}
Its currents are $U_0=\diag(1-\alpha,1+\alpha)$ and $V_0=0$.
Absorbing constant factors at $w=z_i$ into the polarization rows gives
\begin{equation}
 m_i=(\eta_i\mu_i^{(\alpha-1)/2},
                    \mu_i^{-(\alpha+1)/2}),\qquad
 C_i=\ii^{\,1-\alpha}\eta_i.
\end{equation}
For odd $n$, use the generalized seed and rows
\begin{align}
 G_0^{[\alpha]}&=\diag(-\rho^{-\alpha},\rho^\alpha),&
 \Psi_0^{[\alpha]}&=\diag(-Q^{-\alpha/2},Q^{\alpha/2}),
        \label{eq:shifted-seeds-odd}\\
 m_i&=(-\eta_i\mu_i^{\alpha/2},\mu_i^{-\alpha/2}),&
 C_i&=\ii^{-\alpha}\eta_i.\notag
\end{align}
In both cases $\Gamma$ is diagonally congruent to $M^{(\alpha)}$.
After \eqref{eq:bz-normalization}, the determinant lemma yields
\begin{align}
 G_{tt}^{[\alpha]}
 &=s_n^{-1}\rho^{1-\alpha}\Lambda
                       \frac{d_{\alpha-1}}{d_\alpha},&
 G_{\phi\phi}^{[\alpha]}
 &=\frac{s_n\rho^{1+\alpha}}{\Lambda}
                       \frac{d_{\alpha+1}}{d_\alpha},
                       \label{eq:shifted-Gram}\\
 G_{t\phi}^{[\alpha]}
 &=\frac{\ii\rho}{2}
 \frac{\tau_+^{(\alpha+1)}\tau_-^{(\alpha-1)}
       -\tau_+^{(\alpha-1)}\tau_-^{(\alpha+1)}}
      {\tau_+^{(\alpha)}\tau_-^{(\alpha)}}.&&\notag
\end{align}
Thus its metric Ernst pair is
$\widehat{\mathcal F}_{\alpha,\pm}$, proving the field equation.
Evaluating \eqref{eq:dressing-infinity} for these rows gives
\begin{align}
 \chi_t^{[\alpha]}
 &=s_n\rho^{1-\alpha}
   \one^T(M^{(\alpha)})^{-1}(C\lambda^{\alpha-1}),\notag\\
 \chi_\phi^{[\alpha]}
 &=-s_n^{-1}\rho^{1+\alpha}
   (\lambda^{-1})^T(M^{(\alpha)})^{-1}(C\lambda^\alpha).
   \label{eq:shifted-twists}
\end{align}
The two shift identities now give
\eqref{eq:Ernst-neighbors-general}. Finally, simultaneous inversion of
the two components preserves \eqref{eq:Ernst-paired}.
\end{proof}

For real $C_i$ and integer $\alpha$, the reality condition alternates:
\begin{equation}\label{eq:Falpha-reality}
 \begin{cases}
 \widehat{\mathcal F}_{\alpha,-}
       =\overline{\widehat{\mathcal F}_{\alpha,+}},
       &n+\alpha\ \mathrm{odd},\\
 \widehat{\mathcal F}_{\alpha,+},\widehat{\mathcal F}_{\alpha,-}
       \in\R,&n+\alpha\ \mathrm{even}.
 \end{cases}
\end{equation}
The metric center and its neighbors therefore have precisely the
different reality structures required in
\eqref{eq:Ernst-neighbors-flat}. For noninteger shifts a real form
requires a corresponding choice of charge phases. Changing $\alpha$
and interpreting the result as a metric potential generally changes
the seed and the geometry; it need not preserve asymptotic flatness,
positivity or rod regularity. In particular, a shift of the tau indices
must be accompanied by the power $\rho^\alpha$ in \eqref{eq:Falpha}.

\paragraph{Choice of Killing field.}
The oriented exchange
$(\partial_t,\partial_\phi)\mapsto(\partial_\phi,-\partial_t)$
sends the metric Ernst pair to its reciprocal, since
\begin{equation}
 \frac{\rho\pm\ii B}{D}
       =\frac{1}{\mathcal E_{M,\pm}}.
\end{equation}
The two Killing Ernst pairs of a fixed metric, however, are the distinct
neighboring expressions in \eqref{eq:Ernst-neighbors-flat}, not generally
reciprocals of one another.
For $\xi=a\partial_t+b\partial_\phi$ with constant $a,b$,
\begin{equation}\label{eq:Ernst-linear-Killing}
 F_\xi=a^2A+2abB+b^2D,\qquad
 \chi_\xi=a^2\chi_t+ab\chi_{\mathrm{mix}}+b^2\chi_\phi.
\end{equation}
The same dressing coefficient supplies the mixed primitive:
\begin{equation}\label{eq:Ernst-mixed}
 \chi_{\mathrm{mix}}
 =\mathcal W_{22}-\mathcal W_{11}+
 \begin{cases}-2z,&n\ \mathrm{even},\\0,&n\ \mathrm{odd}.
 \end{cases}
\end{equation}
The seed term follows from $U_{0,11}-U_{0,22}=-2$ for the even seed
and $0$ for the odd seed. Additive constants may subsequently be chosen
to impose an asymptotic twist normalization. Rescaling a Killing field
by a constant $c$ multiplies both its norm and its twist by $c^2$.

\section{Geometry of the tau-function}\label{sec:geometry}

The subset sum \eqref{eq:tau} has two complementary uses. Its Pfaffian form explains the bilinear identities and makes insertion of additional solitons transparent. Its subset form makes degenerations and axis limits easy to read. We retain both descriptions, using the same variables $a,b$ as in Section~\ref{sec:three}.

\subsection{Pfaffian identities}
Define the skew matrix
\begin{equation}\label{eq:pf-matrix}
 \mathsf H_{ij}=(1-a_i a_j)K(b_i,b_j)+a_j-a_i.
\end{equation}
Then
\begin{equation}\label{eq:pfaffian}
 \tau_n(a,b)=
 \begin{cases}
 \Pf\mathsf H,&n\ \text{even},\\[2pt]
 \displaystyle\Pf\begin{pmatrix}0&(\one-a)^T\\a-\one&\mathsf H\end{pmatrix},
       &n\ \text{odd}.
 \end{cases}
\end{equation}
The bordering convention fixes the sign for odd $n$; for $n=1$ it gives $1-a_1$. Expansion in the $a_i$ and Schur's Cauchy--Pfaffian identity recover \eqref{eq:tau} \cite{CauchyPfaffian}.

Several determinant squares now have a common explanation. For even $n$, $\det\mathsf H=\tau_n^2$; for odd $n$, the same statement holds for the bordered matrix. This is different from the BZ determinant factorization $d_\alpha=\tau_+^{(\alpha)}\tau_-^{(\alpha)}/\Pi$, which pairs the two charge signs. In particular the determinant identity underlying the off-diagonal metric entry is
\begin{equation}\label{eq:det-square-recalled}
 d_\alpha^2-d_{\alpha+1}d_{\alpha-1}
 =\left(\frac{\one^T\adj(M^{(\alpha)})a^{(\alpha)}}{\Lambda}\right)^2.
\end{equation}
The square does not determine a sign. The signed identity \eqref{eq:antisymmetric-shift} fixes that sign before the metric is reconstructed.

The spectral shifts are small operations on the skew matrix. If $a=C b^\alpha$ and $x_\alpha=C b^\alpha-\one$, then
\begin{equation}\label{eq:pf-rank-two}
 \mathsf H_{\alpha+1}-\mathsf H_\alpha
 =x_{\alpha+1}x_\alpha^T-x_\alpha x_{\alpha+1}^T.
\end{equation}
Likewise charge reversal changes $\mathsf H$ by rank at most two. On an invertible even chart,
\begin{equation}\label{eq:pf-update}
 \Pf(\mathsf H+xy^T-yx^T)
 =\Pf(\mathsf H)\bigl(1+y^T\mathsf H^{-1}x\bigr).
\end{equation}
After clearing denominators this extends across singular charts. It explains why neighbouring tau-functions, rather than unrelated determinants, suffice for the metric and Ernst fields.

\subsection{The insertion recurrence}
Adjoin a slot with charge weight $A$ and spectral coordinate $\zeta$. Put $k_i=K(b_i,\zeta)$ and $L_k=\prod_i k_i$. Splitting the subsets according to whether the new index is occupied gives
\begin{equation}\label{eq:recurrence-one}
 \tau_{n+1}(a,A;b,\zeta)
 =L_k\tau_n(a/k,b)+(-1)^{n+1}A\tau_n(ak,b).
\end{equation}
Products and quotients of weight vectors are componentwise. An empty new slot couples to every empty old slot; an occupied new slot couples to every occupied old slot. The two terms in \eqref{eq:recurrence-one} express precisely this alternative.

For a second new slot $(B,\eta)$, set $h_i=K(b_i,\eta)$, $L_h=\prod_i h_i$ and $c=K(\zeta,\eta)$. Applying the recurrence twice gives
\begin{align}\label{eq:recurrence-two}
 \tau_{n+2}(a,A,B;b,\zeta,\eta)
 ={}&cL_kL_h\tau_n(a/(kh),b)\notag\\
 &+(-1)^{n+1}A L_h\tau_n(ak/h,b)\notag\\
 &+(-1)^{n+2}B L_k\tau_n(ah/k,b)
       -ABc\tau_n(akh,b).
\end{align}
The four terms correspond to the four occupations of the new pair. For $n=0$ this reduces to $(1-AB)c+B-A$, providing a useful sign check. Both formulas are polynomial identities after denominators are cleared. A zero kernel or an infinite affine charge must therefore be treated in the complete expression, or in homogeneous charge coordinates, rather than term by term.

The same idea appears in block elimination. For an invertible even skew block $B_0$,
\begin{equation}\label{eq:pf-schur}
 \Pf\begin{pmatrix}A_0&E\\-E^T&B_0\end{pmatrix}
 =\Pf(B_0)\Pf(A_0+EB_0^{-1}E^T).
\end{equation}
A distant group of solitons can consequently change the effective local skew matrix even when the group itself has left the observation region. This is the algebraic reason to retain the cross terms in a cluster limit. A smaller number of visible poles does not, by itself, give the balance equations of an isolated lower-order solution.

\subsection{The pure-spinor viewpoint}
There is a geometric structure behind these formulas. Let $E$ be a vector space with dual basis $e^1,\ldots,e^n$. The natural pairing on $E\oplus E^*$ is a split quadratic form. A skew form $\mathsf H$ determines an isotropic graph in this space. In the standard exterior-algebra realization of the spin representation, that graph is represented by the pure spinor
\begin{equation}\label{eq:pure-spinor}
 s_{\mathsf H}=\exp\!\left(\frac12\sum_{i,j}\mathsf H_{ij}e^i\wedge e^j\right).
\end{equation}
Its coordinates are the principal Pfaffians of $\mathsf H$, including the empty Pfaffian $1$. The even tau-function is the top coordinate; the bordered odd formula places the odd case in the same even-dimensional construction.

Thus the tau-function is one coordinate of a structured object, rather than an isolated polynomial. The principal Pfaffians obey quadratic relations because the corresponding isotropic planes lie on an orthogonal Grassmannian. The bilinear identities encountered in the metric are concrete specializations of these relations. Rank-two updates change the skew form by a decomposable two-form, while a change of Pfaffian chart changes the coordinate description of the same geometric construction. This viewpoint explains why the algebra closes on a small collection of neighbouring quantities.

Related pure-spinor and Pfaffian coordinates occur in the $B$-Toda and BKP hierarchies \cite{Xie}. No Grassmannian classification or hierarchy of additional times is required below. In particular, identifying this finite Pfaffian geometry does not by itself specify a BKP evolution or a new field equation. Its role here is organizational: it explains the recurrences, the determinant squares and the compatibility of the descendants already used in the BZ metric. The two-component metric and boundary spinors introduced below are descendants of this data; they are not the full Cartan pure spinor \eqref{eq:pure-spinor}.

\subsection{Surviving sectors in the rod limit}
For the remainder of this section take $N=n=2m$ genuine ordered nodes. The axis is divided into rods $I_j=(z_j,z_{j+1})$, with $I_0=(-\infty,z_1)$ and $I_N=(z_N,\infty)$. For the branch \eqref{eq:mu}, the coordinates $\lambda_i$ are of order $\rho$ for $i\le j$ and of order $\rho^{-1}$ for $i>j$. Inverting the latter transfers the limit to a small-coordinate chart, with the corresponding charges inverted. Complementing the occupation subsets gives the equivalent prefix-inversion convention used in \eqref{eq:tau-rods}.

In that chart the term with $|S|=k$ in the unshifted tau-function has order
\begin{equation}\label{eq:sector-order}
 \rho^{\binom{k}{2}+\binom{2m-k}{2}}
 =\rho^{m(m-1)+(k-m)^2}.
\end{equation}
The middle sector $k=m$ survives at leading order. Multiplication of the occupied weights by $b_i^\alpha$ adds $\alpha k$ to the exponent. The $\alpha=2$ descendant therefore selects $k=m-1$; complementary inversion expresses the same information through $k=m+1$. It is the entire sector, not one preferred subset, that survives.

Define the resulting Vandermonde sectors by
\begin{equation}\label{eq:axis-sectors}
 W_k(C,z)=\sum_{|S|=k}(-1)^{\sum_{i\in S}i}C_S
              \Delta_z(S)\Delta_z(S^c),\qquad
 \Delta_z(S)=\prod_{\substack{i<j\\i,j\in S}}(z_i-z_j).
\end{equation}
Here $C_S=\prod_{i\in S}C_i$. They also have an ordinary determinant form. If the $i$th row of $V_k$ is
\[
 \bigl(C_i,C_i z_i,\ldots,C_i z_i^{k-1},
             1,z_i,\ldots,z_i^{N-k-1}\bigr),
\]
then Laplace expansion gives
\begin{equation}\label{eq:axis-vandermonde-det}
 W_k=(-1)^{k^2+(N-k)(N-k-1)/2}\det V_k.
\end{equation}
This form avoids an exponential subset expansion when evaluating the axis data.

A cone-normalized rod vector $\mathsf r_j$ spans the kernel of $G$ on $I_j$ and is scaled so that
\[
 \lim_{\rho\to0}\frac{G(\mathsf r_j,\mathsf r_j)}{\rho^2e^{2\nu}}=1.
\]
Let $C^{[j]}$ be obtained by inverting the first $j$ charges. With the common normalization $k_0$ of \eqref{eq:ffBZ}, a consistent choice is
\begin{equation}\label{eq:tau-rods}
 \mathsf r_j=\sqrt{k_0}\,\gamma_m\left(\prod_{i\le j}C_i\right)
       \begin{pmatrix}-2W_{m+1}(C^{[j]},z)\\ W_m(C^{[j]},z)\end{pmatrix},
 \qquad \gamma_m=(-1)^{m(m+1)/2}2^{m(m-1)}.
\end{equation}
An individual rod may be reversed in sign. Its magnitude cannot be discarded: it is fixed by the cone normalization. In homogeneous coordinates $C_i=P_i/Q_i$, replace $C_S$ by $P_SQ_{S^c}$; the operation on the first $j$ indices is simply $P_i\leftrightarrow Q_i$. The formula then extends through $C_i=0$ and $C_i=\infty$ without changing charts in the final rod equations.

Finally, the metric itself has a two-column factorization. In the even convention set
\begin{equation}\label{eq:metric-spinors}
 \psi_\pm=\frac1{\tau_\pm^{(1)}}
 \begin{pmatrix}\sqrt\Lambda\,\tau_\pm^{(0)}\\
 \pm\ii\rho\Lambda^{-1/2}\tau_\pm^{(2)}\end{pmatrix}.
 \qquad
 G=\frac12(\psi_+\psi_-^T+\psi_-\psi_+^T).
\end{equation}
This is a factorization of the bilinear formula of Section~\ref{sec:three}. Taking its rod limit together with the conformal factor gives \eqref{eq:tau-rods}. Keeping the kinematic prefactors and the common conformal normalization is essential: the projective null direction alone contains no cone angle.

\section{Spinor factorization of the patching matrix}\label{sec:patching}

The bulk tau-function determines the metric, while its axis restriction has a smaller rational representative. This section develops that boundary representation for genuine even data $n=N=2m$, with arbitrary ordered nodes and charges. The exterior Killing frame has unit translational norm at infinity and zero asymptotic twist potential.

\subsection{Positive boundary data}
Put
\[
 D(w)=\prod_{i=1}^N(w-z_i),\qquad Z=\diag(z_1,\ldots,z_N).
\]
The positive real sector is specified by the weights
\begin{equation}\label{eq:positive-weights}
 \varpi_i=\frac{C_i-1}{(C_i+1)D'(z_i)}>0,\qquad
 [P_i:Q_i]=[1+\varpi_iD'(z_i):1-\varpi_iD'(z_i)].
\end{equation}
Here $[P_i:Q_i]$ is the homogeneous polarization coordinate, so $C_i=P_i/Q_i$ when $Q_i\ne0$. In particular $C_i=\infty$ can be an interior chart value. The geometric meaning of the chamber and its walls is discussed in Section~\ref{sec:regularity}.

Let $F=g(\partial_t,\partial_t)$ and $\chi=\chi_t$ be the Killing norm and twist of Section~\ref{sec:three}, with $\chi$ zero at infinity. The exterior-axis norm--twist matrix, analytically continued in the spectral coordinate $w$, is the \emph{patching matrix}
\[
 P(w)=\left.\frac1F
 \begin{pmatrix}1&-\chi\\-\chi&\chi^2-F^2\end{pmatrix}\right|_{\mathrm{exterior\ axis}}.
\]
In the positive BZ sector it has the representation
\begin{equation}\label{eq:patching-residues}
 P(w)=J_0+\sum_{i=1}^N\frac{s_i s_i^T}{w-z_i},\qquad
 J_0=\diag(1,-1),\qquad \det P=-1.
\end{equation}
The real columns $s_i=(u_i,v_i)^T$ are residue spinors, defined up to sign; set $u=(u_1,\ldots,u_N)^T$ and $v=(v_1,\ldots,v_N)^T$. These scalar components are distinct from the primitive lattice rods $\mathsf v_j$ introduced below.

The representation follows by restricting the BZ Ernst fields to an exterior axis. The converse boundary construction starts from rod-adapted solutions of the axis Lax system: matching at each nut gives a rational transition factor, and the ordered product is constrained by symmetry and its value at infinity \cite{LT,Tod,DMT}. This requires the local rod and nut models and the normalized end. Rationality alone is not a global regularity criterion. The explicit Chen--Teo ALF family \cite{ChenTeoALF} and its twistor factorization \cite{DunajskiTod} illustrate the distinction between algebraic reconstruction and the regular subfamily.

Positivity in \eqref{eq:patching-residues} concerns the residues; $P$ itself has determinant $-1$. The divided-difference kernel is
\begin{equation}\label{eq:positive-kernel}
 \frac{P(\eta)-P(w)}{w-\eta}
 =\sum_i\frac{s_i s_i^T}{(w-z_i)(\eta-z_i)}.
\end{equation}
For real sample points away from the poles its block matrix is a positive semidefinite Gram matrix. In particular $-P'(w)$ is positive semidefinite on each pole-free interval.

We use the mass $M$ and angular momentum $J$ normalized by
\[
 G_{tt}=1-\frac{2M}{r}+O(r^{-2}),\qquad
 G_{t\phi}=-\frac{2J\sin^2\vartheta}{r}+O(r^{-2}),
 \quad \rho=r\sin\vartheta,\quad z=r\cos\vartheta.
\]
In the AF frame, expansion of \eqref{eq:patching-residues} gives
\begin{equation}\label{eq:residue-moments}
 \|u\|^2=\|v\|^2=2M,\qquad u^Tv=0,\qquad
 J=-\frac12u^TZv.
\end{equation}
The vanishing mixed zeroth moment is the AF alignment condition. Mass and angular momentum are moments of the same data that determine the normalized rods.

\subsection{Polynomial spinors}
On a nonsingular interpolation chart, let $e$ be monic of degree $m$ and let $\deg o\le m-1$, determined by
\begin{equation}\label{eq:spinor-interpolation}
 e(z_i)+C_i o(z_i)=0,\qquad i=1,\ldots,2m.
\end{equation}
In homogeneous coordinates the equation is $Q_i e(z_i)+P_i o(z_i)=0$.
The letters $e,o$ label the two occupation sectors; they do not specify polynomial parity. Define
\[
 f_+=e+o,\qquad f_-=e-o.
\]
The general boundary representation retains two independent complementary polynomials $g_+,g_-$, with compatible normalization and degree $m$ on the generic genuine chart, satisfying
\begin{equation}\label{eq:spinor-bezout}
 f_+g_++f_-g_-=2D.
\end{equation}
Set
\begin{equation}\label{eq:spinor-AQ}
 \mathcal A=f_+f_-,\qquad
 \mathcal Q=D-f_+g_+,\qquad
 \mathcal C=-g_+g_-.
\end{equation}
Then
\begin{equation}\label{eq:spinor-patching}
 P(w)=\frac1{D(w)}
 \begin{pmatrix}\mathcal A(w)&\mathcal Q(w)\\
 \mathcal Q(w)&\mathcal C(w)\end{pmatrix},
 \qquad \mathcal A\mathcal C-\mathcal Q^2=-D^2.
\end{equation}
Indeed, writing $X=f_+g_+$ and $Y=f_-g_-$, the determinant numerator is
$-XY-(D-X)^2=-D^2$, because $X+Y=2D$. The complementary factorizations are
\begin{equation}\label{eq:spinor-complements}
 D-\mathcal Q=f_+g_+,\qquad D+\mathcal Q=f_-g_-.
\end{equation}

These relations are useful in both directions. Given the charges, \eqref{eq:spinor-interpolation} is a linear system for the coefficients of $e,o$. Given compatible boundary factors, the charges are recovered projectively as $[P_i:Q_i]=[-e(z_i):o(z_i)]$. The factor degrees, pole polarizations and real positive branch must be retained: the polynomial identity \eqref{eq:spinor-bezout} alone does not select a metric. A simultaneous zero of $e,o$ at a node also requires the complementary factors or the full residue; it is not permission to cancel a genuine pole.

The determinant lemma gives a spectral interpretation:
\begin{equation}\label{eq:chiral-spectrum}
 \mathcal A(w)=\det(wI-Z+uu^T),\qquad
 -\mathcal C(w)=\det(wI-Z-vv^T).
\end{equation}
Thus the combined roots of $f_+,f_-$ are the eigenvalues of the real symmetric matrix $Z-uu^T$, and those of $g_+,g_-$ are the eigenvalues of $Z+vv^T$. Each real chiral factor is real-rooted. Rank-one interlacing locates the roots relative to the nodes, with equalities possible at vanishing residue components. No corresponding real-rootedness assertion is needed for $e$ or $o$ separately.

A phantom is treated at the level of the complete matrix. A factor $w-z_p$ can be removed from $D$ only after every numerator entry has been shown to vanish at $z_p$. Degree and moment formulas must then use the reduced divisor. A disappearing polarization chart or a common factor of only one polynomial pair does not suffice.

\subsection{The positive real matrix reconstruction}
Let $V_{ik}=z_i^{k-1}$ for $1\le k\le m$, and let $e_m$ be the last coordinate vector in $\R^m$. Define
\begin{align}\label{eq:boundary-dictionary}
 W_+&=\diag(\sqrt{\varpi_i})V,&
 W_-&=\diag\!\left(\frac1{\sqrt{\varpi_i}D'(z_i)}\right)V,\notag\\
 a&=W_+(W_+^TW_+)^{-1}e_m,&
 c&=W_-(W_-^TW_-)^{-1}e_m,\notag\\
 T&=2W_+(W_+^TW_+)^{-1}W_+^T-I,&
 u&=a+c,\quad v=a-c.
\end{align}
Both Gram matrices are positive definite. Lagrange interpolation makes the two weighted polynomial spaces orthogonal complements and gives
\begin{equation}\label{eq:T-identities}
 T=T^T,\quad T^2=I,\quad Tv=u,\quad
 [Z,T]=uv^T-vu^T,\quad [Z+vv^T,T]=0,
\end{equation}
where $[X,Y]=XY-YX$. These are consequences of the interpolation, not extra field equations.

Let $T_j$ be the leading $j\times j$ block and $v_{[j]}=(v_1,\ldots,v_j)^T$. After exterior normalization, the cone-normalized rod formula becomes
\begin{equation}\label{eq:T-rods}
 \mathsf r_j=(-1)^j
 \left(2v_{[j]}^T\adj(T_j)v_{[j]},\ \det T_j\right),
 \qquad \mathsf r_0=(0,1).
\end{equation}
This adjugate expression continues through $\det T_j=0$, as is necessary when a rod has zero angular component.

For bulk reconstruction put
\[
 R_i=\sqrt{\rho^2+(z-z_i)^2},\qquad
 M_\pm=\diag(R_i\pm(z-z_i)),\qquad
 \mathcal B=M_-+TM_+T.
\]
For $\rho>0$, both $M_\pm$ are strictly positive, hence $\mathcal B>0$. Write the metric as
\[
 g=f(dt-\omega\,d\phi)^2+\rho^2f^{-1}d\phi^2
                  +e^{2k}f^{-1}(d\rho^2+dz^2),
 \qquad e^{2\nu}=e^{2k}/f.
\]
The BZ reconstruction is
\begin{equation}\label{eq:positive-bulk}
 f^{-1}=1+2v^T\mathcal B^{-1}v,\qquad
 \omega=-2v^T\mathcal B^{-1}M_-u,\qquad
 e^{2k}=\frac{\det\mathcal B}{2^N\prod_i R_i}.
\end{equation}
It proves positivity and absence of off-axis poles in this chamber. Axis completion still requires the conditions in the next section.

For axis limits it is preferable to reconstruct the Gram matrix directly.
Put $\mathcal R=\diag(R_i)$ and
\begin{equation}\label{eq:augmented-bulk}
 \widehat{\mathcal B}=\mathcal B+2vv^T
   =\mathcal R+T\mathcal RT+uu^T+vv^T,\qquad w=M_-u.
\end{equation}
Then the same metric is given by
\begin{equation}\label{eq:regular-gram}
 \begin{aligned}
 G_{tt}&=1-2v^T\widehat{\mathcal B}^{-1}v,&
 G_{t\phi}&=2v^T\widehat{\mathcal B}^{-1}w,\\
 G_{\phi\phi}&=\rho^2+2u^TM_-u-2w^T\widehat{\mathcal B}^{-1}w,&
 e^{2\nu}&=\frac{\det\widehat{\mathcal B}}{2^N\prod_iR_i}.
 \end{aligned}
\end{equation}
The first two entries follow from the rank-one inverse formula. For the third,
use $M_-M_+=\rho^2I$ and $T^2=I$ to obtain
\[
 M_--M_-\mathcal B^{-1}M_-=\rho^2T\mathcal B^{-1}T,
\]
and then apply the same rank-one formula. This proves
\eqref{eq:regular-gram} directly from \eqref{eq:positive-bulk}.
Unlike the separate norm and rotation coordinates, these expressions do not
divide by a leading principal minor. Their positivity on the closed axis is
used in Proposition~\ref{prop:canonical-comparison}.

\paragraph{Role in computation.}
Polynomial interpolation and factorization recover charges and asymptotic moments without large subset expansions. The real matrix form supplies stable bulk evaluation and the variation identities used in Section~\ref{sec:analytic}. Both retain the cone normalization. The original homogeneous subset rods remain an independent check of the numerical balance equations.

\subsection{Positive reduced moment ranks}\label{subsec:reduced-reconstruction}
The matrix identities admit a formulation that does not require equal positive and negative ranks. Let $1\le k\le N/2$, choose ordered nodes and positive weights $\omega_i$, and define
\[
 \mathcal W_{ia}=\sqrt{\omega_i}z_i^{a-1}\quad(1\le a\le k),\qquad
 \mathcal P=\mathcal W(\mathcal W^T\mathcal W)^{-1}\mathcal W^T.
\]
With $e_k$ the last coordinate vector in $\R^k$, put
\begin{equation}\label{eq:reduced-dictionary}
 \begin{aligned}
 a&=\mathcal W(\mathcal W^T\mathcal W)^{-1}e_k,&
 c&=(I-\mathcal P)(\sqrt{\omega_i}z_i^k)_{i=1}^N,\\
 T&=2\mathcal P-I,& u&=a+c,\quad v=a-c.
 \end{aligned}
\end{equation}
Polynomial orthogonality gives all identities in \eqref{eq:T-identities}. Common weight scaling sends $a$ to $s^{-1/2}a$ and $c$ to $s^{1/2}c$, so AF alignment $u^Tv=0$ fixes it uniquely. The determinant lemma for prefix Gram updates, followed by the bordered Schur complement, gives the same normalized rod formula \eqref{eq:T-rods}. Thus the positive bulk reconstruction and residue moments retain their form. In these coordinates the reduced slot number is $2(N-k)$ and the same-sign phantom count is $N-2k$; orientation reversal exchanges the two ranks.

The essential positivity property is irreducibility of the weighted polynomial space. If a diagonal matrix preserves it, its diagonal entries are the values of a polynomial of degree below $k$. Preservation of the highest monomial forces that polynomial to be constant, since $2k-2<N$. This observation supplies the equality argument in the reduced-rank cone Hessian and the properness of the direction map. It does not impose a limiting separation from phantom walls in an unreduced charge chart.

\subsection{Comparison with the positive Hankel family}\label{subsec:ADN-dictionary}
The closed family of Andersson--Dahl--Nilsson \cite[Sections~7--8]{ADN} uses the same positive discrete moment data. The following dictionary also fixes the cone normalization. For signed weights $\eta_i$, define
\[
 h_a(\eta)=\sum_i\eta_i z_i^a,\qquad
 \mathsf H_r(\eta)=\det(h_{a+b}(\eta))_{a,b=0}^{r-1},\qquad \mathsf H_0=1.
\]
Let $H=(h_{a+b})_{a,b=0}^{k-1}$, $c=(h_k,\ldots,h_{2k-1})^T$, and $e=e_k$. Put $\Sigma_k=c^T\adj(H)e$ and $s_z=\sum_i z_i$. The bordered determinant identity gives
\begin{equation}\label{eq:ADN-rod-dictionary}
 \det\begin{pmatrix}H&c+e\\(c+e)^T&h_{2k}+s_z\end{pmatrix}
 =\mathsf H_{k+1}-2\Sigma_k-\mathsf H_{k-1}+s_z\mathsf H_k.
\end{equation}
It follows first where $H$ is invertible and then everywhere by polynomial continuation. The last cofactor is $\mathsf H_{k-1}$, so no division by a signed minor is required.

For the signed axis weights $\eta_i=(-1)^{\mathbf1_{i\le j}}\omega_i$, the moment rod used in Appendix~\ref{app:numerics} is therefore
\[
 \mathsf t_j=\bigl(\mathsf H_{k+1}-2\Sigma_k-\mathsf H_{k-1}
                         +s_z\mathsf H_k,\ \mathsf H_k\bigr).
\]
In the notation of \cite[Equation~(5.22)]{ADN}, the shifted determinant is
$\mathsf H_{k+1}^{(a_\infty)}=\mathsf H_{k+1}-2a_\infty\Sigma_k-a_\infty^2\mathsf H_{k-1}$.
Thus our rods are its family vectors $(\mathsf H_k,\mathsf H_{k+1}^{(1)})$ after a common shear and component interchange. Here $a_\infty=1$ is the positive infinity-atom normalization; it is not an asymptotic twist. The flat-seed identification in \cite[Theorem~8.5]{ADN} identifies the corresponding metrics, including the reduced phantom sectors. A common positive factor is fixed by the conformal normalization.

This dictionary permits direct use of the analytic axis and end results of \cite[Propositions~7.5, 8.7--8.9 and 11.4]{ADN}. Positive separated data whose complete normalized rods form one primitive unimodular lattice have smooth complete AF or ALF completion, with the end type determined by the exterior rods. Accordingly, an exact root of the full vector equations in the appendices proves regularity; a direction-only root does not. This completion result does not establish a root's existence or the collision comparison needed for an entire AF parameter interval.

\section{Regularity of the rod structure}\label{sec:regularity}

Rod regularity means smooth collapse of primitive circles in one period lattice, including at their intersections \cite{ChenTeoRods,Nilsson}. The positive-family completion result in Section~\ref{subsec:ADN-dictionary} makes these normalized rod conditions sufficient for the AF and ALF ends considered here.

Fix the marked period lattice
\begin{equation}\label{eq:lattice}
 \Lambda_{\mathrm{per}}=\mathbb Z A+\mathbb Z B,\qquad
 A=(0,1),\quad B=(L,\beta),\quad L>0,
\end{equation}
where a lattice vector generates a flow of period $2\pi$ in the Killing coordinates $(t,\phi)$. We use the positively oriented basis $(A,-B)$, of covolume $L$, and write
\[
 \mathsf v_j=p_jA-q_jB=(-q_jL,p_j-q_j\beta).
\]
A vector is primitive when $\gcd(p_j,q_j)=1$. Signs do not change the vanishing circle. Choose signed lifts with
$\det_{\mathbb Z}(\mathsf v_{j-1},\mathsf v_j)=1$, starting from
$\mathsf v_0=A$ and $\mathsf v_1=-B$.
The subscript $\mathbb Z$ denotes the determinant of coefficient pairs in $(A,-B)$; the physical determinant is $L$ times this integer.

\subsection{The balance equations}
A regular gravitational instanton in this paper is a smooth, complete, oriented, positive definite Ricci-flat four-manifold with
$\int_M|\mathrm{Rm}|^2\,d\mathrm{vol}_g<\infty$, where $\mathrm{Rm}$ is the Riemann curvature tensor. In the finite-rod AF setting, regularity requires:
\begin{enumerate}
\item a positive nonsingular metric for $\rho>0$, distinct genuine nodes and positive finite Weyl intervals;
\item primitive vanishing vectors in one period lattice, with transverse angle $2\pi$ on every rod;
\item adjacent primitive vectors forming a lattice basis at every nut;
\item the smooth rod and nut normal forms and the prescribed complete AF end with finite curvature energy.
\end{enumerate}
The positive reconstruction \eqref{eq:positive-bulk} gives the bulk requirement. The finite lattice and cone conditions, together with the local BZ expansions established in Proposition~\ref{prop:canonical-comparison}, give smooth axis completion by Corollary~\ref{cor:smooth-completion}. The end estimates supply completeness and curvature decay \cite{KL,LiSun}.

For a primitive rod generator define its cone coefficient and angle by
\begin{equation}\label{eq:cone-definition}
 \kappa_j^2=\lim_{\rho\to0}
 \frac{G(\mathsf v_j,\mathsf v_j)}{\rho^2e^{2\nu}},
 \qquad \theta_j=2\pi\kappa_j,\qquad
 b_j=\log\frac{2\pi}{\theta_j}.
\end{equation}
The cone-normalized vector is
\[
 \mathsf r_j=\kappa_j^{-1}\mathsf v_j=e^{b_j}\mathsf v_j.
\]
It is the vector obtained in \eqref{eq:tau-rods} and \eqref{eq:T-rods}.
Regularity therefore requires the full vector equations $\mathsf r_j=\mathsf v_j$, with consistent signs, rather than only equality of their directions.

The signed rod sequence is determined by its self-intersection word
$W=(d_1,\ldots,d_{N-1})$:
\begin{equation}\label{eq:rod-recurrence}
 \mathsf v_{j+1}=-\mathsf v_{j-1}-d_j\mathsf v_j,\qquad
 (p_0,q_0)=(1,0),\quad(p_1,q_1)=(0,1).
\end{equation}
Here $d_j$ is the self-intersection of the sphere over the $j$th finite rod. At a smooth balanced metric,
\begin{equation}\label{eq:rod-area}
 \ell_j=z_{j+1}-z_j,\qquad
 \operatorname{Area}(S_j)=2\pi L\ell_j,\qquad
 x_j=\ell_j/L.
\end{equation}
Thus $x_j$ is the sphere area in units $2\pi L^2$. It is a geometric length coordinate, not a proper distance along the axis, and topology alone does not bound it away from zero.

\subsection{The positive charge chamber}
For $N=2m$ ordered nodes, $\operatorname{sign}D'(z_i)=(-1)^i$. Hence positivity of the weights \eqref{eq:positive-weights} is equivalent to
\begin{equation}\label{eq:alternating-chamber}
 |C_i|<1\quad(i\ \text{odd}),\qquad
 |C_i|>1\quad(i\ \text{even}).
\end{equation}
The adjacent-rod determinant has a square factor. For example, where $T_{i-1}$ is invertible, write
\[
 T_i=\begin{pmatrix}T_{i-1}&t_i\\t_i^T&T_{ii}\end{pmatrix},
 \qquad
 \eta_i=v_i-t_i^TT_{i-1}^{-1}v_{[i-1]}.
\]
The rod formula gives
\begin{equation}\label{eq:nut-square}
 \det(\mathsf r_{i-1},\mathsf r_i)
   =2(\det T_{i-1})^2\eta_i^2>0.
\end{equation}
For $i=1$, take $\det T_0=1$ and $\eta_1=v_1$.
Clearing the inverse gives the deleted-minor identity, which extends through singular principal blocks. The strict sign is supplied by the positive occupation sums in the genuine chamber.

The charge values $C_i=\pm1$ are phantom walls: adjacent rod directions can become parallel there \cite{Chen}. An infinite affine charge is not such a wall. Membership in the open chamber is a statement about an individual metric, not a uniform separation estimate. At fixed distinct directions, separated nodes give charge compactness in Section~\ref{sec:analytic}. Colliding nodes can invalidate that argument; simultaneous charge and length degeneration must be controlled separately.

\subsection{The block determinant restriction}
\begin{theorem}[Block determinant bound]\label{thm:block-bound}
For a balanced metric in the genuine positive even BZ chamber,
\begin{equation}\label{eq:block-bound}
 |\det_{\mathbb Z}(\mathsf v_a,\mathsf v_b)|<b-a
 \quad(b-a\ge2).
\end{equation}
Equivalently, the integral determinant is at most $b-a-1$ in absolute value.
\end{theorem}
\begin{proof}[Proof mechanism]
Let $\rho_i^{\mathrm L}$ and $\rho_i^{\mathrm R}$ be the ratios of the deleted minor at nut $i$ after the prefix flips to rods $a,b$, respectively, to its value at that nut. Polarization of the adjacent square identity gives
\[
 \det(\mathsf r_a,\mathsf r_b)
 =\sum_{i=a+1}^{b}\det(\mathsf r_{i-1},\mathsf r_i)
       \rho_i^{\mathrm L}\rho_i^{\mathrm R}.
\]
The unflipped, sign-adjusted minor is a sum of positive occupation weights. A nonempty partial flip produces a signed sum in which both signs occur. Its absolute value is therefore strictly smaller than the unflipped sum. Consequently each ratio has absolute value at most one, and in every summand of a block containing at least two nuts at least one ratio has absolute value strictly below one. At balance, all adjacent determinants equal $L$. The triangle inequality gives $|\det(\mathsf r_a,\mathsf r_b)|<(b-a)L$.
\end{proof}

For two nuts this yields
\begin{equation}\label{eq:ternary}
 d_j\in\{-1,0,1\}.
\end{equation}
Longer blocks give additional restrictions. For example, the four successive entries $(-1,1,-1,1)$, or their negatives, produce absolute boundary determinant five across five nuts and are excluded. Equal neighbouring nonzero entries are not excluded in general. The full block test is the admissibility condition.

\subsection{The admissible twist interval}
Every singular value of a principal block of the orthogonal matrix $T$ is at most one. Equality $|\det T_j|=1$ for a proper prefix would make that coordinate subspace invariant. The weighted-polynomial projector in \eqref{eq:boundary-dictionary} is irreducible, so
\begin{equation}\label{eq:angular-bound}
 |\mathsf r_j^\phi|<1,\qquad
 \theta_j>2\pi|\mathsf v_j^\phi|\qquad(0<j<N).
\end{equation}
Thus balance requires
\begin{equation}\label{eq:twist-window}
 |p_j-q_j\beta|<1\qquad(0<j<N).
\end{equation}
An interior rod parallel to the exterior rod is impossible: primitiveness would give $(p_j,q_j)=(\pm1,0)$. Other interior directions may repeat.

For a word $W$, define its necessary twist interval
\[
 \mathcal I_W=\{\beta\in(-1,0): |p_j-q_j\beta|<1
                   \text{ for }1\le j<N\}.
\]
It is the intersection of explicit open intervals with rational endpoints. For $N\ge3$, orientation can be chosen so that $d_1=-1$; the first two finite rods then give $-1<\beta<0$. Exterior closure is
$\mathsf v_N=\sigma\mathsf v_0$, $\sigma\in\{-1,1\}$. Applying the block bound at both ends gives
\[
 |q_j|\le\min(j-1,N-j-1)\quad(2\le j\le N-2),
 \qquad |q_1|=|q_{N-1}|=1.
\]
In particular, for even $N\ge4$, reduced wall denominators are at most $N/2-1$. Walls may lie strictly inside $(-1,0)$.

\subsection{Winding determines the reduced sector}\label{subsec:phantoms}
The number of dressing slots is not determined by the number of nuts alone. Consider a positive even representation reduced by removing redundant opposite-sign phantom pairs. Let $h$ be the remaining number of same-sign phantoms and put
\[
 n=N+h,\qquad k=\frac{N-h}{2}.
\]
The following extension uses the positive reduced moment representation; it is not an assertion about arbitrary rational patching matrices.

For ordered genuine nodes and positive weights $\omega_i$, introduce prefix moments
\[
 \mu_q^{[j]}=\sum_{i=1}^N\epsilon_i^{[j]}\omega_i z_i^q,
 \qquad
 \epsilon_i^{[j]}=\begin{cases}-1&i\le j,\\1&i>j.\end{cases}
\]
For $k\ge1$, let
\[
 \mathsf M_j=(\mu_{a+b}^{[j]})_{a,b=0}^{k-1},\qquad
 c_j=(\mu_k^{[j]},\ldots,\mu_{2k-1}^{[j]})^T,\qquad
 e_k=(0,\ldots,0,1)^T.
\]
After normalizing the positive parameter at infinity to one, the two rod components, before a common frame change and common scale, are
\begin{equation}\label{eq:reduced-moment-rods}
 \mathsf t_j=\left(
 \det\begin{pmatrix}
 \mathsf M_j&c_j+e_k\\
 (c_j+e_k)^T&\mu_{2k}^{[j]}+\sum_i z_i
 \end{pmatrix},
 \det\mathsf M_j\right).
\end{equation}
Here $\mathsf t_j$ is a two-component rod, not the prefix column $t_i$ used in \eqref{eq:nut-square}. These formulas follow by taking the same-sign phantom limit of the homogeneous axis sectors before removing the common phantom factors.

The polarized identity in these coordinates is
\[
 \det(\mathsf t_a,\mathsf t_b)
   =2\sum_{i=a+1}^b\omega_i\Psi_a(z_i)\Psi_b(z_i),\qquad
 \Psi_j(x)=\det\begin{pmatrix}
 \mathsf M_j&c_j+e_k\\
 (1,x,\ldots,x^{k-1})&x^k
 \end{pmatrix}.
\]
At nut $i$, the two adjacent values of $\Psi$ agree. Its sign-adjusted expansion contains positive terms of subset sizes $k-1$ and $k$. For $k\ge1$ both partial-flip parities occur, giving the same strict block bound. Likewise Cauchy--Binet gives
\[
 \det\mathsf M_j=\sum_{|S|=k}
 (-1)^{|S\cap\{1,\ldots,j\}|}
 \prod_{i\in S}\omega_i\prod_{\substack{a<b\\a,b\in S}}(z_a-z_b)^2.
\]
Both signs occur for a proper prefix, so $|\det\mathsf M_j|<\det\mathsf M_0$. The angular bound and the ternary restriction therefore persist in these reduced sectors.

The maximal case $k=0$ is different: the reduced rods are
$\mathsf t_j=(\mu_0^{[j]}+\sum_i z_i,1)$, and
\[
 \det(\mathsf t_0,\mathsf t_N)=2\sum_i\omega_i>0.
\]
Their exterior lines do not agree, so this case is excluded by AF closure. It must not be excluded from the ALF problem by an AF argument.

\begin{proposition}[AF phantom restriction]\label{prop:phantom-bound}
In the positive reduced representation above, a regular AF metric with $N\ge3$ genuine nuts satisfies
\[
 h\equiv N\pmod2,\qquad 3h\le N.
\]
Equality occurs only for $N=3$, $h=1$, with word $(-1,-1)$ up to orientation. In particular $h<N/3$ for $N>3$.
\end{proposition}
\begin{proof}[Winding argument]
Let $\mathfrak w$ be the projective winding of the signed rod path: the sum of its successive positive angles, each in $(0,\pi)$, is $\pi\mathfrak w$. AF closure makes $\mathfrak w$ an integer. Writing $f_i=(1,z_i,\ldots,z_i^{k-1})^T$ gives
\[
 \mathsf M_j=\mathsf M_0-2\sum_{i\le j}\omega_i f_i f_i^T,\qquad
 \mathsf M_0>0,\quad \mathsf M_N=-\mathsf M_0.
\]
Continuous prefix flips give exactly $k$ eigenvalue crossings. The two determinants in \eqref{eq:reduced-moment-rods} are affine during each rank-one flip, so the rod follows the straight segment between its endpoints. Positive adjacent determinant prevents it from passing through zero. Hence $\mathfrak w=k$ in this orientation; orientation reversal gives $N-k$. Therefore $h=|N-2\mathfrak w|$.

Write $\mathsf v_N=\sigma\mathsf v_0$ and
$\mathsf v_{N-1}=(p,-\sigma)$ in lattice coordinates, and put $c=\sigma p$.
The angular bound gives $c\in\{0,1\}$. Closing the signed sequence by its $\sigma$-multiple and applying the unimodular rotation formula \cite{HigashitaniMasuda} gives
\[
 3h=\left|\sum_{j=1}^{N-1}d_j-c\right|\le N.
\]
Equality forces every $d_j=-1$ and $c=1$. That recurrence reaches
$\mathsf v_3=-\mathsf v_0$, which is a forbidden interior return unless $N=3$. Parity follows from $h=N-2k$.
\end{proof}

For $N=4,6$ this leaves only $h=0$. The three-nut Chen--Teo word has $h=1$ and is therefore included, not excluded. These are necessary sector restrictions. The reduced-rank reconstruction in Section~\ref{subsec:reduced-reconstruction} also supplies the matrix identities used in the cone-convexity theorem.

\subsection{Enumeration of admissible words}\label{subsec:classification}
We call a word \emph{admissible} when it passes the necessary tests in this section in its stated positive BZ sector. This term does not assert metric existence. Enumeration proceeds as follows:
\begin{enumerate}
\item Generate the ternary words with $d_1=-1$, reconstruct their signed primitive rods by \eqref{eq:rod-recurrence}, and impose AF exterior closure.
\item Check every block determinant, exclude interior returns to the exterior direction, and retain the words with $\mathcal I_W\ne\varnothing$.
\item Compute the projective winding and the allowed reduced dressing sector. For the genuine even sector require $\mathfrak w=N/2$.
\item Identify axis-reversal and orientation equivalences, carrying along their transformations of the marking and of $\beta$.
\end{enumerate}
The last step is a quotient of the candidate list, not an additional regularity inequality or a symmetry imposed on a metric. For an oriented classification, retain the orientation labels. With $d_1=-1$, normalized reversal takes $W$ to $\operatorname{reverse}W$ and $\beta$ to $-1-\beta$ when $d_{N-1}=-1$; when $d_{N-1}=1$, it takes $W$ to $-\operatorname{reverse}W$ and leaves $\beta$ unchanged.

For example, the six-nut tests leave precisely
\[
 (-1,0,0,0,1),\qquad
 (-1,1,0,-1,1),\qquad
 (-1,1,1,1,-1),
\]
all with necessary interval $(-1,0)$ and $h=0$. These three words satisfy the collision criterion of Theorem~\ref{thm:collision-existence}; their certified balance data are included in the appendix.

The enumeration classifies integral candidates in the specified positive BZ sectors. The distinction between admissibility and existence is essential: the collision criterion in Section~\ref{sec:analytic} is sufficient for balance, whereas the boundary-gradient obstruction excludes certain admissible words. Both tests depend on the normalized cone family, in addition to the integral recurrence. A common self-intersection block need not have the same obstruction in different ambient words, because its exterior cone normalization must be matched to the surviving geometry.

\section{The cone potential}\label{sec:analytic}

The regularity problem becomes a variational problem once the charges have been solved as functions of the rod lengths and directions. We first establish these coordinates and the positive Hessian, then analyze finite collision faces and large lengths. The algebraic statements use the positive reduced reconstruction of Section~\ref{subsec:reduced-reconstruction}; applications of canonical cone limits have the additional hypotheses stated in Section~\ref{subsec:comparison}.

\subsection{Unique direction matching at separated nodes}
Fix ordered distinct nodes. Multiplying every positive weight $\omega_i$ by the same scalar $s>0$ leaves $T$ unchanged. In \eqref{eq:reduced-dictionary}, it sends $a\mapsto s^{-1/2}a$ and $c\mapsto s^{1/2}c$. The AF alignment $u^Tv=0$ therefore determines this common scale uniquely. The AF charge domain is parametrized by $N-1$ logarithmic weight ratios and is connected.

Fix the first rod direction and orient the normalized rods so that adjacent determinants are positive. Let $\alpha_i$ denote the counterclockwise angle from rod $i-1$ to rod $i$. These successive angular increments satisfy
\begin{equation}\label{eq:angle-simplex}
 0<\alpha_i<\pi,\qquad \sum_{i=1}^{N}\alpha_i=k\pi.
\end{equation}
The prefix moment inertia calculation in Section~\ref{subsec:phantoms} gives winding $k$. For the genuine even sector, $N=2m$ and $k=m$. The set \eqref{eq:angle-simplex} is a convex domain of dimension $N-1$.

\begin{proposition}[Direction coordinates]\label{prop:direction-coordinates}
At fixed ordered nodes, the AF positive-charge domain maps analytically and bijectively onto the direction domain \eqref{eq:angle-simplex}. The inverse is analytic. In particular, for fixed directions in this domain and fixed $L$, the direction-matched BZ data form one analytic family over the full positive length orthant.
\end{proposition}
\begin{proof}
Lemma~\ref{lem:paired-transport} supplies the coupled triangular transport, including when a principal block of $T$ is singular. The positive identity proved below, applied with $\delta Z=0$, then makes every fixed-direction charge variation vanish. This gives local invertibility. Properness has a particularly useful geometric proof. The identity $[Z+vv^T,T]=0$ and AF alignment give $\|v\|^2\le z_N-z_1$, so residue data have convergent subsequences on compact sets of separated nodes.

A normalized rod cannot vanish in such a limit. Otherwise a nonzero vector $h$ supported in a prefix would satisfy $T_jh=0$ and $v_{[j]}^Th=0$. Its image $Th$ is supported in the complementary suffix and has the same norm. Commutation with $Z+vv^T$ would then equate a quadratic form bounded above by $z_j\|h\|^2$ with one bounded below by $z_{j+1}\|h\|^2$, contradicting separation.

If relative weights escaped, the positive Cauchy--Binet sums for the polynomial projector would force a diagonal projector entry to $0$ or $1$. The corresponding row of $T$ decouples, the displacement identity forces the residue to vanish, and adjacent limiting rods become parallel. This contradicts the prescribed distinct directions. The AF scalar normalization controls the remaining common weight scale. Thus the map is proper. A proper local diffeomorphism onto the connected convex domain \eqref{eq:angle-simplex} is a one-sheeted covering, proving the claim.
\end{proof}

A marked word belongs to this direction domain precisely when its positive signed lift has projective winding $\mathfrak w=k$. Orientation reversal may be used to choose $k\le N/2$. No restriction on individual positive lengths is imposed.

The proposition prevents charge escape on compact sets of separated positive lengths. Collision limits require the additional matching stated below.

\subsection{Convexity from the mass differential}
By Lemma~\ref{lem:paired-transport}, every fixed-direction variation of the paired boundary data $(T,v)$ has a unique lower triangular transport matrix $K$ satisfying
\begin{equation}\label{eq:triangular-transport}
 \delta T=KT+TK^T,\qquad \delta v=Kv.
\end{equation}
If $\mathcal L$ is a lower triangular change of boundary frame and $\mathcal L_j$ its leading block, congruence multiplies both components of the $j$th normalized rod by $(\det\mathcal L_j)^2$. Let $K_j$ denote the leading block of the infinitesimal transport $K$. With $S=K+K^T$,
\begin{equation}\label{eq:cone-transport}
 ST+TS=0,\qquad \delta u=-K^Tu,\qquad
 \delta b_j=2\tr K_j.
\end{equation}
The coupled transport is essential. The map $K\mapsto KT+TK^T$ on lower triangular matrices need not be invertible when a leading principal minor vanishes. As Andersson pointed out \cite{AnderssonCorrection}, this occurs at a finite rod with zero angular component, $p_j-\beta q_j=0$, and cannot be excluded by bulk positivity. Appendix~\ref{app:paired-transport} instead transports the auxiliary pencil $T+c vv^T$, for a fixed generic real $c$. The nonzero adjacent-rod determinants then force the same transport to carry $v$, recovering \eqref{eq:triangular-transport} for the original data. Thus the argument below also applies at these special twists.

Since $Z+vv^T$ commutes with $T$ and $S$ anticommutes with it, $\tr((Z+vv^T)S)=0$. Using \eqref{eq:residue-moments} and summation by parts gives the first law at fixed $L,\beta$:
\begin{equation}\label{eq:mass-cone-law}
 2\,\delta M=v^TSv=-2\sum_i z_iK_{ii}
                   =\sum_{j=1}^{N-1}\ell_j\,\delta b_j.
\end{equation}
It follows that
\begin{equation}\label{eq:cone-potential}
 \Phi_W(\ell;\beta)=\sum_j\ell_jb_j-2M,\qquad
 d_\ell\Phi_W=\sum_jb_j\,d\ell_j.
\end{equation}
The subscript in $d_\ell$ means that only the lengths vary; the period parameters $L$ and $\beta$ are held fixed. The identity makes $D_\ell b$ symmetric.

\begin{theorem}[Strict cone convexity]\label{thm:convexity}
Fix the marked directions in \eqref{eq:angle-simplex}, the asymptotic scale and the exterior gauge. On the full positive length orthant, $\Phi_W$ is real analytic and strictly convex. Its Hessian $H=D_\ell b$ satisfies
\begin{equation}\label{eq:positive-Hessian}
 \delta\ell^TH\delta\ell
 =\sum_{i>j}(z_i-z_j)K_{ij}^2+\|(K+K^T)v\|^2>0
\end{equation}
for every nonzero length variation.
\end{theorem}
\begin{proof}[The positive identity]
Differentiate $[Z+vv^T,T]=0$ and set
\[
 F=\delta Z+[Z,K]+vv^TS.
\]
One obtains $FT=TF^T$. Together with $TST=-S$, this implies $\tr(FS)=0$. Expanding the trace and using lower triangularity gives
\[
 -2\tr(\delta ZK)
 =\sum_{i>j}(z_i-z_j)K_{ij}^2+\|Sv\|^2.
\]
By \eqref{eq:cone-transport}, the left side is $\sum_j\delta\ell_j\delta b_j$. All coefficients on the right are positive. In the equality case $K$ is diagonal, anticommutes with $T$, and satisfies $Kv=0$. Hence $K^2$ is a diagonal matrix commuting with $T$. Irreducibility from Section~\ref{subsec:reduced-reconstruction} makes $K^2$ scalar, and $v\ne0$ makes that scalar zero. Thus $K=0$. The remaining displacement identity forces $\delta Z$ to be scalar. This is translation, already removed by passing to the lengths. Thus equality occurs only for $\delta\ell=0$.
\end{proof}

All $N-1$ independent lengths are varied. Fixing $L$ is essential: if scale is allowed to vary, homogeneity gives the null direction of a common rescaling of $\ell$ and $L$.

\subsection{Mass as the dual potential}
Strict convexity on a convex domain gives at most one length vector with any prescribed attainable cone vector. In particular, a balanced point is unique and is the global minimum of $\Phi_W$.

The Legendre transform has a simple geometric interpretation. On the realized cone image,
\begin{equation}\label{eq:mass-dual}
 \Phi_W^*(b)=2M,\qquad \nabla_b(2M)=\ell,\qquad
 D_b^2(2M)=H^{-1}>0.
\end{equation}
Lengths and logarithmic cone amplitudes are dual coordinates. This does not assert that every cone vector is realized. It also does not say that the mass is minimized as a function of lengths. In fact,
\begin{equation}\label{eq:mass-length-gradient}
 2\nabla_\ell M=H\ell,\qquad
 \ell\cdot\nabla_\ell M=\tfrac12\ell^TH\ell>0.
\end{equation}
The mass minima on the balanced twist curves in the appendix concern a different variation.

\subsection{Canonical identification}\label{subsec:comparison}
The positive BZ reconstruction solves the canonical conical problem at separated nodes. This provides the geometric interpretation needed to apply the cone-limit results of \cite{KL,LiSun}.

\begin{proposition}[Canonical identification]\label{prop:canonical-comparison}
In every positive reduced moment sector of Section~\ref{subsec:reduced-reconstruction}, the AF-aligned BZ metric with distinct genuine nodes is the canonical conical harmonic-map metric for its marked rod directions and normalized end. For the direction lift in Proposition~\ref{prop:direction-coordinates}, the identification holds over the full positive length orthant, including rods with zero angular component.
\end{proposition}
\begin{proof}[Proof outline]
The augmented bulk matrix stays positive at the axis and nuts. Its analytic rod and nut expansions place the harmonic map at bounded distance from the canonical local models. The AF expansion fixes the normalization at infinity, so uniqueness of the canonical harmonic map identifies the metrics. Appendix~\ref{app:canonical-proof} gives the expansions and comparison argument. The same identification follows from the positive-family dictionary in Section~\ref{subsec:ADN-dictionary} and \cite[Propositions~8.9 and 10.2]{ADN}.
\end{proof}

\begin{corollary}[Smooth completion at balance]\label{cor:smooth-completion}
A balanced positive reduced solution with distinct genuine nodes, primitive rods and unit adjacent lattice determinants defines a smooth complete AF gravitational instanton.
\end{corollary}
Cone balance removes the axis singularities, adjacent unimodularity removes the nut quotients, and the normalized end gives completeness and finite curvature energy. This also follows from \cite[Proposition~11.4]{ADN}. These are separated-node statements; they do not identify every collision limit with a smaller BZ model.

\subsection{Opening a collision}\label{subsec:collision-matching}
Consider a consecutive cluster word $C$ whose two exterior primitive directions are independent. Solve its homogeneous moment model with both exterior cone angles fixed to $2\pi$. The model has no prescribed finite asymptotic length scale. Write $\xi_j>0$ for its finite gaps and $b_{C,j}$ for its internal logarithmic cone errors. The model and its moment normalization are specified in Appendix~\ref{app:homogeneous-model}.

The homogeneous first law gives the convex potential
\begin{equation}\label{eq:cluster-potential}
 \Psi_C(\xi)=\sum_j\xi_jb_{C,j}(\xi),\qquad
 d\Psi_C=\sum_jb_{C,j}\,d\xi_j,\qquad
 \mu(C)=\inf_{\xi_j>0,\ \sum_j\xi_j=1}\Psi_C(\xi).
\end{equation}
The normalization removes the common scale; the infimum allows further collisions among the cluster's own rods. The value of $\mu(C)$ belongs to this marked homogeneous problem. Using it inside a larger geometry additionally requires matching the exterior cone normalization.

Only two facts about this invariant are needed here. A negative value supplies a positive shape with negative opening cost. Its sign can be certified at any positive direction-matched shape $\xi$ with $\sum \xi_j=1$, because convexity and homogeneity give
\begin{equation}\label{eq:cluster-trial-bounds}
 \min_j b_{C,j}(\xi)\le\mu(C)\le\sum_j\xi_jb_{C,j}(\xi).
\end{equation}
Thus one negative weighted trial proves negativity, whereas a trial with every error positive proves positivity. At a positive constrained stationary point all errors are equal, and their common value is $\mu(C)$.

\begin{assumption}[Homogeneous collision matching]\label{hyp:comparison}
For each nonparallel collision block used in the existence criterion, assume the following.
\begin{enumerate}
\item At a collision with balanced surviving rods, the cluster cones converge to the homogeneous model with the inherited primitive periods and both exterior angles equal to $2\pi$. If $q$ is the boundary length vector and $\iota_C \xi$ inserts the cluster gaps into the full length vector, then
\[
 b_j(q+\varepsilon\iota_C \xi)\longrightarrow b_{C,j}(\xi),\qquad j\in C,
\]
locally uniformly in positive shapes and interior twist. Simultaneous disjoint blocks use the same normalization.
\item At fixed distinct cluster nodes, the homogeneous direction map is locally invertible in positive weight ratios after fixing the exterior frame.
\end{enumerate}
\end{assumption}
The first condition concerns the amplitudes of the rods, not just their projective directions. The fixed-node surjectivity theorem of \cite[Section~9]{ADN} does not prove either collision condition. Proposition~\ref{prop:canonical-comparison} identifies the separated AF family but does not by itself prove this normalized collision matching.

\subsection{Existence from boundary descent}\label{subsec:collision-existence}
Set $L=1$ and extend $\Phi_W$ lower-semicontinuously to the closed length orthant. A finite minimum on a boundary stratum has balanced surviving rods, since their positive lengths can vary in both directions. Under Hypotheses~\ref{hyp:comparison}, opening a nonparallel block in shape $\xi_C$ has cost
\begin{equation}\label{eq:opening-cost}
 \lim_{\varepsilon\downarrow0}\sum_{j\in C}\xi_j b_j(q+\varepsilon\iota_C \xi_C)=\Psi_C(\xi_C).
\end{equation}
A negative trial therefore rules out that boundary minimum.

\begin{theorem}[Existence from collision descent]\label{thm:collision-existence}
Let $W$ be an admissible AF word in a positive reduced sector, and assume Hypotheses~\ref{hyp:comparison}. If every nonparallel proper consecutive block $C$ has $\mu(C)<0$, then at every $\beta\in\mathcal I_W$ and fixed $L>0$ there is exactly one balanced positive length vector. It defines a smooth complete AF metric and depends real-analytically on the interior parameters.
\end{theorem}
\begin{proof}[Proof outline]
Negative opening costs exclude nonparallel collision minima. For parallel blocks a cone error tends to $-\infty$, while the other errors remain bounded above at bounded total span. The large-span cone limit and the strict angular inequalities exclude escape to infinity. Minimizers on closed truncated orthants therefore converge to an interior minimum. Strict convexity gives uniqueness, and the positive Hessian gives analytic dependence. The bounds and compactness argument are recorded in Appendix~\ref{app:collision-proof}.
\end{proof}
Consequently the balanced lengths have positive lower and finite upper bounds on every compact subinterval of $\mathcal I_W$. No such bound at the twist endpoints is asserted. Finite supporting-plane certificates can give effective enclosures independently of the collision analysis; see Appendix~\ref{app:finite-certificates}.

\subsection{Refinement by zero pairs}\label{subsec:zero-insertion}
Inserting two adjacent zero self-intersections replaces the consecutive signed directions $\mathsf v_-,\mathsf v_+$ by
\begin{equation}\label{eq:zero-rod-operation}
 \mathsf v_-,\quad\mathsf v_+,\quad-\mathsf v_-,\quad-\mathsf v_+,
\end{equation}
and reverses the signs of the remaining suffix. The set of vanishing circles is unchanged, while the ordered walk acquires two genuine turning points. This operation is different from adding phantom slots.

\begin{lemma}[Positive zero-pair insertion]\label{lem:zero-insertion}
Let $D$ be a nonparallel homogeneous word with positive direction-matched data, and let $D^\sharp$ be obtained by a zero-pair insertion. Under the local direction-coordinate condition in Hypotheses~\ref{hyp:comparison}, the insertion has positive exactly direction-matched realizations at small separation and
\[
 \mu(D^\sharp)\le\mu(D).
\]
At an interior nut, a symmetric opening of total length $\varepsilon$, compensated by shortening the two adjacent rods equally, satisfies
\begin{equation}\label{eq:zero-opening-derivative}
 \Psi_{D^\sharp}(\xi_\varepsilon)=\Psi_D(\xi)-\varepsilon\log2+o(\varepsilon).
\end{equation}
\end{lemma}
The proof splits one positive moment node into three and evaluates the two consecutive Hankel determinants by a Schur complement; see Appendix~\ref{app:insertion-proof}. The inequality preserves a negative sign under refinement. It does not preserve positivity, and does not prove that every sufficiently long refinement becomes negative.

\subsection{Families obtained by extending zero strings}\label{subsec:families}
The following extension preserves both the direction arrangement and the asymptotic twist interval.
\begin{theorem}[Zero-string extension]\label{thm:zero-extension}
Suppose an admissible AF word $W$ has negative invariant for every nonparallel proper consecutive block. Replace any selected zero entries by $0^{2r+1}$, independently, where $r\ge0$ and $0^a$ denotes a string of $a$ zeros. The resulting word remains admissible, has the same strict twist interval and reduced phantom count, and again has negative invariant for every nonparallel proper block. Its full-interval existence follows from Theorem~\ref{thm:collision-existence} under Hypotheses~\ref{hyp:comparison} for that family.
\end{theorem}
\begin{proof}
A zero pair repeats existing directions; it increases the nut count by two and the projective winding by one, leaving $h=|N-2\mathfrak w|$ unchanged. Block determinants and angular bounds are preserved. Delete added zero pairs from any proper block. Unless it consists entirely of zeros, the result is a proper old block, so Lemma~\ref{lem:zero-insertion} preserves negativity. Even all-zero blocks have a negative Weyl trial, while odd all-zero blocks have parallel exteriors. The elementary negative tests are proved in Appendix~\ref{app:negative-tests}.
\end{proof}

\begin{theorem}[Explicit infinite families]\label{thm:families}
For all nonnegative integers $r,s$, the words in Table~\ref{tab:analytic-families} satisfy the proper-block collision criterion. Their balanced positive realizations exist uniquely throughout the displayed intervals under Hypotheses~\ref{hyp:comparison}.
\end{theorem}
\begin{table}[htbp]\centering\small
\begin{tabular}{@{}lccc@{}}\toprule
Self-intersection word & $N$ & $\mathcal I_W$ & $h$\\\midrule
$(-1,0^{2r+1},1)$ & $2r+4$ & $(-1,0)$ & $0$\\
$(-1,1,0^{2r+1},-1,1)$ & $2r+6$ & $(-1,0)$ & $0$\\
$(-1,0^{2r+2},1,0^{2s+1},-1,1)$ & $2r+2s+8$ & $(-1,0)$ & $0$\\
$(-1,1,-1,0^{2r+1},1,-1,1)$ & $2r+8$ & $(-1,-1/3)$ & $0$\\
$(-1,0^{2r+1},-1,-1,-1,0^{2s+1},-1)$ & $2r+2s+8$ & $(-1,0)$ & $2$\\\bottomrule
\end{tabular}
\caption{Zero-string families. Here $N$ counts genuine turning points and $h$ is the reduced number of same-sign phantom slots.}\label{tab:analytic-families}
\end{table}
Their finite seed checks and a direct reduction for the second family appear in Appendix~\ref{app:family-proof}. The first family is the Li--Sun family, whose canonical AF existence also follows from \cite[Theorem~6.4]{LiSun}. The additional families use the explicit matching hypotheses stated here.

\subsection{A positive boundary obstruction}\label{subsec:boundary-obstruction}
The potential also provides a nonexistence test that will distinguish walks on the same root arrangement.
\begin{proposition}[Positive boundary gradient]\label{prop:boundary-obstruction}
Suppose a direction-matched family $q_\varepsilon>0$ satisfies
\[
 q_\varepsilon\to q_0,\qquad b(q_\varepsilon)\to g,\qquad
 g\ge0,\quad g\ne0,\quad g\cdot q_0=0.
\]
Then the positive length orthant contains no balanced point.
\end{proposition}
\begin{proof}
If $b(x)=0$, convex gradient monotonicity gives $b(q_\varepsilon)\cdot(q_\varepsilon-x)\ge0$. The limit implies $g\cdot x\le0$, contradicting $x>0$ and the stated properties of $g$.
\end{proof}
The condition refers to a matched degeneration of the full word. A positive invariant for an isolated subword is not sufficient to apply it in an arbitrary embedding.

\subsection{Variation along a balanced family}
Let $\overline M=M/L$ and $\overline J=J/L^2$ be the dimensionless mass and angular momentum. The extended first law, derived in Appendix~\ref{app:twist-derivatives}, gives
\[
 \frac{d\overline M}{d\beta}=2\overline J
\]
on a balanced branch at fixed $L$. A zero of angular momentum is therefore a stationary point of its mass curve. Strict convexity of the length potential does not alone make this stationary point a minimum as $\beta$ varies. The numerical curves in Appendix~\ref{app:mass} concern this separate variation.

\section{Classification through generalized root systems}\label{sec:classification-roots}

The rod restrictions organize the candidates by a finite set of lattice directions. This section identifies that set with a rank-two crystallographic arrangement, explains its polygon description, and tests what this description retains about cone balance. The AF and ALF catalogues test which information this organization retains; Section~\ref{sec:endpoint-composition} gives the complementary endpoint description. The ordered walk remains part of the data: distinct metrics can revisit the same directions in different orders and with different multiplicities.

\subsection{From rod directions to a root arrangement}
Use coefficient vectors $v_j=(p_j,q_j)\in\mathbb Z^2$ for the primitive rods in the marked basis $(A,-B)$. The associated signed direction set is
\[
 \mathcal R(W)=\{\pm v_j:0\le j\le N\},
\]
with repetitions removed. In two dimensions, a \emph{crystallographic arrangement} can be described by a finite centrally symmetric set of primitive vectors for which every pair of angularly consecutive rays forms an integral basis. It determines a finite generalized root system of a Weyl groupoid \cite{CuntzArrangements}. Here the simple-root basis changes with the chosen chamber; the groupoid records the reflection maps between these bases. A single reflection group for one fixed inner product need not preserve the entire arrangement.

\begin{proposition}[The direction arrangement]\label{prop:root-arrangement}
Let a primitive unimodular rod walk start with $d_1=-1$ and have $d_j\in\{-1,0,1\}$. Then its distinct signed directions form a rank-two crystallographic arrangement. The conclusion concerns the integral directions and does not require cone balance.
\end{proposition}
\begin{proof}
The first three directions are $v_0,v_1,v_1-v_0$, which form a complete unimodular arrangement after adjoining their negatives. A zero entry in the recurrence repeats a line; an entry $\pm1$ produces a sum or difference of the preceding basis vectors. Choose signs so the proposed new ray is $u+v$ inside the unimodular cone generated by $u,v$.

Any nontrivial unimodular subdivision of that cone already contains $u+v$. Indeed, in coordinates $u=(1,0)$, $v=(0,1)$, a subdivision omitting $(1,1)$ would have adjacent rays crossing the diagonal. Their determinant is at least two, contradicting unimodularity. Thus either $u+v$ is already present, or the cone contains no old ray and insertion of $u+v$ preserves its two adjacent determinants. Adjoin the opposite ray at the same time and proceed by induction.
\end{proof}

The hypothesis $d_1=-1$ is the catalogue convention. The two-direction all-zero models belong to the reducible case and are not among these AF candidates. The proposition explains the route from the continuous construction to the discrete organization: positive BZ minors impose the ternary recurrence, and that recurrence generates the crystallographic arrangement. Purity of an arbitrary spinor, without the positivity and integral-period conditions, does not imply this conclusion.

\subsection{Triangulated polygons}
Let $s$ be the number of unoriented directions in $\mathcal R(W)$ and choose representatives $\alpha_1,\ldots,\alpha_s$ in angular order over a half-turn. Continue them by $\alpha_0=-\alpha_s$ and $\alpha_{s+1}=-\alpha_1$. Consecutive determinants are one, so there are positive integers $c_i$ satisfying
\begin{equation}\label{eq:polygon-cycle}
 \alpha_{i-1}+\alpha_{i+1}=c_i\alpha_i,\qquad
 c_i=\det(\alpha_{i-1},\alpha_{i+1}).
\end{equation}
The cyclic sequence $(c_1,\ldots,c_s)$ records the numbers of triangles incident at the vertices of a triangulated convex $s$-gon. The correspondence with finite rank-two Weyl groupoids is described in \cite{CuntzHeckenberger}. Adding the sum of two neighboring directions corresponds to attaching an ear to the triangulation.

\begin{figure}[htbp]\centering
\includegraphics[width=\textwidth]{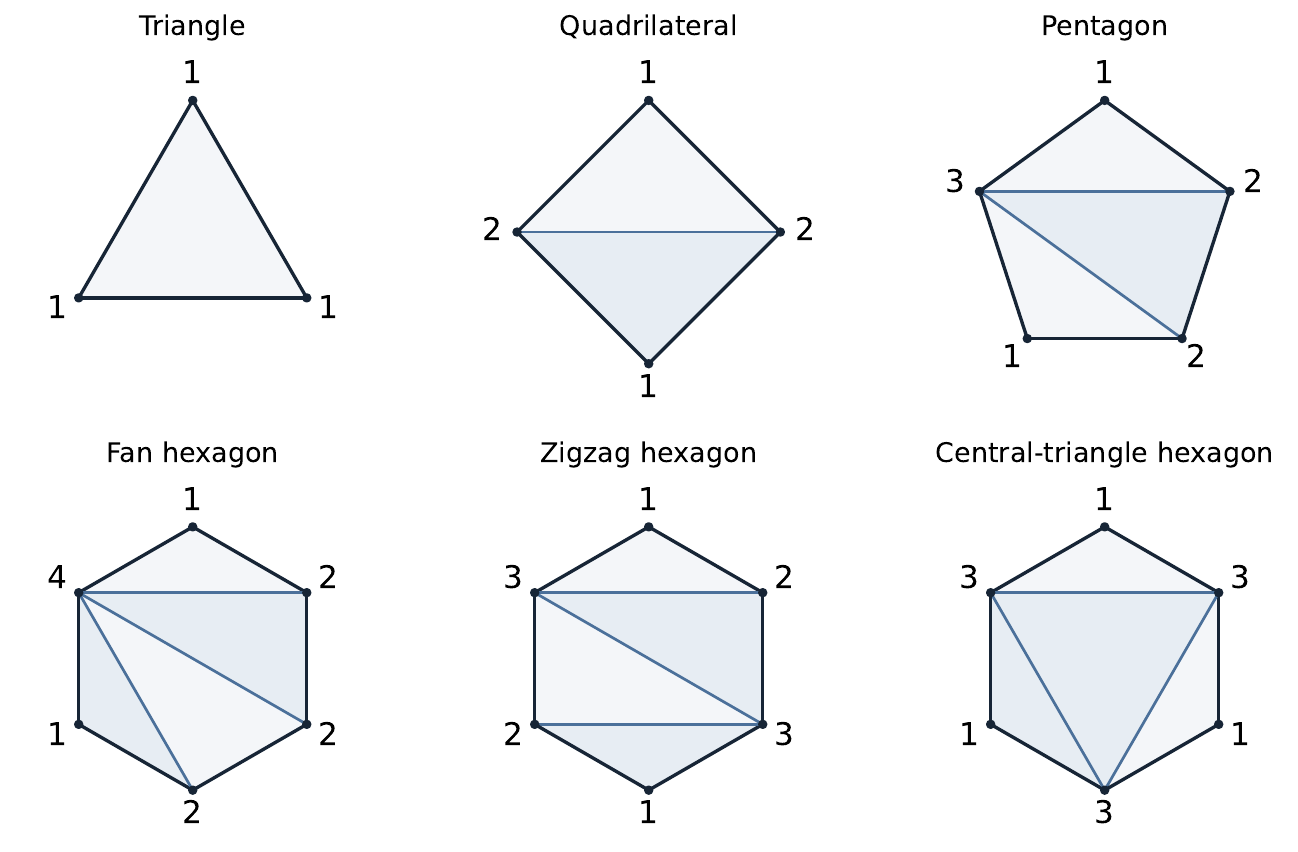}
\caption{The six triangulation types in the AF census through $N=12$. Vertex labels are the coefficients $c_i$ of the sorted direction arrangement. The fan, zigzag and central-triangle hexagons have the same number of directions but different triangulations.}\label{fig:root-polygons}
\end{figure}

The coefficients $c_i$ are not the self-intersections $d_j$. They use angularly sorted distinct directions, whereas the $d_j$ use the full sequence along the orbit-space axis. In particular, the triangulated polygon in Figure~\ref{fig:root-polygons} is a combinatorial encoding; it is neither the orbit space nor the rod polygon obtained by following the original ordered vectors.

\subsection{The AF census by arrangement}
Table~\ref{tab:AF-root-types} groups all 359 candidates through twelve genuine turning points by the complete direction set, including the exterior. Exactly 356 have certified midpoint balance roots; the other three are excluded by Proposition~\ref{prop:boundary-obstruction}. The analytic block tests verify the conditional full-interval criterion for the 356 retained words.

\begin{table}[htbp]\centering\small
\begin{tabular}{@{}lcrrr@{}}\toprule
Triangulation & Cycle $(c_i)$ & Candidates & Balance roots & Excluded\\\midrule
Triangle & $(1,1,1)$ &10&10&0\\
Quadrilateral & $(1,2,1,2)$ &187&187&0\\
Pentagon & $(1,2,2,1,3)$ &152&152&0\\
Fan hexagon & $(1,2,2,2,1,4)$ &6&6&0\\
Zigzag hexagon & $(1,2,3,1,2,3)$ &3&0&3\\
Central-triangle hexagon & $(1,3,1,3,1,3)$ &1&1&0\\\midrule
Total&&359&356&3\\\bottomrule
\end{tabular}
\caption{AF candidates for $3\le N\le12$, modulo the stated reversal and orientation equivalences. A balance-root entry denotes an interval certificate at one interior twist.}\label{tab:AF-root-types}
\end{table}

The triangle is the Li--Sun direction arrangement; within this census it has one word at each $N$. The triangle and quadrilateral rows are exactly the 197 full-window words with $\mathcal I_W=(-1,0)$. The remaining 159 retained words have restricted windows. This separation is an observed equality in the finite census, not an all-order classification theorem.

The usual root-system names depend on whether the exterior is included. The triangle gives $A_2$, the quadrilateral gives $B_2$, and the central-triangle hexagon gives $G_2$. The six signed interior vectors $\pm B$, $\pm(A+B)$ and $\pm(A+2B)$ form an $A_2$ subsystem; adjoining the exterior $\pm A$ completes the $B_2$ arrangement. Generalized root systems supply a common description beyond these classical examples.

Fixing an arrangement also turns word enumeration into a finite walk problem. Its states are ordered root pairs $(u,v)$ with determinant one, and the allowed transition is $(u,v)\mapsto(v,-u-dv)$ with $d\in\{-1,0,1\}$, subject to the end conditions. This gives exact linear recurrences for word counts. For the marked $G_2$ root set, the canonical counts at $N=6,8,10,12$ are respectively $3,9,40,199$, agreeing with independent enumeration. These count admissible words before the cone-balance test.

An independent one-point search at $N=14$ found convergent representatives for 1,159 of 1,225 candidates, including all 362 quadrilateral and 686 pentagon words. The 66 remaining searches were zigzag words. Independent high-precision moment evaluation checked the convergent points, but this search supplied neither interval certificates for the whole set nor nonexistence proofs for the failures.

\subsection{Different walks on the same polygon}
The three excluded twelve-nut words are \texttt{12-065}, \texttt{12-066} and \texttt{12-145}. Their matched degenerations contain the six-nut block
\[
 E=(-1,0,-1,1,-1),\qquad \mu(E)=0.09003150665\ldots>0.
\]
Appendix~\ref{app:exceptional} proves the positive full boundary gradient required for exclusion. The common zigzag polygon does not by itself cause the obstruction.

To see this, define $E^\sharp=(-1,0,0,0,-1,1,-1)$ by extending the zero string, and let $E^R$ denote reversal of a word. The two AF words
\[
 W_0=(E,0,-E^R),\qquad
 W_1=(E^\sharp,0,-(E^\sharp)^R)
\]
have the same zigzag arrangement and necessary window $(-3/5,-1/5)$. The first is the excluded word \texttt{12-066}. The second has sixteen genuine nuts and an interval-certified positive balance root at $\beta=-2/5$. Its cluster satisfies $\mu(E^\sharp)=-0.01811463686\ldots<0$. Zero insertion changes the cone problem while preserving the set of directions.

\begin{table}[htbp]\centering\small
\begin{tabular}{@{}lrrrrrr@{}}\toprule
Direction & $A$ & $B$ & $A+B$ & $A+2B$ & $A+3B$ & $2A+5B$\\\midrule
Occurrences in $W_0$ &2&4&2&2&2&1\\
Occurrences in $W_1$ &2&6&4&2&2&1\\\bottomrule
\end{tabular}
\caption{The same six unoriented directions with different multiplicities. The exterior rods are included, so the totals are $N+1$.}\label{tab:zigzag-walks}
\end{table}
The root certificate for $W_1$ and its original-BZ consistency check are included in the supplement. Generic initial guesses failed for this word; a seed constructed by splitting the conical $W_0$ data converged. This provides a concrete reason not to interpret failed searches as exclusions.

\subsection{A realizability conjecture}
The arrangement supplies a finite alphabet of directions. An ordered walk specifies the rod topology, and the asymptotic marking fixes which period data are to be balanced. The remaining question is which of these marked walks have an interior minimum of the cone potential.

\begin{conjecture}[AF polygon realizability]\label{conj:polygon-realizability}
Every rank-two crystallographic arrangement occurring in an admissible AF word in the positive reduced sectors occurs in an admissible marked word with a regular positive BZ realization at some interior twist.
\end{conjecture}
The realizing walk may differ from the word that first exhibits the arrangement. The conjecture does not assert that every walk is regular, nor that every sufficiently long zero extension works. Zero-pair insertion provides a candidate construction because it preserves the arrangement and can lower the opening cost. An arbitrary arrangement may support many words, and a general realization theorem remains open.

There is a related geometric description of the direction data: the ordered matrix $(\alpha_1\ \cdots\ \alpha_s)$ has adjacent minors equal to one. The normalized minors belong to the Grassmannian framework underlying the polygon correspondence \cite[Section~4]{CuntzHeckenberger}. This describes the integral boundary directions; it does not identify the moduli space of regular metrics with that Grassmannian. Cone balance still imposes the amplitudes of the boundary spinor descendants in one period lattice.

\section{Classification by endpoint compositions}\label{sec:endpoint-composition}

The polygon records which directions occur. An endpoint composition retains how successive groups of rods fit together when the AF twist approaches an endpoint of its admissible interval. Each component carries its own ordered word, normalized exterior rods and asymptotic circle marking. For an ALF component, this includes the primitive lattice direction of the circle whose length stays bounded at infinity; Appendix~\ref{app:ALF-system} explains the marking and its possible multiplicity. These additional data distinguish assemblies that have the same polygon. We first give exact compatibility rules, then separate their finite catalogue consequence from the analytic opening problem.

\subsection{Marked components and the joining rule}
Write the normalized rods of a component $C$ as $R_0,\ldots,R_m$, where $m$ is its number of genuine turning points. All components are expressed in a common Killing frame and period scale $L$. Put $\phi_i=R_i^\phi$ and normalize the two exterior values to $|\phi_0|=|\phi_m|=1$. Define its left and right boundary parameters by
\begin{equation}\label{eq:endpoint-boundary-parameters}
 \beta_L(C)=-\frac{\phi_1}{\phi_0},\qquad
 \beta_R(C)=-\frac{\phi_{m-1}}{\phi_m}.
\end{equation}
A component is \emph{marked} when these parameters, its orientation, moment rank and primitive period lattice are retained. The name of its underlying metric alone is insufficient.

For its residue spinors $u_C,v_C$, let
\begin{equation}\label{eq:endpoint-twist}
 \gamma_C=u_C^Tv_C,\qquad g_C=\frac{\gamma_C}{L}.
\end{equation}
The ALF residue identities give $\|u_C\|^2=\|v_C\|^2=2M_C$; the mixed moment need not vanish. The difference between the right and left exterior slopes $R^t/R^\phi$ is $2\gamma_C$. Thus $g_C$ records the component's contribution to the total asymptotic twist. It is unrelated to the positive infinity-atom parameter in the Hankel representation of \cite{ADN}.

Suppose components $C,D$ meet across a separating rod. Its adjacent angular components are $-\beta_R(C)\phi$ and $-\beta_L(D)\phi$, where $\phi=\pm1$ is the angular component of that separating rod. Substitution in the signed rod recurrence gives
\begin{equation}\label{eq:endpoint-join}
 d_{\rm join}=\beta_R(C)+\beta_L(D).
\end{equation}
Hence an admissible join must be integral and belong to $\{-1,0,1\}$. The parent self-intersection word is obtained by concatenating the component words with these joining entries. Its number of turning points and its moment rank are the sums of the corresponding component data.

For an ordered assembly $(C_1,\ldots,C_s)$ define the partial twist sums
\[
 G_j=\sum_{\alpha=1}^j g_{C_\alpha}.
\]
AF alignment requires $G_s=0$. At a separating rod the normalized slope is $2LG_j$. In the parent lattice its physical vector is $(-Lq_i,p_i-q_i\beta)$; therefore $q_i=-2\phi_iG_j$ at the endpoint. Differentiating its angular square gives
\[
 \left.\frac{d}{d\beta}(p_i-q_i\beta)^2\right|_{\beta=\beta_*}=4G_j.
\]
The strict angular bound opens into the allowed interval precisely when
\begin{equation}\label{eq:endpoint-prefix}
 \begin{aligned}
 G_j&<0 &&\text{at a left endpoint},\\
 G_j&>0 &&\text{at a right endpoint},
 \end{aligned}
 \qquad 1\le j<s.
\end{equation}
These are necessary first-order angular conditions. They do not impose cone balance inside the opened parent.

\subsection{Exact reconstruction of the AF catalogue}
Table~\ref{tab:endpoint-generators} gives a small component library. The letters $\mathsf N,\mathsf B,\mathsf S$ denote Taub--NUT, Taub--bolt and Schwarzschild components, respectively \cite{GH,Page}. The signs on $\mathsf N$ and $\mathsf B$ specify the sign of $g$, not a self-intersection. The four-nut component $C$ is $\mathcal C_4$; the eight-nut component $D$ is \texttt{ALF8-026} in Appendix~\ref{app:ALF-system}. A bar reverses the component orientation, changing the signs of its word, boundary parameters and $g$, and replacing rank $k$ by $N-k$.

\begin{table}[htbp]\centering\small
\begin{tabular}{@{}lrrlrr@{}}\toprule
Component & $N$ & $k$ & Internal word & $\beta_L=\beta_R$ & $g$\\\midrule
$\mathsf N_-$&1&0&$\varnothing$&$-1$&$-1/2$\\
$\mathsf N_+$&1&1&$\varnothing$&$1$&$1/2$\\
$\mathsf B_-$&2&1&$(1)$&$0$&$-1/2$\\
$\mathsf B_+$&2&1&$(-1)$&$0$&$1/2$\\
$\mathsf S$&2&1&$(0)$&$0$&$0$\\
$C$&4&2&$(-1,1,-1)$&$-1/3$&$3/2$\\
$\overline C$&4&2&$(1,-1,1)$&$1/3$&$-3/2$\\
$D$&8&3&$(-1,1,-1,-1,-1,1,-1)$&$-2/3$&$-3/2$\\
$\overline D$&8&5&$(1,-1,1,1,1,-1,1)$&$2/3$&$3/2$\\\bottomrule
\end{tabular}
\caption{Nine marked components sufficient to generate the retained AF words through twelve turning points. The library generates words; it is not the list of all component types appearing at both ends.}\label{tab:endpoint-generators}
\end{table}

\begin{proposition}[Finite endpoint reconstruction]\label{prop:endpoint-census}
Enumerate ordered assemblies from Table~\ref{tab:endpoint-generators} with total turning-point number $3\le N\le12$. Impose \eqref{eq:endpoint-join}, $G_s=0$, \eqref{eq:endpoint-prefix}, and an endpoint twist in $[-1,0]$. Retain the AF end conditions and admissibility tests of Appendix~\ref{app:numerics}, and quotient by its reversal and orientation conventions. The resulting set is exactly the 356 retained AF words in Table~\ref{tab:general-counts}; none of the three excluded words is generated.
\end{proposition}
\begin{proof}[Finite verification]
The component words and boundary parameters are integral or rational, so every operation is exact. Exhaustive generation gives, for $N=3,\ldots,12$,
\[
 1,\ 1,\ 1,\ 3,\ 3,\ 9,\ 15,\ 40,\ 78,\ 205.
\]
Comparison of canonical words gives equality of sets, not only equality of counts. Independently, concatenating the components in all 712 recorded endpoint assemblies reconstructs their parent words. The joining formula is checked at all 3,607 joins, and all partial twist sums have the required signs. The generator and exact checks are supplied with the catalogue.
\end{proof}

The five classical marked components alone generate 355 of the 356 words. The remaining word, \texttt{12-135}, is
\[
 (-1,1,-1,-1,-1,1,-1,-1,-1,1,-1).
\]
At its left endpoint $\beta_*=-2/3$, the assembly $(D,C)$ has joining entry $-2/3-1/3=-1$ and partial sum $-3/2$. At its right endpoint $\beta_*=-1/3$, the assembly $(C,D)$ has the same joining entry and partial sum $3/2$. Both reconstruct the displayed word. This example shows why nonclassical ALF components enter even at low order.

Of the 356 words, 197 admit classical assemblies at both endpoints, 158 at one endpoint, and one at neither. Generation from one classical endpoint does not make the other endpoint classical: the endpoint atlas also contains five- and six-nut ALF components and marked Kerr components. A fully marked assembly determines the parent word, hence its necessary interval $\mathcal I_W$. Cutting that word where $|p_i-q_i\beta|=1$ at the opposite end predicts the opposite partition. The metric limit is a further assertion, not part of this arithmetic reconstruction.

At fourteen turning points the earlier component library generated 1,141 words. The 18 additional numerically convergent words require fourteen further component types, three with nine nuts and eleven with ten. After these components are extracted, the same rules generate exactly the 1,159 convergent words and none of the 66 unconverged candidates. This extended agreement is a reconstruction using components identified from the data, not an independent existence test or a proof excluding the 66 failures. No fixed finite component library is asserted to suffice at arbitrary order.

\subsection{The analytic opening problem}
The next step is to turn compatible boundary data into a regular interior metric. Even the leading separation law retains information absent from the polygon. Put $m_\alpha=M_{C_\alpha}/L$ and $g_\alpha=g_{C_\alpha}$, and define
\[
 K_{\alpha\delta}=8m_\alpha m_\delta-2g_\alpha g_\delta\ge0.
\]
The inequality follows from $|\gamma_C|\le2M_C$. Assume that a balanced family has separated components of bounded size, and that their centres satisfy $a_\alpha/L=X_\alpha/\sqrt\varepsilon+O(1)$ as $\beta=\beta_*+\sigma\varepsilon$, where $\sigma=1$ at a left endpoint and $\sigma=-1$ at a right endpoint. The leading balance equations are
\begin{equation}\label{eq:endpoint-leading-opening}
 \sum_{\alpha\le j<\delta}
 \frac{K_{\alpha\delta}}{(X_\delta-X_\alpha)^2}
 =-4\sigma G_j,\qquad 1\le j<s.
\end{equation}
Indeed, $T^2=I$ gives $T_{[j]}^2+T_{[j],>j}T_{[j],>j}^T=I$, while the displacement identity in \eqref{eq:T-identities} gives the cross-component entries. Their squared leading norm is $L^2K_{\alpha\delta}/(a_\delta-a_\alpha)^2$. Expanding the squared angular determinant and comparing it with its derivative $4G_j$ yields \eqref{eq:endpoint-leading-opening}.

For gaps $y_j=X_{j+1}-X_j>0$, these equations are the critical-point equations of
\[
 V(y)=\sum_{\alpha<\delta}
 \frac{K_{\alpha\delta}}{y_\alpha+\cdots+y_{\delta-1}}
 -4\sigma\sum_{j=1}^{s-1}G_jy_j.
\]
If every adjacent $K_{j,j+1}>0$ and $\sigma G_j<0$, this function is strictly convex, diverges at each collision face and is coercive at infinity. It therefore has a unique minimum. This determines leading relative separations under the assumed degeneration; it does not construct the full family or control the remainder.

A gluing theorem must also solve the internal component equations, retain the common period normalization, and bound the correction to \eqref{eq:endpoint-leading-opening}. The collision hypotheses of Section~\ref{subsec:collision-matching} concern precisely the comparison needed to continue such a family globally. The rational recursions of \cite{DMT} and the positive moment representation of \cite{ADN} provide explicit coordinates for this problem, but neither a rational patching matrix nor compatibility of its endpoint factors enforces cone balance by itself.

\section*{Acknowledgements}

I thank Lars Andersson for pointing out a gap in an earlier version of the triangular-transport argument at rods with zero angular component.

The Belinski--Zakharov factorization presented here was obtained by the author in 2022. ChatGPT (OpenAI) was subsequently used to explore its consequences, assist with Mathematica and Python code development, and support the writing and overall presentation of the manuscript. Computations for the catalogue were performed using Mathematica locally by the author and Python within ChatGPT's code-execution environment.

\appendix

\section{AF catalogue through twelve turning points}\label{app:numerics}

This appendix specifies the finite algebraic system and records the 356 surviving AF word classes for $3\le N\le12$, modulo axis reversal and orientation. The interval computations use Arb ball arithmetic \cite{Arb}. Every listed class satisfies the proper-block collision criterion and has an interval-certified balance root at the midpoint of its strict twist interval. Numerical continuation supplies mass, angular momentum and normalized lengths; exact rod arithmetic determines the endpoint partitions. The three words excluded by the boundary-gradient argument are treated separately in Appendix~\ref{app:exceptional} and do not occur in the catalogue or its endpoint atlas.

The main distinctions are the 197 full-window words with $\mathcal I_W=(-1,0)$, the 159 restricted-window words, and the different marked walks carried by the same polygon. Every retained word is reproduced by the endpoint composition rules of Proposition~\ref{prop:endpoint-census}. The tables therefore serve both as a list of certified solutions and as data for testing the two classification approaches. The mass curves satisfy the first-law relation to angular momentum; a zero-spin point is stationary, with the observed minima recorded below.

The certificates, continuation samples and endpoint partitions have distinct roles. A root certificate proves the balance equations at one parameter. Full-interval existence uses Theorem~\ref{thm:collision-existence} and Hypotheses~\ref{hyp:comparison}, rather than interpolation between samples. Endpoint partitions specify marked candidate components; convergence to their metrics requires a separate degeneration analysis.

\subsection{Enumeration conventions}
Here $N$ is the number of genuine turning points, $n$ is the even number of dressing slots, and $h=n-N$ is the number of same-sign phantoms after redundant opposite-sign pairs have been removed. The two polynomial degrees are $\{k,N-k\}$, where
\[
 k=\frac{N-h}{2},\qquad n=N+h=2(N-k),\qquad 1\le k\le N/2.
\]
The familiar two-nut Kerr family is a separate base case. The tables start with the three-nut Chen--Teo word, represented by $[4-1-0]$, and include odd as well as even $N$.

Write $A=(0,1)$, $B=(L,\beta)$, with $L>0$. A signed lattice rod is
\[
 \mathsf v_j=p_jA-q_jB=(-q_jL,p_j-q_j\beta).
\]
For a self-intersection word $d=(d_1,\ldots,d_{N-1})$, the integral coefficient vectors are reconstructed from
\begin{equation}\label{eq:app-recurrence}
 (p_0,q_0)=(1,0),\quad(p_1,q_1)=(0,1),\qquad
 (p_{j+1},q_{j+1})=-(p_{j-1},q_{j-1})-d_j(p_j,q_j).
\end{equation}
Enumeration fixes $d_1=-1$ and applies the following tests:
\[
 \begin{gathered}
 d_j\in\{-1,0,1\},\qquad (p_N,q_N)=(\sigma,0),\quad \sigma\in\{-1,1\},\\
 q_j\ne0\quad(1\le j\le N-1),\qquad
 |p_aq_b-q_ap_b|<b-a\quad(b-a\ge2),\\
 \mathcal I_W=\{\beta\in(-1,0):|p_j-q_j\beta|<1\text{ for every finite rod}\}\ne\varnothing.
 \end{gathered}
\]
All interval endpoints and integer tests are computed exactly. The projective winding $\mathfrak w$ of the signed rod sequence determines $h=|N-2\mathfrak w|$ and $k=\min(\mathfrak w,N-\mathfrak w)$. The resulting sectors satisfy $h\equiv N\pmod2$ and $3h<N$ for $N>3$; the three-nut word has $h=1$.

The final enumeration step is a quotient, not a symmetry constraint on the unknown metric. For $d_{N-1}=-1$, normalized reversal sends $d$ to its reverse and $\beta$ to $-1-\beta$. For $d_{N-1}=1$, it sends $d$ to minus its reverse and leaves $\beta$ unchanged. We choose the lexicographically smaller word and retain the parameter maps in the data. An oriented classification would retain the orientation labels. No equality between different node positions, gaps or charges is imposed in the solver.

\begin{table}[htbp]\centering
\begin{tabular}{@{}rrlr@{}}\toprule
$N$ & Classes & Sectors ($h$: number) & Certified roots\\\midrule
3&1&$1:1$&1\\
4&1&$0:1$&1\\
5&1&$1:1$&1\\
6&3&$0:3$&3\\
7&3&$1:3$&3\\
8&9&$0:8$, $2:1$&9\\
9&15&$1:15$&15\\
10&40&$0:35$, $2:5$&40\\
11&78&$1:77$, $3:1$&78\\
12&205&$0:161$, $2:44$&205\\
\midrule
Total &356&&356\\\bottomrule
\end{tabular}
\caption{Surviving AF word classes and certified midpoint balance roots. The certificates comprise 208 genuine even-sector roots and 148 roots in reduced phantom sectors.}\label{tab:general-counts}
\end{table}

\subsection{The full system solved}
We use the outer-node normalization $z_1=-1$, $z_N=1$; the period scale $L$ is recovered afterwards. The unknowns are the $N$ positive effective charge weights $u_i$ and the $N-2$ interior nodes. Set
\[
 D(x)=\prod_{i=1}^N(x-z_i),\qquad
 \omega_i=\frac{u_i}{|D'(z_i)|}>0,\qquad
 \mu_a^{[j]}=\sum_{i=1}^N\epsilon_i^{[j]}\omega_i z_i^a,
 \quad \epsilon_i^{[j]}=\begin{cases}-1&i\le j,\\1&i>j.\end{cases}
\]
For $j=0,\ldots,N$, define
\[
 H_j=(\mu_{a+b}^{[j]})_{a,b=0}^{k-1},\quad
 c_j=(\mu_k^{[j]},\ldots,\mu_{2k-1}^{[j]})^T,\quad
 e=(0,\ldots,0,1)^T
\]
and the two-component homogeneous rod
\begin{equation}\label{eq:app-moment-rods}
 \mathsf t_j=\left(
 \det\begin{pmatrix}H_j&c_j+e\\(c_j+e)^T&\mu_{2k}^{[j]}+\sum_i z_i\end{pmatrix},
 \det H_j\right).
\end{equation}
These are the reduced moment rods. They incorporate the cone normalization, so equality of the complete vectors, rather than only their directions, is essential.

If $\mathfrak w\ne k$, the calculation uses the orientation-reversed word $d^{\rm s}=-d$ and parameter $\beta_{\rm s}=-\beta$; otherwise $d^{\rm s}=d$ and $\beta_{\rm s}=\beta$. Write $\varepsilon=\pm1$ for this recorded orientation change and reconstruct $(p_j^{\rm s},q_j^{\rm s})$ from $d^{\rm s}$. The $2N-2$ scalar equations are
\begin{equation}\label{eq:app-balance}
 \begin{aligned}
 \mathsf t_j&=p_j^{\rm s}\mathsf t_0+q_j^{\rm s}\mathsf t_1,
                         &&j=2,\ldots,N-1,\\
 t_N^t&=\sigma_{\rm s}t_0^t,
                         &&\sigma_{\rm s}=(-1)^k,\\
 t_1^\phi+\beta_{\rm s}t_0^\phi&=0.
 \end{aligned}
\end{equation}
The other exterior component is automatic: $t_N^\phi=(-1)^kt_0^\phi$. Each equation is divided by the nonzero unflipped angular factor when certified. After the common exterior shear and normalization,
\begin{equation}\label{eq:app-units}
 L=\frac{\det(\mathsf t_0,\mathsf t_1)}{(t_0^\phi)^2}>0,
 \qquad x_j=\frac{z_{j+1}-z_j}{L},\qquad
 \overline M=\frac ML,\quad\overline J=\frac J{L^2}.
\end{equation}
The resulting rods are exactly $(-q_j^{\rm s}L,p_j^{\rm s}-q_j^{\rm s}\beta_{\rm s})$ at balance. The reported mass and lengths are orientation-independent; the reported spin is multiplied by $\varepsilon$ on returning to the canonical marking.

The numerical search uses $\log u_i$ and logarithmic ratios of the first $N-2$ gaps to the last gap. Positivity and ordering are therefore built into the chart. All $2N-2$ variables are independent. Positive occupation sums evaluate the rods and their analytic Jacobian without expanding a high-degree elimination polynomial. The initial search is at the exact midpoint of $\mathcal I_W$, with additional starts when needed. Original homogeneous BZ determinants provide a second evaluation after convergence.

\paragraph{Reinserting phantoms.}
A root of the reduced system has an explicit even-slot representative. Choose $h$ auxiliary nodes $a_s=3+s$, $s=0,\ldots,h-1$, to the right of the genuine nodes and set
\[
 y_i=\frac{u_i}{\prod_s(a_s-z_i)},\qquad
 P_i=1-(-1)^{i-1}y_i,\quad Q_i=1+(-1)^{i-1}y_i.
\]
The genuine affine BZ charge is $C_i=P_i/Q_i$ when $Q_i\ne0$. Each auxiliary slot has $(P,Q)=(-1,1)$, hence $C=-1$. The homogeneous chart includes affine charge infinity. Removing the common phantom factors gives \eqref{eq:app-moment-rods}. Across each phantom the homogeneous rod changes sign, while its vanishing circle is unchanged. The independent checks use this sign and compare all genuine full rods and the recovered scale $L$ with the original $n$-slot determinants. The auxiliary nodes are representation choices, not additional geometric moduli.

\subsection{What a balance certificate proves}
For each accepted midpoint, Newton refinement at 90 decimal digits produces a centre in the variables
\[
 (u_1,\ldots,u_N,z_2,\ldots,z_{N-1}).
\]
Outward-rounded ball arithmetic then treats a box of radius $10^{-32}$ in every coordinate. For the rational equation map $F$ in \eqref{eq:app-balance}, a fixed decimal matrix $B$ approximates the inverse Jacobian at the centre $x_0$. The stored certificate checks
\[
 \sup_{x\in X}\|I-B\,DF(x)\|_\infty\le q<1,\qquad
 \|BF(x_0)\|_\infty+q\,10^{-32}<10^{-32}.
\]
It also verifies positive weights, strictly ordered nodes and every denominator used in the evaluation. The contraction map $x\mapsto x-BF(x)$ consequently has a unique zero in that box, and its Jacobian is invertible. This is local uniqueness and gives a local parameter branch by the implicit function theorem.

All 356 midpoint roots passed. The largest contraction bound is below $4.49\times10^{-17}$. At their high-precision centres, the largest original-BZ full-rod discrepancy is below $3.40\times10^{-74}$; the latter is an independent consistency check, not the existence argument. The files retain the box centres, approximate inverses and original homogeneous charges, so the checks can be rerun.

The certificates establish exact balanced data in the positive reduced representation. Bulk positivity follows from the matrix reconstruction, and the primitive rods and unit adjacent lattice determinants are exact properties of the prescribed word. Corollary~\ref{cor:smooth-completion} therefore gives smooth complete metrics at the certified parameters; this conclusion does not require the homogeneous collision hypotheses. The finite contraction calculation supplies the exact balanced data, while the analytic reconstruction supplies completion. Independent whole-axis interval audits are available for the five types \texttt{6-002}, \texttt{8-002}, \texttt{8-006}, \texttt{8-007} and \texttt{8-009}. They are additional checks, not an audit of every row. Global uniqueness in the positive reduced sector follows from Theorem~\ref{thm:convexity}.

\subsection{Explicit low-order solutions}\label{app:explicit-low-order}
The finite moment system gives exact solutions without requiring a radical parametrization: a polynomial system, a specified real branch and a reconstruction map are sufficient. The following examples make this description concrete. Here $N$ always counts genuine turning points. Write
\[
 \kappa=\frac{z_N-z_1}{2}
\]
for half the outer node span. This coordinate scale differs from the asymptotic period scale $L$. In particular $M/\kappa$ and $M/L$ must not be compared without converting the normalization. The symmetric ansatz used in the four- and six-nut examples is a simplification of those examples; it is not imposed on the catalogue.

\paragraph{Four turning points: the Teo family.}
For $d=(-1,0,1)$, corresponding to \texttt{4-001}, take
\begin{equation}\label{eq:low4-ansatz}
 z=\kappa(-1,-r,r,1),\qquad
 C=(a,b,-b^{-1},-a^{-1}),\qquad 0<r<1.
\end{equation}
The reciprocal charges satisfy AF alignment. The chamber containing the regular branch has $0<a<1<b$. Introduce the charge polynomial
\[
 \varphi(a,b)=(a-1)^2(1+b^2)+4ab,
 \qquad U=\varphi(a,b),\quad V=\varphi(b,a),
\]
and its two derivatives
\[
 N_c=\tfrac12\partial_a\varphi=(a-1)(1+b^2)+2b,
 \qquad D_c=\tfrac12\partial_b\varphi=b(a-1)^2+2a.
\]
The normalized rod equations reduce to
\begin{equation}\label{eq:low4-balance}
 (1+r^2)U-2rV=0,\qquad N_c-rD_c=0.
\end{equation}
Since $D_c>0$ in this chamber, elimination of $r$ gives
\begin{equation}\label{eq:low4-charge-curve}
 U(D_c^2+N_c^2)-2N_cD_cV=0,\qquad r=N_c/D_c.
\end{equation}
This is a compact charge-curve description of the four-nut AF solution of \cite{Teo}. The inequalities $0<r<1$, positivity and the normalized rod relations select the physical component of the algebraic curve.

A useful alternative chart is
\[
 a=\frac{1-u}{1+u},\qquad b=\frac{1+v}{1-v},\qquad
 r=\frac{1-t}{1+t},\qquad 0<u<v<1,\quad0<t<1.
\]
Equation~\eqref{eq:low4-balance} becomes the polynomial pair
\begin{equation}\label{eq:low4-cayley}
 \begin{aligned}
 u^2-v^2+t^2+3u^2v^2t^2&=0,\\
 u^2\{1+v(1-t)+3v^2t\}+2uv^2(1+t)+(v-1)(v+t)&=0.
 \end{aligned}
\end{equation}
For example, setting $v=m$ and $u=m\gamma$ determines
\[
 t=\frac{m\sqrt{1-\gamma^2}}{\sqrt{1+3m^4\gamma^2}},\qquad 0<\gamma<1,
\]
and leaves one equation $A+tB=0$, where
\[
 \begin{aligned}
 A&=m^2\gamma^2(1+m)+2m^3\gamma+m(m-1),\\
 B&=m^3\gamma^2(3m-1)+2m^3\gamma+m-1.
 \end{aligned}
\]
The eliminated equation $A^2(1+3m^4\gamma^2)-m^2(1-\gamma^2)B^2=0$ must retain the unsquared sign condition $A+tB=0$.

A numerical zero-spin representative is
\[
 (a,b,r)\simeq(0.5228834958296844,\ 3.4674923457823892,\ 0.3930322274582969).
\]
The residue formulas give $M/\kappa\simeq0.6415956930225703$, $L/\kappa\simeq1.6003207199276681$, and hence $M/L\simeq0.4009169443557347$ at $\beta\simeq-0.4448362860057288$. These values agree with Table~\ref{tab:zero-general}. The numerical location of this spin zero is separate from the certified midpoint existence result.

\paragraph{Five turning points: a one-phantom representative.}
The sole five-nut word in the census is
\[
 d=(-1,0,0,-1),\qquad
 (v_0,\ldots,v_5)=((1,0),(0,1),(-1,1),(0,-1),(1,-1),(1,0)).
\]
Its reduced rank is $k=2$, so it has a six-slot BZ representative with one phantom. At $\beta=-1/2$, set $\kappa=1$ and use the positive occupation variables of \eqref{eq:app-moment-rods}:
\[
 z=(-1,-r,0,r,1),\qquad (u_1,\ldots,u_5)=(a_5,b_5,c_5,b_5,a_5).
\]
Here $a_5,b_5,c_5$ are positive occupations, not affine BZ charges. Reflection of the signed moments gives $t_{5-j}^t=-t_j^t$ and $t_{5-j}^\phi=t_j^\phi$. Consequently the complete balance equations reduce exactly to four scalar equations,
\begin{equation}\label{eq:low5-system}
 t_0^t=0,\qquad
 t_1^\phi=\tfrac12t_0^\phi,\qquad
 t_2^t=t_1^t,\qquad
 t_2^\phi=-\tfrac12t_0^\phi.
\end{equation}
Each $\mathsf t_j$ is the explicit $2\times2$ and bordered $3\times3$ determinant in \eqref{eq:app-moment-rods}. The certificate for \texttt{5-001} selects the positive solution near
\[
 \begin{aligned}
 r&=0.5221304683685786,&a_5&=0.4994219409569822,\\
 b_5&=0.7106895457957642,&c_5&=0.9800704111573928.
 \end{aligned}
\]
Its dimensionless invariants are $M/L\simeq0.4533846765481244$ and $J=0$. In these symmetric moment coordinates, $a$ is odd and $c$ is even in \eqref{eq:reduced-dictionary}, and $u_i v_i=a_i^2-c_i^2$ is even.

To recover ordinary BZ data, place the phantom at $z_6=3$ and put $y_i=u_i/(3-z_i)$. The genuine charges are
\[
 C_i=\frac{1-(-1)^{i-1}y_i}{1+(-1)^{i-1}y_i},\quad1\le i\le5,
 \qquad C_6=-1.
\]
Thus the reduced determinant solution gives explicit six-slot data as well as the five genuine rods. The auxiliary sixth position changes the representation, not the geometry.

\paragraph{Six turning points beyond the Li--Sun word.}
The six-nut census contains three words,
\[
 (-1,0,0,0,1),\qquad (-1,1,0,-1,1),\qquad(-1,1,1,1,-1).
\]
We give a charge chart for the second, \texttt{6-002}. With nodes in ascending order, take
\begin{equation}\label{eq:low6-ansatz}
 z=\kappa(-1,-r,-s,s,r,1),\qquad
 C=(a,b,c,-c^{-1},-b^{-1},-a^{-1}),\qquad0<s<r<1.
\end{equation}
The relevant positive chamber has $0<a<1$, $b>1$, $-1<c<1$, with $c\ne0$ in the displayed affine chart. Put
\[
 R=\frac{1+r}{1-r},\qquad S=\frac{1+s}{1-s},\qquad 1<S<R,
\]
and retain $U=\varphi(a,b)$, $V=\varphi(b,a)$ from the four-nut example. Define two charge vectors by
\begin{align*}
 \mathbf f_1={}&\begin{pmatrix}
 -\tfrac12(1+c^2)U\\
 -\tfrac12(1+c^2)V\\
 \tfrac12(1-c^2)(a^2-1)(b^2+1)+c\{V-2a(b+1)^2\}\\
 \tfrac12(1-c^2)(a-1)^2(b^2-1)-2c(a^2b+ab^2-2ab+a+b)
 \end{pmatrix},\\
 \mathbf f_2={}&\begin{pmatrix}
 a\{(1+c^2)(b+1)^2+4bc\}\\
 (1+c^2)(a^2+b)(b+1)+2c(a^2+b^2)\\
 2a^2b+c(1+a^2)(b+1)^2+2bc^2\\
 a\{(1+c^2)(b^2+1)+2c(b+1)^2\}
 \end{pmatrix}.
\end{align*}
The spectral vector $\mathbf w=(W_0,W_1,W_2,W_3)^T$ is
\[
 \begin{aligned}
 W_0&=(R-S)^2(1+R^2S^2)+(RS-1)^2(R^2+S^2),\\
 W_1&=-2S(R^2-1)(R-S)(RS-1),\\
 W_2&=2R(R-S)(S^2-1)(RS-1),\\
 W_3&=-2RS(R^2-1)(S^2-1).
 \end{aligned}
\]
For the derivative equations introduce
\[
 \begin{aligned}
 \Omega_a&=4RS\bigl((1+R^2)(1+S^2)-4RS,\ 0,\ 0,\ -(R^2-1)(S^2-1)\bigr)^T,\\
 \Omega_b&=-\frac{4R(R-1)(R-S)(RS-1)}{R+1}
               \bigl(1+S^2,\ 0,\ -(S^2-1),\ 0\bigr)^T,\\
 \Omega_c&=\frac{4S(S-1)(R-S)(RS-1)}{S+1}
               \bigl(1+R^2,\ -(R^2-1),\ 0,\ 0\bigr)^T.
 \end{aligned}
\]
Substitution into the normalized rod equations gives the four conditions
\begin{equation}\label{eq:low6-system}
 \mathbf w^T\mathbf f_i=0,\qquad
 \sum_{q\in\{a,b,c\}}\Omega_q^T\partial_q\mathbf f_i=0,
 \qquad i=1,2.
\end{equation}
The derivatives act only on the charge variables. The denominators $R+1$ and $S+1$ are strictly positive in the stated node chamber, so the system can be made polynomial without introducing a vanishing denominator. These equations describe a one-dimensional algebraic solution set before the positivity and branch conditions are imposed. They exhibit a concrete relation between charge derivatives and the rod equations without requiring elimination to one high-degree polynomial.

The signed primitive rods of the intended branch are
\[
 (1,0),\ (0,1),\ (-1,1),\ (1,-2),\ (1,-1),\ (0,1),\ (-1,0).
\]
A numerical zero-spin representative is
\[
 \begin{aligned}
 a&\simeq0.5220893460504066,& b&\simeq13.0668692937716688,\\
 c&\simeq-0.3477147817372575,& r&\simeq0.7010419372452687,\\
 s&\simeq0.4049395967930934.
 \end{aligned}
\]
It has $M/\kappa\simeq0.6976106882824007$, $L/\kappa\simeq1.2838780930652332$, and therefore $M/L\simeq0.5433620933720189$ at $\beta\simeq-0.5033316240992153$. Independent positive-moment reconstruction checks the full rods at this point and at two nonzero-spin points with $r=0.65$ and $r=0.75$; the residuals are below $10^{-60}$. These checks validate the reconstruction conventions. Rigorous midpoint existence and the local family follow from the certificate for \texttt{6-002}, while global continuation remains subject to the hypotheses in Section~\ref{subsec:collision-matching}.

The exact endpoint partitions for this word are, from left to right,
\[
 \begin{array}{c|c}
 \beta=-1&\varnothing,\ (1),\ (-1),\ \varnothing\\
 \beta=0&(-1),\ \varnothing,\ \varnothing,\ (1).
 \end{array}
\]
Here $\varnothing$ is a one-nut component. Thus the two endpoints have different ordered assemblies of nuts and bolts despite using the same classical component types. The marking and joining rules of Section~\ref{sec:endpoint-composition} retain this ordering.

\subsection{Asymptotic invariants along the branches}\label{app:mass}
For every certified word, continuation was corrected at
\[
 \beta=\beta_-+s(\beta_+-\beta_-),\qquad
 s\in\{0.01,0.05,0.20,0.35,0.50,0.65,0.80,0.95,0.99\}.
\]
All 3204 standard samples passed the full-rod test. The maximum residual in the integral rod coordinates was below $4.18\times10^{-11}$, and the smallest sampled reciprocal Jacobian condition number was greater than $7.10\times10^{-5}$ in the logarithmic search coordinates. These are finite sample observations. Neither a uniform bound between samples nor continuation to an endpoint is certified by this sweep.

The invariants are evaluated directly from the positive moment data. Put $H=H_0$, let $V_{ia}=z_i^a$ for $a=0,\ldots,k-1$, and set
\[
 W=\diag(\sqrt{\omega_i})V,\quad
 b=V^T(\omega_i z_i^k)_{i=1}^N,\quad
 a=WH^{-1}e,\quad
 c=(\sqrt{\omega_i}z_i^k)_{i=1}^N-WH^{-1}b.
\]
With residue columns $u^{\rm res}=a+c$, $v^{\rm res}=a-c$ and $Z=\diag(z_i)$, the asymptotic moments are
\begin{equation}\label{eq:app-mass-spin}
 M=\frac12\|u^{\rm res}\|^2,\qquad
 J=-\frac12(u^{\rm res})^TZv^{\rm res}.
\end{equation}
The superscript distinguishes residue columns from the positive scalar weights $u_i$. Here $a\perp c$ and AF alignment is $\|a\|^2=\|c\|^2$.

For a stable second evaluation, define for $r=k-1,k,k+1$
\[
 A_r=\sum_{|I|=r}\prod_{i\in I}\omega_i
                \prod_{i<j,\ i,j\in I}(z_j-z_i)^2,
 \qquad
 S_r=\frac1{A_r}\sum_{|I|=r}\left(\sum_{i\in I}z_i\right)
          \prod_{i\in I}\omega_i\prod_{i<j,\ i,j\in I}(z_j-z_i)^2.
\]
For the empty subset take $A_0=1$, $S_0=0$. Then
\[
 M=\frac{A_{k-1}+A_{k+1}}{2A_k},\qquad
 J=\frac{A_{k+1}(S_{k+1}-S_k)-A_{k-1}(S_k-S_{k-1})}{2A_k}.
\]
These identities follow from the Gram determinants of consecutive monic orthogonal polynomials. They use positive sums, remain valid when a signed prefix determinant vanishes, and agree with the independent high-precision matrix evaluation at every certified centre: the largest double-precision differences in $\overline M$ and $\overline J$ were $1.20\times10^{-14}$ and $4.70\times10^{-14}$.

\begin{figure}[htbp]\centering
\includegraphics[width=\textwidth]{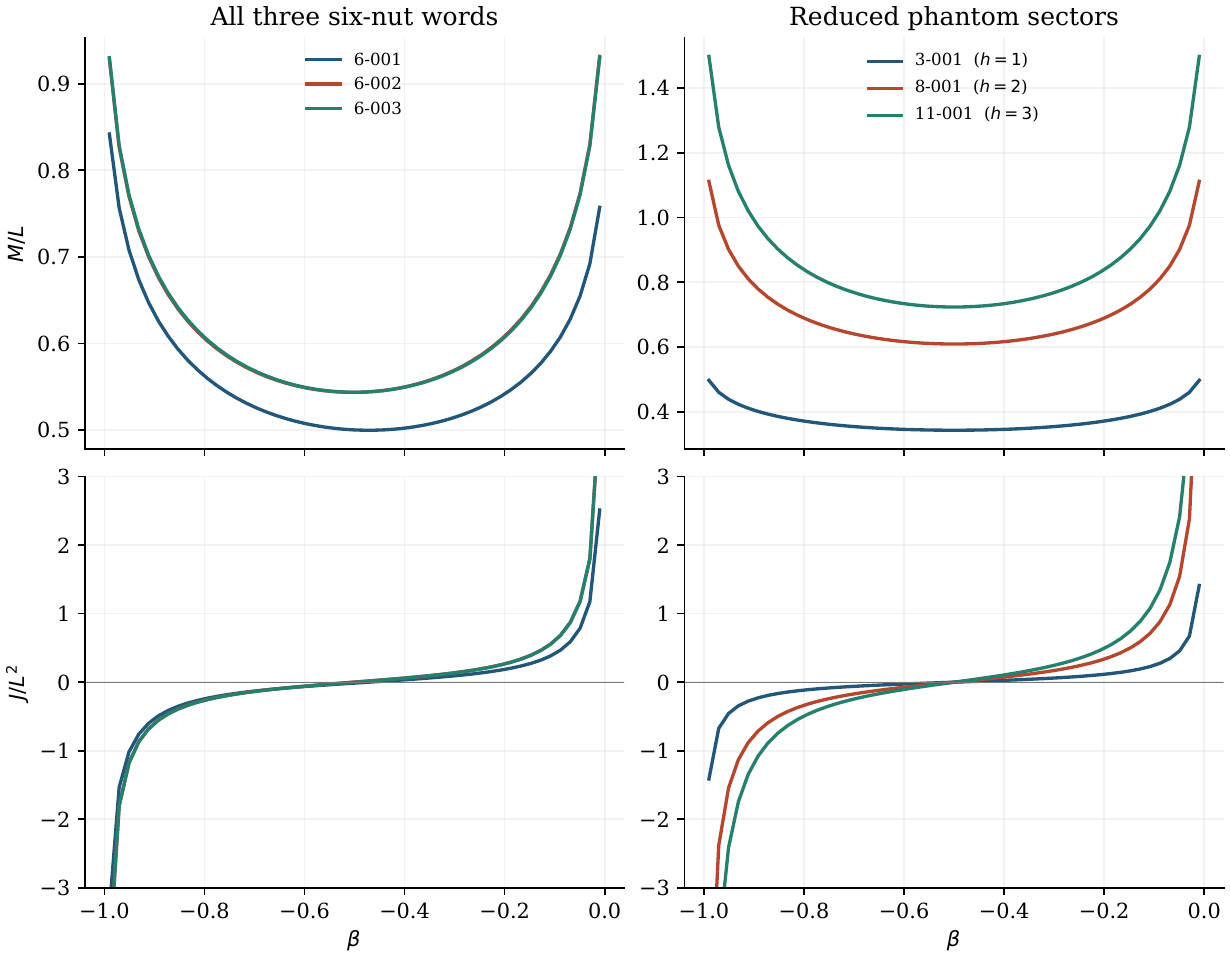}
\caption{Balanced AF invariant curves. Left: all three six-nut words. Right: the three-nut one-phantom, eight-nut two-phantom and eleven-nut three-phantom words. Each line joins 51 corrected samples; endpoint values are not extrapolated. The spin panels are clipped to $[-3,3]$ to show the interior.}\label{fig:general-mass-spin}
\end{figure}
\FloatBarrier

As checks of scale and marking, the Chen--Teo midpoint gives $\overline M=6-4\sqrt2$, with the normalization in \eqref{eq:app-units}. Central differences on five representatives, spanning the genuine and phantom sectors, agree with $d\overline M/d\beta=2\overline J$ to within $3.0\times10^{-10}$ at the tested parameters. This independently checks the implementation of the first law \eqref{eq:full-first-law}.

Each of the 356 sampled curves has negative spin at its leftmost standard sample and positive spin at its rightmost sample. A corrected one-dimensional solve locates one zero-spin representative for each word. Table~\ref{tab:zero-general} gives all words through $N=8$ and the first three-phantom example. The full zero-spin data are included electronically. Uniqueness of a spin zero over the entire branch is not asserted. A vanishing total angular momentum does not imply that the twist field vanishes, and hence does not imply staticity.

\begingroup\small\renewcommand{\arraystretch}{1.08}
\setlength{\tabcolsep}{5pt}
\begin{longtable}{@{}lcrr@{}}
\caption{Numerical representatives with vanishing total angular momentum.}\label{tab:zero-general}\\
\toprule
ID & Self-intersection word & $\beta$ at $J=0$ & $M/L$\\\midrule
\endfirsthead
\toprule
ID & Self-intersection word & $\beta$ at $J=0$ & $M/L$\\\midrule
\endhead
\midrule\multicolumn{4}{r}{Continued on next page}\\
\endfoot\bottomrule\endlastfoot
\texttt{3-001}&$(-1,-1)$&-0.500000000&0.343145751\\
\texttt{4-001}&$(-1,0,1)$&-0.444836286&0.400916944\\
\texttt{5-001}&$(-1,0,0,-1)$&-0.500000000&0.453384677\\
\texttt{6-001}&$(-1,0,0,0,1)$&-0.471359683&0.499576553\\
\texttt{6-002}&$(-1,1,0,-1,1)$&-0.503331624&0.543362093\\
\texttt{6-003}&$(-1,1,1,1,-1)$&-0.500000000&0.543451044\\
\texttt{7-001}&$(-1,0,-1,-1,-1,1)$&-0.484009754&0.577018376\\
\texttt{7-002}&$(-1,0,-1,0,1,-1)$&-0.481056693&0.576771901\\
\texttt{7-003}&$(-1,0,0,0,0,-1)$&-0.500000000&0.543163587\\
\texttt{8-001}&$(-1,0,-1,-1,-1,0,-1)$&-0.500000000&0.609363826\\
\texttt{8-002}&$(-1,0,-1,0,1,0,1)$&-0.469907537&0.608293839\\
\texttt{8-003}&$(-1,0,0,0,0,0,1)$&-0.481632417&0.583490954\\
\texttt{8-004}&$(-1,0,0,1,0,-1,1)$&-0.504310347&0.612529009\\
\texttt{8-005}&$(-1,0,0,1,1,1,-1)$&-0.501649632&0.612732215\\
\texttt{8-006}&$(-1,1,-1,0,1,-1,1)$&-0.658066183&0.759768242\\
\texttt{8-007}&$(-1,1,0,0,0,-1,1)$&-0.530521354&0.628663322\\
\texttt{8-008}&$(-1,1,0,0,1,1,-1)$&-0.527773830&0.629241621\\
\texttt{8-009}&$(-1,1,1,0,-1,-1,1)$&-0.474950467&0.629782822\\
\texttt{11-001}&$(-1,0,-1,-1,-1,-1,-1,-1,0,-1)$&-0.500000000&0.723634123\\
\end{longtable}\endgroup

Two examples illustrate the positive phantom sectors. At $\beta=-1/2$, the word
\[
 \texttt{8-001}:\quad d=(-1,0,-1,-1,-1,0,-1),\qquad [n-h-0]=[10-2-0]
\]
has $\overline M\simeq0.609363826117$ and
\[
 \begin{aligned}
 (x_1,\ldots,x_7)\simeq{}&(0.247301498,0.277587993,0.253411943,\\
 &0.470134767,0.253411943,0.277587993,0.247301498).
 \end{aligned}
\]
The word
\[
 \texttt{11-001}:\quad d=(-1,0,-1,-1,-1,-1,-1,-1,0,-1),\qquad [14-3-0]
\]
has $\overline M\simeq0.723634122631$. Its ten normalized gaps are approximately
\[
 \begin{gathered}
 (0.211393072,0.237170999,0.218740736,0.336329201,0.266762863,\\
  0.266762863,0.336329201,0.218740736,0.237170999,0.211393072).
 \end{gathered}
\]
At these parameters the solution has equal paired lengths, although the system treats all gaps as independent variables.

\subsection{Endpoint components}\label{app:endpoints}
The finite rod spheres have the tridiagonal intersection matrix
\[
 Q_W=\begin{pmatrix}
 d_1&1&&0\\1&d_2&\ddots&\\&\ddots&\ddots&1\\0&&1&d_{N-1}
 \end{pmatrix}.
\]
Exact rational congruence reduction was performed for every candidate. Each matrix has a one-dimensional radical and signature of absolute value $h$. In the moment orientation its inertia is
\[
 (b_+,b_-,b_0)=(k-1,N-k-1,1).
\]
For the standard simply connected toric AF completion, the rod topology gives $\chi=N$ and $b_2=N-1$. These are topological invariants of a putative completion, not separate existence tests. The complete signed rods, intersection inertia and canonical orientation are stored for each word. The higher phantom sectors therefore add topology with nonzero signature as well as new rod orderings.

At either rational wall $\beta_*$ of $\mathcal I_W$, take all indices satisfying
\[
 |p_j-q_j\beta_*|=1,
\]
including $j=0,N$. If $a<b$ are consecutive such indices, the associated candidate component has
\begin{equation}\label{eq:app-end-block}
 \text{nuts}=b-a,\quad
 d^{\rm block}=(d_{a+1},\ldots,d_{b-1}),\quad
 p^{\rm block}=|p_aq_b-q_ap_b|.
\end{equation}
This construction uses only the ordered marked rods. To obtain a pointed metric limit, one must additionally establish the divergence of separating lengths, finite nondegenerate internal geometry, and the matching of periods and local twist. The nine-point sweep above does not establish those limiting assertions.

The 33 rod types occurring for surviving words are listed in Table~\ref{tab:block-types}. A number in the endpoint atlas refers to that table; repeated numbers denote successive blocks. Type 1 is the one-nut Taub--NUT type, type 3 the two-nut Taub--bolt type, and type 2 the zero-square two-nut AF type. The latter is Kerr in general and Schwarzschild only when its inherited local twist is zero. Types 4 and 5 have the three- and four-nut AF rod data. Type 6 is the component $\mathcal C_4$; types 7 and 8 have the rod data of $\mathcal C_5^a,\mathcal C_5^b$, and types 9--13 those of $\mathcal C_6^a,\ldots,\mathcal C_6^e$. Reusing a rod type does not assert uniqueness of its metric or determine its inherited scale and marking. A nonzero boundary index specifies the order of the boundary lattice quotient, not the full oriented lens-space type.

\begingroup\small\renewcommand{\arraystretch}{1.08}
\setlength{\tabcolsep}{5pt}
\begin{longtable}{@{}rclr@{}}
\caption{Rod types occurring in the endpoint partitions, modulo reversal and orientation.}\label{tab:block-types}\\
\toprule
Type & Nuts & Internal self-intersections & Boundary index\\\midrule
\endfirsthead
\toprule
Type & Nuts & Internal self-intersections & Boundary index\\\midrule
\endhead
\midrule\multicolumn{4}{r}{Continued on next page}\\
\endfoot\bottomrule\endlastfoot
1&1&$\varnothing$&1\\
2&2&$(0)$&0\\
3&2&$(-1)$&1\\
4&3&$(-1,-1)$&0\\
5&4&$(-1,0,1)$&0\\
6&4&$(-1,1,-1)$&3\\
7&5&$(-1,0,-1,1)$&3\\
8&5&$(-1,1,1,-1)$&3\\
9&6&$(-1,0,-1,-1,1)$&3\\
10&6&$(-1,0,-1,0,-1)$&3\\
11&6&$(-1,0,0,1,-1)$&3\\
12&6&$(-1,1,-1,-1,1)$&4\\
13&6&$(-1,1,0,1,-1)$&4\\
15&7&$(-1,0,-1,-1,0,-1)$&3\\
16&7&$(-1,0,-1,0,0,1)$&3\\
17&7&$(-1,0,0,0,-1,1)$&3\\
18&7&$(-1,0,0,1,1,-1)$&3\\
19&7&$(-1,1,0,-1,-1,1)$&3\\
20&7&$(-1,1,0,0,1,-1)$&3\\
21&7&$(-1,1,1,1,1,-1)$&3\\
22&8&$(-1,0,-1,-1,-1,-1,1)$&3\\
23&8&$(-1,0,-1,-1,0,0,1)$&3\\
24&8&$(-1,0,-1,-1,0,1,-1)$&3\\
25&8&$(-1,0,-1,0,0,-1,1)$&3\\
26&8&$(-1,0,-1,0,0,0,-1)$&3\\
27&8&$(-1,0,-1,0,1,1,-1)$&3\\
28&8&$(-1,0,0,0,-1,-1,1)$&3\\
29&8&$(-1,0,0,0,0,1,-1)$&3\\
30&8&$(-1,0,0,1,0,0,-1)$&3\\
31&8&$(-1,1,-1,-1,-1,1,-1)$&3\\
32&8&$(-1,1,0,-1,-1,-1,1)$&3\\
33&8&$(-1,1,0,-1,0,1,-1)$&3\\
34&8&$(-1,1,1,1,1,1,-1)$&3\\
\end{longtable}\endgroup

\paragraph{Endpoint partitions of the zero-string family.}
For $W_r=(-1,1,0^{2r+1},-1,1)$, the recurrence gives the following ordered internal words at the two walls:
\[
\begin{aligned}
 \beta=-1:&\quad \varnothing,\ (1),\
 \underbrace{(0),\ldots,(0)}_{r\text{ blocks}},\ (-1),\ \varnothing,\\
 \beta=0:&\quad (-1),\
 \underbrace{\varnothing,\ldots,\varnothing}_{2r+2\text{ blocks}},\ (1).
\end{aligned}
\]
An empty word denotes a one-nut component. Thus adding two zeros introduces one zero-square two-nut block at the left wall and two one-nut blocks at the right wall. These are exact marked rod partitions. The one-nut and $\pm1$ two-nut entries have Taub--NUT and Taub--bolt rod types; a zero-square block has two-nut AF rod data. Its inherited twist is needed to distinguish Kerr from Schwarzschild. The insertion theorem does not establish endpoint metric convergence, noncollapse of these components, or rates of separation.

\paragraph{Metric evidence at selected endpoints.}
High-precision tails and local comparisons are available for \texttt{4-001}, \texttt{6-001}, \texttt{6-002}, \texttt{8-002}, \texttt{8-003}, \texttt{8-006}, \texttt{8-007} and \texttt{8-009}. For the latter four eight-nut branches, the tails reach within $1.08\times10^{-9}$ of the stated walls, with rod residuals below $7.3\times10^{-39}$. These comparisons concern the specified branches; equality of endpoint rod partitions alone does not identify metric limits.

For example, the right wall of \texttt{8-006} is $-1/3$ and its partition is $(6,6)$ in the component numbering, namely two four-nut $\mathcal C_4$ blocks. The independent exact component is given by
\begin{equation}\label{eq:C4-exact}
 \begin{gathered}
 r=\frac{7-2\sqrt{10}}3,\qquad
 u=\frac{r}{3+4\sqrt{3r}},\quad v=3u,\qquad
 a=\frac{1-u}{1+u},\quad c=\frac{1-v}{1+v},\\
 s_0=\frac{12ru}{(1-r^2)(1+u)^2},\qquad
 z=s_0(-1,-r,r,1),\quad C=(a,c^{-1},c,a^{-1}).
 \end{gathered}
\end{equation}
Its three finite lengths in the retained local normalization are
\[
 s_0(1-r),\quad2s_0r,\quad s_0(1-r)
 \simeq 0.073603373799,\quad0.042773778540,\quad0.073603373799.
\]
The parent tail approaches these values while its central separating gap diverges. Independently computed balanced points for component types 9, 11, 12 and 13 have original-rod residuals between $1.43\times10^{-84}$ and $1.86\times10^{-83}$; these are numerical constructions without interval-certified ALF completion. Appendix~\ref{app:ALF-system} records the later balance certificates for seven- and eight-nut components. Matching these separate solutions to the parent endpoint limits still requires metric convergence.

\paragraph{Effective length enclosures.}
The finite supporting-plane calculation for \texttt{6-002} proves a balanced solution for every $|\beta+0.5033316240992|\le9\times10^{-7}$, with $\min_jx_j>0.2002$. For \texttt{8-006}, it gives $\min_jx_j>0.03690$ for every $|\beta+0.4|\le9\times10^{-11}$. These results control the full length space, including cluster collapse, on their stated intervals. These quantitative bounds complement the non-effective compact-interval bounds following Theorem~\ref{thm:collision-existence}.

\subsection{The three excluded twelve-nut words}\label{app:exceptional}
Let $E=(-1,0,-1,1,-1)$ and let $E^R$ denote its reversal. The three excluded candidates are
\begin{center}\small
\begin{tabular}{@{}llccc@{}}\toprule
ID & Word & $k$ & $h$ & Necessary interval\\\midrule
\texttt{12-065}&$(E,-1,E)$&5&2&$(-3/5,-1/5)$\\
\texttt{12-066}&$(E,0,-E^R)$&6&0&$(-3/5,-1/5)$\\
\texttt{12-145}&$(E^R,0,-E)$&6&0&$(-4/5,-2/5)$\\\bottomrule
\end{tabular}
\end{center}
Their exclusion uses a certified conical model for $E$ and exact matching of two such models to a surviving central rod. The argument concerns a chosen path in the conical family, and does not assume that a balanced branch approaches an endpoint.

\paragraph{A positive common-cone model.}
The coefficient rods of $E$ are
\[
 (1,0),(0,1),(-1,1),(0,-1),(1,-2),(-1,3),(-2,5).
\]
For six nodes $\xi_i$, positive occupation variables $\zeta_i$ and moment weights $\omega_i=\zeta_i/\prod_{a\ne i}|\xi_i-\xi_a|$, put
$\mathsf t_j=(H_3^{[j]},H_2^{[j]})$ as in \eqref{eq:cluster-moments}. With $\xi_1=-1$, $\xi_6=1$ and $\zeta_6=1$, solve
\begin{equation}\label{eq:exception-common-cone}
 5\mathsf t_j=\lambda\bigl[(5p_j+2q_j)\mathsf t_0+q_j\mathsf t_6\bigr],
 \qquad j=1,\ldots,5.
\end{equation}
An interval certificate isolates an exact root near
\[
 \begin{aligned}
 \xi={}&(-1,-0.6043803602402567,\ 0.2317344053019192,\\
       &\hspace{2.6em}0.4639277115162177,\ 0.6405233684366971,1),\\
 \log\zeta={}&(0.6561651153410297,1.213438511134147,\\
       &0.1801659598474378,1.0700871944828563,1.372320385032065,0).
 \end{aligned}
\]
The decisive bounds are
\begin{equation}\label{eq:exception-positive-error}
 \begin{gathered}
 1.09420875801306906<\lambda_*<1.09420875801306907,\\
 0.09003150664870304<c_*:=\log\lambda_*
                         <0.09003150664870305.
 \end{gathered}
\end{equation}
The certificate uses a box of radius $10^{-30}$ about a 75-digit decimal centre, evaluated with 90-digit outward-rounded ball arithmetic. Its contraction bound is below $2.244\times10^{-23}$ and its centre displacement below $4.168\times10^{-75}$. Positive gaps, positive weights, nonvanishing denominators and the exterior frame are certified. The fixed-node five-variable direction Jacobian is separately certified nonsingular. Reversal and the particle-hole transformation \eqref{eq:particle-hole} give the required equivalent words, with the same multiplier and local invertibility.

\paragraph{Matching the two clusters.}
Place the left nodes at $-1+\varepsilon(\xi_i+1)$ and the right nodes at $1+\varepsilon(\eta_i-1)$. Each cluster has span $2\varepsilon$; the surviving central gap tends to two. It is convenient to use the Vandermonde occupation sums
\[
 Z_m^{[j]}=\sum_{|I|=m}(-1)^{|I\cap\{1,\ldots,j\}|}
 \Delta(\xi_I)\Delta(\xi_{I^c})\prod_{i\in I}\zeta_i.
\]
Here each Vandermonde uses positive differences in the ordered subset. These sums differ from the corresponding Hankel determinants by one common positive factor. On the left write $(Y,X)=(Z_3,Z_2)$; on the right write $(Q,P)=(Z_3,Z_2)$ for rank five and $(Q,P)=(Z_4,Z_3)$ for rank six. Reversal or reciprocal occupations implement the words in the table.

For rank five, give the two clusters occupations $\varepsilon A\zeta_i$ and $\varepsilon B\widetilde\zeta_i$. After removing their common leading factor, the limiting raw rods are
\begin{equation}\label{eq:exception-rank5}
 t^t=ABYQ-4XP-2BXQ+2AYP,\qquad
 t^\phi=BXQ+AYP.
\end{equation}
For rank six, the occupations are $\varepsilon A\zeta_i$ and $\varepsilon^{-1}B\widetilde\zeta_i$, and the corresponding formula is
\begin{equation}\label{eq:exception-rank6}
 t^t=ABYQ-XP-BXQ,\qquad
 t^\phi=AYP+\tfrac14BXQ.
\end{equation}
These identities follow by retaining the leading occupations in each cluster; the cross-cluster factor of a term with degrees $r,s$ is $2^{rs+(6-r)(6-s)}$. On a left rod the right sums take their initial values; on a right rod the left sums take their terminal values.

Put $\alpha=-2-5\beta$ for the first two words and $\alpha=-3-5\beta$ for the third, and set $s=\sqrt{1-\alpha^2}$. The interval in each row is precisely $|\alpha|<1$. Define $h_L=AY_0/X_0$ and $h_R=BQ_0/P_0$. AF alignment and central direction matching are solved by
\[
 (h_L,h_R)=
 \begin{cases}
 (2(1-\alpha)/s,\ 2(1+\alpha)/s),&k=5,\\
 ((1-\alpha)/(2s),\ 2(1+\alpha)/s),&k=6.
 \end{cases}
\]
Dividing by $t_0^\phi$ and subtracting the exterior time component gives
\[
 \mathsf r_0=(0,1),\qquad
 \mathsf r_6=(-2s(1+s),\alpha),\qquad
 L=\frac{2s(1+s)}5.
\]
The surviving rods are balanced in this lattice, while all ten internal cluster rods have multiplier $\lambda_*$. This is because \eqref{eq:exception-rank5} and \eqref{eq:exception-rank6} restrict to fixed linear frame maps of each cluster's two-component rods, preserving \eqref{eq:exception-common-cone}. The contracted nuts have lattice index five; no smoothness of the contracted background is required.

After extracting the common powers of $\varepsilon$, the direction equations extend analytically to zero. Their limiting Jacobian is block triangular: two certified five-variable internal blocks and one background direction derivative. The latter is nonzero because $\alpha/[s(1+s)]$ is strictly increasing on $(-1,1)$. The implicit-function theorem therefore supplies positive, exactly direction-matched data at small positive $\varepsilon$, with the period scale restored by a common rescaling.

In normalized length coordinates the resulting path satisfies
\begin{equation}\label{eq:exception-limit}
 \begin{aligned}
 q_\varepsilon&\longrightarrow
 q_0=(0,0,0,0,0,\,x_{6,*},\,0,0,0,0,0),
 &x_{6,*}&=\frac5{s(1+s)},\\
 b(q_\varepsilon)&\longrightarrow
 c_*(1,1,1,1,1,\,0,\,1,1,1,1,1).
 \end{aligned}
\end{equation}
Proposition~\ref{prop:boundary-obstruction} excludes balance for every twist in the respective interval. More quantitatively, put $S(x)=\sum_{j\ne6}x_j$. Gradient monotonicity and $c_*>9/100$ give
\begin{equation}\label{eq:exception-mismatch}
 \max_j|b_j(x)|\ge\frac9{100}
 \frac{S(x)}{S(x)+|x_6-x_{6,*}|}>0
 \qquad(x>0).
\end{equation}
This estimate is pointwise on the full positive orthant; it need not remain uniformly positive at its boundary. Neither the index five nor the occurrence of $E$ alone is the obstruction. The decisive feature is the matched limiting gradient \eqref{eq:exception-limit}.

\subsection{Verification of collision-block coverage}\label{app:collision-coverage}
Every proper consecutive block of the 356 remaining words is enumerated exactly, with reversal and orientation identified. The resulting counts are
\begin{center}
\begin{tabular}{@{}lr@{}}\toprule
Block test & Number\\\midrule
Nonparallel, $D(C)>1$ & 1595\\
Nonparallel, elementary reduction to a nonempty GH chain & 147\\
Nonempty even all-zero word & 4\\
Parallel exterior rods & 335\\\midrule
Total & 2081\\\bottomrule
\end{tabular}
\end{center}
Thus all 1,746 nonparallel blocks have negative invariant by Appendix~\ref{app:negative-tests}. Independent interval trials also verify a normalized opening cost at most $-\log(16/15)$ for each of these blocks. The limiting case of that certificate is the exact four-nut trial
\[
 C=(-1,1,-1),\quad z=(-1,-1/5,1/5,1),\quad
 \zeta=(1,3,3,1),\qquad a_1=a_2=a_3=15/16.
\]
The uniform opening-cost bound is a finite-catalogue statement, not a lower bound for the rod lengths. The existence proof needs only the strict signs, which the exact reduction census supplies without interval evaluation of every block.

\subsection{Organization by zero strings}\label{app:seed-organization}
Compress each nonempty zero run to one zero if its length is odd and to two zeros if its length is even, then quotient by reversal and orientation. The 356 surviving words have exactly 250 distinct compressed seeds, all themselves surviving catalogue members. Their counts by turning-point order are
\[
\begin{array}{c|rrrrrrrrrr}
N&3&4&5&6&7&8&9&10&11&12\\\hline
\text{seeds}&1&1&1&2&2&7&11&29&57&139.
\end{array}
\]
Theorem~\ref{thm:zero-extension} applies to every such seed. Ten seeds have no zero entries and hence have no nontrivial extension by this operation. The other 240 generate infinite families, with one nonnegative insertion parameter for each zero run. This is an exact organization of the finite catalogue; it is not a bounded-seed theorem for arbitrary order.

A coarser combinatorial compression deletes every even zero run entirely and retains one zero from each odd run. It produces 126 parity cores, again all surviving words in the catalogue. The inverse operation can create a zero run where none existed. Theorem~\ref{thm:zero-extension} only extends existing runs; preservation of all proper-block signs under creation of a new run requires an additional argument. The 126-core organization is therefore recorded as combinatorial data, not as a stronger existence theorem. Exact membership maps and the verification script accompany the catalogue.

\subsection{A higher-order test}\label{app:higher-order}
The AF counterexample \texttt{16-3583} has a separate 100-digit interval certificate at $\beta=-2/5$, with contraction bound below $1.913\cdot10^{-5}$. Its moment and original BZ evaluations agree independently. This tests the polygon-only exclusion beyond the exhaustive twelve-nut census.

\subsection{The electronic catalogue}
In the word tables, the letters $m,z,p$ mean $-1,0,+1$. Each code gives the \emph{entire} self-intersection word, not a half-word. The ID starts with the genuine nut count. The integer $h$ is the reduced phantom count, and the displayed masses and shortest lengths are evaluated at $\beta=(\beta_-+\beta_+)/2$. Every row has its own interval balance certificate. Identifiers are retained from the enumeration, so the omitted excluded cases leave gaps in the numbering. Six decimals are displayed for readability, while full centres and charges are provided in the data. Small displayed lengths are values at these representatives, not branch-wide lower bounds.

The electronic catalogue contains the complete words and rods, exact rational windows and equivalence maps, dressing degrees, certified midpoint data, nine-point continuations, zero-spin sections, asymptotic invariants, and exact component partitions. The solver and certificate scripts are included. The catalogue CSV provides a compact index; JSON files retain the full numerical state and evidence. The supplementary verification also includes the proper-block census and the exact zero-string insertion checks.

\begingroup\small\renewcommand{\arraystretch}{1.08}
\setlength{\tabcolsep}{5pt}
\begin{longtable}{@{}llrcrr@{}}
\caption{AF $N=3$ catalogue. Values refer to the midpoint of the displayed necessary interval.}\label{tab:catalogue-3}\\
\toprule
ID & Full word & $h$ & $\mathcal I_W$ & $M/L$ & $\min x_j$\\\midrule
\endfirsthead
\toprule
ID & Full word & $h$ & $\mathcal I_W$ & $M/L$ & $\min x_j$\\\midrule
\endhead
\midrule\multicolumn{6}{r}{Continued on next page}\\
\endfoot\bottomrule\endlastfoot
\texttt{3-001}&\texttt{mm}&1&$(-1,0)$&0.343146&0.485281\\
\end{longtable}\endgroup

\begingroup\small\renewcommand{\arraystretch}{1.08}
\setlength{\tabcolsep}{5pt}
\begin{longtable}{@{}llrcrr@{}}
\caption{AF $N=4$ catalogue. Values refer to the midpoint of the displayed necessary interval.}\label{tab:catalogue-4}\\
\toprule
ID & Full word & $h$ & $\mathcal I_W$ & $M/L$ & $\min x_j$\\\midrule
\endfirsthead
\toprule
ID & Full word & $h$ & $\mathcal I_W$ & $M/L$ & $\min x_j$\\\midrule
\endhead
\midrule\multicolumn{6}{r}{Continued on next page}\\
\endfoot\bottomrule\endlastfoot
\texttt{4-001}&\texttt{mzp}&0&$(-1,0)$&0.401899&0.401626\\
\end{longtable}\endgroup

\begingroup\small\renewcommand{\arraystretch}{1.08}
\setlength{\tabcolsep}{5pt}
\begin{longtable}{@{}llrcrr@{}}
\caption{AF $N=5$ catalogue. Values refer to the midpoint of the displayed necessary interval.}\label{tab:catalogue-5}\\
\toprule
ID & Full word & $h$ & $\mathcal I_W$ & $M/L$ & $\min x_j$\\\midrule
\endfirsthead
\toprule
ID & Full word & $h$ & $\mathcal I_W$ & $M/L$ & $\min x_j$\\\midrule
\endhead
\midrule\multicolumn{6}{r}{Continued on next page}\\
\endfoot\bottomrule\endlastfoot
\texttt{5-001}&\texttt{mzzm}&1&$(-1,0)$&0.453385&0.360446\\
\end{longtable}\endgroup

\begingroup\small\renewcommand{\arraystretch}{1.08}
\setlength{\tabcolsep}{5pt}
\begin{longtable}{@{}llrcrr@{}}
\caption{AF $N=6$ catalogue. Values refer to the midpoint of the displayed necessary interval.}\label{tab:catalogue-6}\\
\toprule
ID & Full word & $h$ & $\mathcal I_W$ & $M/L$ & $\min x_j$\\\midrule
\endfirsthead
\toprule
ID & Full word & $h$ & $\mathcal I_W$ & $M/L$ & $\min x_j$\\\midrule
\endhead
\midrule\multicolumn{6}{r}{Continued on next page}\\
\endfoot\bottomrule\endlastfoot
\texttt{6-001}&\texttt{mzzzp}&0&$(-1,0)$&0.499962&0.329479\\
\texttt{6-002}&\texttt{mpzmp}&0&$(-1,0)$&0.543368&0.231512\\
\texttt{6-003}&\texttt{mpppm}&0&$(-1,0)$&0.543451&0.232052\\
\end{longtable}\endgroup

\begingroup\small\renewcommand{\arraystretch}{1.08}
\setlength{\tabcolsep}{5pt}
\begin{longtable}{@{}llrcrr@{}}
\caption{AF $N=7$ catalogue. Values refer to the midpoint of the displayed necessary interval.}\label{tab:catalogue-7}\\
\toprule
ID & Full word & $h$ & $\mathcal I_W$ & $M/L$ & $\min x_j$\\\midrule
\endfirsthead
\toprule
ID & Full word & $h$ & $\mathcal I_W$ & $M/L$ & $\min x_j$\\\midrule
\endhead
\midrule\multicolumn{6}{r}{Continued on next page}\\
\endfoot\bottomrule\endlastfoot
\texttt{7-001}&\texttt{mzmmmp}&1&$(-1,0)$&0.577182&0.219730\\
\texttt{7-002}&\texttt{mzmzpm}&1&$(-1,0)$&0.577002&0.217740\\
\texttt{7-003}&\texttt{mzzzzm}&1&$(-1,0)$&0.543164&0.308133\\
\end{longtable}\endgroup

\begingroup\small\renewcommand{\arraystretch}{1.08}
\setlength{\tabcolsep}{5pt}
\begin{longtable}{@{}llrcrr@{}}
\caption{AF $N=8$ catalogue. Values refer to the midpoint of the displayed necessary interval.}\label{tab:catalogue-8}\\
\toprule
ID & Full word & $h$ & $\mathcal I_W$ & $M/L$ & $\min x_j$\\\midrule
\endfirsthead
\toprule
ID & Full word & $h$ & $\mathcal I_W$ & $M/L$ & $\min x_j$\\\midrule
\endhead
\midrule\multicolumn{6}{r}{Continued on next page}\\
\endfoot\bottomrule\endlastfoot
\texttt{8-001}&\texttt{mzmmmzm}&2&$(-1,0)$&0.609364&0.247301\\
\texttt{8-002}&\texttt{mzmzpzp}&0&$(-1,0)$&0.608935&0.245152\\
\texttt{8-003}&\texttt{mzzzzzp}&0&$(-1,0)$&0.583700&0.290401\\
\texttt{8-004}&\texttt{mzzpzmp}&0&$(-1,0)$&0.612542&0.206960\\
\texttt{8-005}&\texttt{mzzpppm}&0&$(-1,0)$&0.612734&0.209088\\
\texttt{8-006}&\texttt{mpmzpmp}&0&$(-1,-\tfrac{1}{3})$&0.759899&0.062720\\
\texttt{8-007}&\texttt{mpzzzmp}&0&$(-1,0)$&0.629377&0.190574\\
\texttt{8-008}&\texttt{mpzzppm}&0&$(-1,0)$&0.629833&0.191120\\
\texttt{8-009}&\texttt{mppzmmp}&0&$(-1,0)$&0.630264&0.196721\\
\end{longtable}\endgroup

\begingroup\small\renewcommand{\arraystretch}{1.08}
\setlength{\tabcolsep}{5pt}
\begin{longtable}{@{}llrcrr@{}}
\caption{AF $N=9$ catalogue. Values refer to the midpoint of the displayed necessary interval.}\label{tab:catalogue-9}\\
\toprule
ID & Full word & $h$ & $\mathcal I_W$ & $M/L$ & $\min x_j$\\\midrule
\endfirsthead
\toprule
ID & Full word & $h$ & $\mathcal I_W$ & $M/L$ & $\min x_j$\\\midrule
\endhead
\midrule\multicolumn{6}{r}{Continued on next page}\\
\endfoot\bottomrule\endlastfoot
\texttt{9-001}&\texttt{mzmmmzzp}&1&$(-1,0)$&0.643428&0.236463\\
\texttt{9-002}&\texttt{mzmmzzmp}&1&$(-1,0)$&0.659191&0.183966\\
\texttt{9-003}&\texttt{mzmmzppm}&1&$(-1,0)$&0.659606&0.189330\\
\texttt{9-004}&\texttt{mzmzzmmp}&1&$(-1,0)$&0.658594&0.188724\\
\texttt{9-005}&\texttt{mzmzzzpm}&1&$(-1,0)$&0.658110&0.182314\\
\texttt{9-006}&\texttt{mzmzpzzm}&1&$(-1,0)$&0.642996&0.234548\\
\texttt{9-007}&\texttt{mzmpzmpm}&1&$(-\tfrac{2}{3},0)$&0.776002&0.061100\\
\texttt{9-008}&\texttt{mzzzmmmp}&1&$(-1,0)$&0.647936&0.199665\\
\texttt{9-009}&\texttt{mzzzmzpm}&1&$(-1,0)$&0.647714&0.197092\\
\texttt{9-010}&\texttt{mzzzzzzm}&1&$(-1,0)$&0.622162&0.276686\\
\texttt{9-011}&\texttt{mpmmmzmp}&1&$(-1,-\tfrac{1}{3})$&0.775219&0.061168\\
\texttt{9-012}&\texttt{mpmmmppm}&1&$(-1,-\tfrac{1}{3})$&0.777742&0.060961\\
\texttt{9-013}&\texttt{mpzmmmmp}&1&$(-1,0)$&0.668120&0.186840\\
\texttt{9-014}&\texttt{mpzmmzpm}&1&$(-1,0)$&0.667889&0.187013\\
\texttt{9-015}&\texttt{mppppppm}&1&$(-1,0)$&0.668359&0.189605\\
\end{longtable}\endgroup

\begingroup\small\renewcommand{\arraystretch}{1.08}
\setlength{\tabcolsep}{5pt}
\begin{longtable}{@{}llrcrr@{}}
\caption{AF $N=10$ catalogue. Values refer to the midpoint of the displayed necessary interval.}\label{tab:catalogue-10}\\
\toprule
ID & Full word & $h$ & $\mathcal I_W$ & $M/L$ & $\min x_j$\\\midrule
\endfirsthead
\toprule
ID & Full word & $h$ & $\mathcal I_W$ & $M/L$ & $\min x_j$\\\midrule
\endhead
\midrule\multicolumn{6}{r}{Continued on next page}\\
\endfoot\bottomrule\endlastfoot
\texttt{10-001}&\texttt{mzmmmmmmp}&2&$(-1,0)$&0.696431&0.183188\\
\texttt{10-002}&\texttt{mzmmmmzpm}&2&$(-1,0)$&0.696179&0.180211\\
\texttt{10-003}&\texttt{mzmmmzzzm}&2&$(-1,0)$&0.677374&0.227383\\
\texttt{10-004}&\texttt{mzmmzzmzm}&2&$(-1,0)$&0.686915&0.213704\\
\texttt{10-005}&\texttt{mzmmzppzp}&0&$(-1,0)$&0.687871&0.219696\\
\texttt{10-006}&\texttt{mzmmpppmp}&0&$(-1,-\tfrac{1}{3})$&0.804364&0.058675\\
\texttt{10-007}&\texttt{mzmzmmmpm}&2&$(-\tfrac{2}{3},0)$&0.794166&0.059545\\
\texttt{10-008}&\texttt{mzmzzzpzp}&0&$(-1,0)$&0.685742&0.211773\\
\end{longtable}\endgroup

\begingroup\small\renewcommand{\arraystretch}{1.08}
\setlength{\tabcolsep}{5pt}
\begin{longtable}{@{}llrcrr@{}}
\caption{AF $N=11$ catalogue. Values refer to the midpoint of the displayed necessary interval.}\label{tab:catalogue-11}\\
\toprule
ID & Full word & $h$ & $\mathcal I_W$ & $M/L$ & $\min x_j$\\\midrule
\endfirsthead
\toprule
ID & Full word & $h$ & $\mathcal I_W$ & $M/L$ & $\min x_j$\\\midrule
\endhead
\midrule\multicolumn{6}{r}{Continued on next page}\\
\endfoot\bottomrule\endlastfoot
\texttt{11-001}&\texttt{mzmmmmmmzm}&3&$(-1,0)$&0.723634&0.211393\\
\texttt{11-002}&\texttt{mzmmmmzpzp}&1&$(-1,0)$&0.723059&0.208808\\
\texttt{11-003}&\texttt{mzmmmzzzzp}&1&$(-1,0)$&0.710725&0.219056\\
\texttt{11-004}&\texttt{mzmmmzpzmp}&1&$(-1,0)$&0.729618&0.174839\\
\texttt{11-005}&\texttt{mzmmmzpppm}&1&$(-1,0)$&0.729871&0.177888\\
\texttt{11-006}&\texttt{mzmmzmzpmp}&1&$(-1,-\tfrac{1}{3})$&0.828231&0.056270\\
\texttt{11-007}&\texttt{mzmmzzmzzp}&1&$(-1,0)$&0.717287&0.212727\\
\texttt{11-008}&\texttt{mzmmzzzzmp}&1&$(-1,0)$&0.728696&0.162797\\
\end{longtable}\endgroup

\begingroup\small\renewcommand{\arraystretch}{1.08}
\setlength{\tabcolsep}{5pt}
\begin{longtable}{@{}llrcrr@{}}
\caption{AF $N=12$ catalogue. Values refer to the midpoint of the displayed necessary interval.}\label{tab:catalogue-12}\\
\toprule
ID & Full word & $h$ & $\mathcal I_W$ & $M/L$ & $\min x_j$\\\midrule
\endfirsthead
\toprule
ID & Full word & $h$ & $\mathcal I_W$ & $M/L$ & $\min x_j$\\\midrule
\endhead
\midrule\multicolumn{6}{r}{Continued on next page}\\
\endfoot\bottomrule\endlastfoot
\texttt{12-001}&\texttt{mzmmmmmmzzp}&2&$(-1,0)$&0.753053&0.204872\\
\texttt{12-002}&\texttt{mzmmmmmzzmp}&2&$(-1,0)$&0.764987&0.160098\\
\texttt{12-003}&\texttt{mzmmmmmzppm}&2&$(-1,0)$&0.765392&0.165346\\
\texttt{12-004}&\texttt{mzmmmmzzmmp}&2&$(-1,0)$&0.764016&0.164460\\
\texttt{12-005}&\texttt{mzmmmmzzzpm}&2&$(-1,0)$&0.763530&0.157976\\
\texttt{12-006}&\texttt{mzmmmmzpzzm}&2&$(-1,0)$&0.752414&0.204483\\
\texttt{12-007}&\texttt{mzmmmmpzmpm}&2&$(-\tfrac{2}{3},0)$&0.868379&0.053916\\
\texttt{12-008}&\texttt{mzmmmzzmmmp}&2&$(-1,0)$&0.762081&0.173217\\
\end{longtable}\endgroup

\FloatBarrier
\subsection{Endpoint-partition atlas}
The left and right walls are those in the corresponding word table. Every entry is an exact rod-data prediction by \eqref{eq:app-end-block} for a surviving parent word. Component numbers refer to Table~\ref{tab:block-types}; the atlas does not assign metric existence or convergence to a component number.

\begingroup\small\renewcommand{\arraystretch}{1.08}
\setlength{\tabcolsep}{5pt}
\begin{longtable}{@{}lcc@{}}
\caption{Exact saturated-rod partitions for $N=3$. Component numbers are defined in Table~\ref{tab:block-types}; these are candidate limit data.}\label{tab:ends-3}\\
\toprule
ID & Left wall partition & Right wall partition\\\midrule
\endfirsthead
\toprule
ID & Left wall partition & Right wall partition\\\midrule
\endhead
\midrule\multicolumn{3}{r}{Continued on next page}\\
\endfoot\bottomrule\endlastfoot
\texttt{3-001}&$1,3$&$3,1$\\
\end{longtable}\endgroup

\begingroup\small\renewcommand{\arraystretch}{1.08}
\setlength{\tabcolsep}{5pt}
\begin{longtable}{@{}lcc@{}}
\caption{Exact saturated-rod partitions for $N=4$. Component numbers are defined in Table~\ref{tab:block-types}; these are candidate limit data.}\label{tab:ends-4}\\
\toprule
ID & Left wall partition & Right wall partition\\\midrule
\endfirsthead
\toprule
ID & Left wall partition & Right wall partition\\\midrule
\endhead
\midrule\multicolumn{3}{r}{Continued on next page}\\
\endfoot\bottomrule\endlastfoot
\texttt{4-001}&$1,2,1$&$3,3$\\
\end{longtable}\endgroup

\begingroup\small\renewcommand{\arraystretch}{1.08}
\setlength{\tabcolsep}{5pt}
\begin{longtable}{@{}lcc@{}}
\caption{Exact saturated-rod partitions for $N=5$. Component numbers are defined in Table~\ref{tab:block-types}; these are candidate limit data.}\label{tab:ends-5}\\
\toprule
ID & Left wall partition & Right wall partition\\\midrule
\endfirsthead
\toprule
ID & Left wall partition & Right wall partition\\\midrule
\endhead
\midrule\multicolumn{3}{r}{Continued on next page}\\
\endfoot\bottomrule\endlastfoot
\texttt{5-001}&$1,2,3$&$3,2,1$\\
\end{longtable}\endgroup

\begingroup\small\renewcommand{\arraystretch}{1.08}
\setlength{\tabcolsep}{5pt}
\begin{longtable}{@{}lcc@{}}
\caption{Exact saturated-rod partitions for $N=6$. Component numbers are defined in Table~\ref{tab:block-types}; these are candidate limit data.}\label{tab:ends-6}\\
\toprule
ID & Left wall partition & Right wall partition\\\midrule
\endfirsthead
\toprule
ID & Left wall partition & Right wall partition\\\midrule
\endhead
\midrule\multicolumn{3}{r}{Continued on next page}\\
\endfoot\bottomrule\endlastfoot
\texttt{6-001}&$1,2,2,1$&$3,2,3$\\
\texttt{6-002}&$1,3,3,1$&$3,1,1,3$\\
\texttt{6-003}&$1,3,1,3$&$3,1,3,1$\\
\end{longtable}\endgroup

\begingroup\small\renewcommand{\arraystretch}{1.08}
\setlength{\tabcolsep}{5pt}
\begin{longtable}{@{}lcc@{}}
\caption{Exact saturated-rod partitions for $N=7$. Component numbers are defined in Table~\ref{tab:block-types}; these are candidate limit data.}\label{tab:ends-7}\\
\toprule
ID & Left wall partition & Right wall partition\\\midrule
\endfirsthead
\toprule
ID & Left wall partition & Right wall partition\\\midrule
\endhead
\midrule\multicolumn{3}{r}{Continued on next page}\\
\endfoot\bottomrule\endlastfoot
\texttt{7-001}&$1,2,1,3,1$&$3,3,1,3$\\
\texttt{7-002}&$1,2,1,1,3$&$3,3,3,1$\\
\texttt{7-003}&$1,2,2,3$&$3,2,2,1$\\
\end{longtable}\endgroup

\begingroup\small\renewcommand{\arraystretch}{1.08}
\setlength{\tabcolsep}{5pt}
\begin{longtable}{@{}lcc@{}}
\caption{Exact saturated-rod partitions for $N=8$. Component numbers are defined in Table~\ref{tab:block-types}; these are candidate limit data.}\label{tab:ends-8}\\
\toprule
ID & Left wall partition & Right wall partition\\\midrule
\endfirsthead
\toprule
ID & Left wall partition & Right wall partition\\\midrule
\endhead
\midrule\multicolumn{3}{r}{Continued on next page}\\
\endfoot\bottomrule\endlastfoot
\texttt{8-001}&$1,2,1,3,3$&$3,3,1,2,1$\\
\texttt{8-002}&$1,2,1,1,2,1$&$3,3,3,3$\\
\texttt{8-003}&$1,2,2,2,1$&$3,2,2,3$\\
\texttt{8-004}&$1,2,3,3,1$&$3,2,1,1,3$\\
\texttt{8-005}&$1,2,3,1,3$&$3,2,1,3,1$\\
\texttt{8-006}&$1,3,1,1,3,1$&$6,6$\\
\texttt{8-007}&$1,3,2,3,1$&$3,1,1,1,1,3$\\
\texttt{8-008}&$1,3,2,1,3$&$3,1,1,1,3,1$\\
\texttt{8-009}&$1,3,1,1,3,1$&$3,1,2,1,3$\\
\end{longtable}\endgroup

\begingroup\small\renewcommand{\arraystretch}{1.08}
\setlength{\tabcolsep}{5pt}
\begin{longtable}{@{}lcc@{}}
\caption{Exact saturated-rod partitions for $N=9$. Component numbers are defined in Table~\ref{tab:block-types}; these are candidate limit data.}\label{tab:ends-9}\\
\toprule
ID & Left wall partition & Right wall partition\\\midrule
\endfirsthead
\toprule
ID & Left wall partition & Right wall partition\\\midrule
\endhead
\midrule\multicolumn{3}{r}{Continued on next page}\\
\endfoot\bottomrule\endlastfoot
\texttt{9-001}&$1,2,1,3,2,1$&$3,3,1,2,3$\\
\texttt{9-002}&$1,2,1,2,3,1$&$3,3,1,1,1,3$\\
\texttt{9-003}&$1,2,1,2,1,3$&$3,3,1,1,3,1$\\
\texttt{9-004}&$1,2,1,1,1,3,1$&$3,3,2,1,3$\\
\texttt{9-005}&$1,2,1,1,1,1,3$&$3,3,2,3,1$\\
\texttt{9-006}&$1,2,1,1,2,3$&$3,3,3,2,1$\\
\texttt{9-007}&$7,6$&$3,3,1,1,3,1$\\
\texttt{9-008}&$1,2,2,1,3,1$&$3,2,3,1,3$\\
\texttt{9-009}&$1,2,2,1,1,3$&$3,2,3,3,1$\\
\texttt{9-010}&$1,2,2,2,3$&$3,2,2,2,1$\\
\texttt{9-011}&$1,3,1,3,3,1$&$6,7$\\
\texttt{9-012}&$1,3,1,3,1,3$&$6,8$\\
\texttt{9-013}&$1,3,3,1,3,1$&$3,1,1,3,1,3$\\
\texttt{9-014}&$1,3,3,1,1,3$&$3,1,1,3,3,1$\\
\texttt{9-015}&$1,3,1,3,1,3$&$3,1,3,1,3,1$\\
\end{longtable}\endgroup

\begingroup\small\renewcommand{\arraystretch}{1.08}
\setlength{\tabcolsep}{5pt}
\begin{longtable}{@{}lcc@{}}
\caption{Exact saturated-rod partitions for $N=10$. Component numbers are defined in Table~\ref{tab:block-types}; these are candidate limit data.}\label{tab:ends-10}\\
\toprule
ID & Left wall partition & Right wall partition\\\midrule
\endfirsthead
\toprule
ID & Left wall partition & Right wall partition\\\midrule
\endhead
\midrule\multicolumn{3}{r}{Continued on next page}\\
\endfoot\bottomrule\endlastfoot
\texttt{10-001}&$1,2,1,3,1,3,1$&$3,3,1,3,1,3$\\
\texttt{10-002}&$1,2,1,3,1,1,3$&$3,3,1,3,3,1$\\
\texttt{10-003}&$1,2,1,3,2,3$&$3,3,1,2,2,1$\\
\texttt{10-004}&$1,2,1,2,3,3$&$3,3,1,1,1,2,1$\\
\texttt{10-005}&$1,2,1,2,1,2,1$&$3,3,1,1,3,3$\\
\texttt{10-006}&$1,2,1,3,1,3,1$&$9,6$\\
\texttt{10-007}&$10,6$&$3,3,3,1,3,1$\\
\texttt{10-008}&$1,2,1,1,1,1,2,1$&$3,3,2,3,3$\\
\end{longtable}\endgroup

\begingroup\small\renewcommand{\arraystretch}{1.08}
\setlength{\tabcolsep}{5pt}
\begin{longtable}{@{}lcc@{}}
\caption{Exact saturated-rod partitions for $N=11$. Component numbers are defined in Table~\ref{tab:block-types}; these are candidate limit data.}\label{tab:ends-11}\\
\toprule
ID & Left wall partition & Right wall partition\\\midrule
\endfirsthead
\toprule
ID & Left wall partition & Right wall partition\\\midrule
\endhead
\midrule\multicolumn{3}{r}{Continued on next page}\\
\endfoot\bottomrule\endlastfoot
\texttt{11-001}&$1,2,1,3,1,3,3$&$3,3,1,3,1,2,1$\\
\texttt{11-002}&$1,2,1,3,1,1,2,1$&$3,3,1,3,3,3$\\
\texttt{11-003}&$1,2,1,3,2,2,1$&$3,3,1,2,2,3$\\
\texttt{11-004}&$1,2,1,3,3,3,1$&$3,3,1,2,1,1,3$\\
\texttt{11-005}&$1,2,1,3,3,1,3$&$3,3,1,2,1,3,1$\\
\texttt{11-006}&$1,2,1,2,1,1,3,1$&$15,6$\\
\texttt{11-007}&$1,2,1,2,3,2,1$&$3,3,1,1,1,2,3$\\
\texttt{11-008}&$1,2,1,2,2,3,1$&$3,3,1,1,1,1,1,3$\\
\end{longtable}\endgroup

\begingroup\small\renewcommand{\arraystretch}{1.08}
\setlength{\tabcolsep}{5pt}
\begin{longtable}{@{}lcc@{}}
\caption{Exact saturated-rod partitions for $N=12$. Component numbers are defined in Table~\ref{tab:block-types}; these are candidate limit data.}\label{tab:ends-12}\\
\toprule
ID & Left wall partition & Right wall partition\\\midrule
\endfirsthead
\toprule
ID & Left wall partition & Right wall partition\\\midrule
\endhead
\midrule\multicolumn{3}{r}{Continued on next page}\\
\endfoot\bottomrule\endlastfoot
\texttt{12-001}&$1,2,1,3,1,3,2,1$&$3,3,1,3,1,2,3$\\
\texttt{12-002}&$1,2,1,3,1,2,3,1$&$3,3,1,3,1,1,1,3$\\
\texttt{12-003}&$1,2,1,3,1,2,1,3$&$3,3,1,3,1,1,3,1$\\
\texttt{12-004}&$1,2,1,3,1,1,1,3,1$&$3,3,1,3,2,1,3$\\
\texttt{12-005}&$1,2,1,3,1,1,1,1,3$&$3,3,1,3,2,3,1$\\
\texttt{12-006}&$1,2,1,3,1,1,2,3$&$3,3,1,3,3,2,1$\\
\texttt{12-007}&$22,6$&$3,3,1,3,1,1,3,1$\\
\texttt{12-008}&$1,2,1,3,2,1,3,1$&$3,3,1,2,3,1,3$\\
\end{longtable}\endgroup

\FloatBarrier

\section{ALF catalogue through eight turning points}\label{app:ALF-system}\label{sec:ALF}

This appendix records 52 certified isolated regular solutions on 49 admissible rod words through eight turning points. Three eight-nut words support two different asymptotic circle markings. All 49 words have negative homogeneous cluster invariant, while the sixteen remaining candidates have positive invariant; the general implication from this sign to ALF existence remains open. The catalogue tests the polygon organization of Section~\ref{sec:classification-roots} and supplies marked components for the endpoint construction in Section~\ref{sec:endpoint-composition}.

The nonparallel exterior problem retains the positive moment reconstruction and changes the asymptotic marking. Here ALF is used in the finite-cover sense: the end has a circle direction of bounded length, and exceptional orbits of that circle are permitted. This is also the convention used in \cite{DMT,ADN}; the census is not restricted to principal circle bundles over a smooth two-sphere.

\subsection{The bounded circle at infinity}
Keep the lattice marking of Section~\ref{sec:regularity}. Let the last coefficient vector be $v_N=(p_N,q_N)$ and let $k$ be the lower degree of a positive reduced sector. Its normalized exterior angular component is $\sigma=(-1)^k$. For nonparallel exteriors, $q_N\ne0$, and exterior matching forces
\begin{equation}\label{eq:ALF-beta}
 p_N-q_N\beta=\sigma,\qquad
 \beta=\frac{p_N-\sigma}{q_N}.
\end{equation}
Thus $\beta$ is fixed by the word and sector. The AF parameter survives because its exterior directions are parallel; that argument does not supply a continuous twist family in the nonparallel problem.

The primitive bounded circle is represented in coefficient coordinates by
\begin{equation}\label{eq:ALF-fibre}
 f_\infty=\operatorname{prim}(v_N-\sigma v_0),
\end{equation}
where $\operatorname{prim}$ divides an integral vector by the greatest common divisor of its components and its overall sign is immaterial. Its physical angular component vanishes by \eqref{eq:ALF-beta}. The pair consisting of the ordered rod word and this marked fibre distinguishes sectors that the direction polygon alone forgets.

For example, the word
\[
 (-1,0,0,0,-1,1,-1),\qquad v_N=(2,-5),
\]
has certified balance roots in sectors $(k,h)=(3,2)$ and $(4,0)$. Their twist values are $-3/5$ and $-1/5$, with primitive fibre vectors $(3,-5)$ and $(1,-5)$, respectively. The same ordered word therefore supports two distinct asymptotic circle markings. Conversely, a polygon and a fibre need not determine the word: the five-nut words $(-1,0,-1,1)$ and $(-1,1,1,-1)$ share both but have different ordered walks.

\subsection{The finite census}
The ALF enumeration imposes the ternary rule, all block determinant inequalities, nonparallel exteriors, the strict angular bounds at \eqref{eq:ALF-beta}, and the prefix winding test in each positive reduced sector $k\ge1$. Reversal and orientation are quotiented with the marking retained. Table~\ref{tab:ALF-census} is exhaustive for these tests through eight genuine turning points.

\begin{table}[htbp]\centering\small
\begin{tabular}{@{}rrrr@{}}\toprule
$N$ & Candidate words & Words with certified roots & Other words\\\midrule
4&1&1&0\\
5&2&2&0\\
6&6&5&1\\
7&16&11&5\\
8&40&30&10\\\midrule
Total&65&49&16\\\bottomrule
\end{tabular}
\caption{The positive reduced ALF census. Three eight-nut words have roots in two sectors, giving 52 certified word--sector points. The last column does not assert nonexistence.}\label{tab:ALF-census}
\end{table}

Each certificate proves an exact root of the normalized rod equations, positive weights, separated nodes and nonzero required frame factors. Its invertible normalized Jacobian proves local isolation after translation and scale are removed. The dictionary in Section~\ref{subsec:ADN-dictionary} and \cite[Proposition~11.4]{ADN} then give smooth complete ALF metrics in the stated finite-cover sense. Local isolation does not imply global uniqueness in a sector. The equations and scope of the verification are given below.

These isolated roots do not imply that ALF moduli spaces are discrete. Collinear multi--Taub--NUT metrics retain continuous centre-separation parameters; they occupy the excluded maximal-phantom rank-zero sector. Their $-2$ chains \cite{GH} also explain why the ternary restriction cannot be extended to arbitrary phantom limits or to all ALF geometries.

The two five-nut words are $(-1,0,-1,1)$ and $(-1,1,1,-1)$. Up to reversal and orientation, the two five-node ALF examples in \cite[Section~7]{DMT} have these same words, and their printed nodes and residues agree with our points to the stated precision. That work checks regularity numerically; here the normalized balance equations have interval certificates and completion follows from the positive-family theorem above. The four-nut component also appears in explicit radical form in \cite[Section~11.5]{ADN}. These overlaps provide independent comparisons of the twistor, BZ and moment descriptions.

\subsection{The polygon structure of the ALF data}
All 65 words give crystallographic arrangements by the exact rod audit. Their polygon distribution is summarized in Table~\ref{tab:ALF-polygons}. The original four- and five-nut endpoint components have the pentagon arrangement.

\begin{table}[htbp]\centering\small
\begin{tabular}{@{}lrr@{}}\toprule
Triangulation & Candidate words & Words with certified roots\\\midrule
Pentagon &25&25\\
Fan hexagon &15&15\\
Central-triangle hexagon &1&1\\
Zigzag hexagon &14&3\\
Fan heptagon &5&5\\
Non-fan heptagon, $(1,2,2,3,1,2,4)$ &5&0\\\bottomrule
\end{tabular}
\caption{ALF direction arrangements through $N=8$. Polygon type alone does not determine balance.}\label{tab:ALF-polygons}
\end{table}

The homogeneous sign test provides an additional distinction. All 49 words with balance roots have certified $\mu<0$, while the other 16 have certified $\mu>0$. The sign belongs to the word's normalized auxiliary homogeneous model and does not distinguish the two fibre markings of a doubled word. This agreement suggests an ALF existence criterion based on $\mu$, but its necessity and sufficiency require a general comparison with the finite ALF family. The neutral rank-zero GH case would need a separate statement.

\subsection{The finite equations}
Use the reduced moment rods $\mathsf t_j$ defined in \eqref{eq:app-moment-rods}. Fix the outer nodes at $z_1=-1$, $z_N=1$ and vary all $N$ positive occupation weights and all $N-2$ interior nodes independently. For an admissible marked sector the balance equations are
\begin{equation}\label{eq:ALF-balance-system}
 \mathsf t_j=p_j\mathsf t_0+q_j\mathsf t_1,
 \qquad j=2,\ldots,N.
\end{equation}
There are $2N-2$ scalar equations for $2N-2$ variables. The last angular equation determines the ratio of the first two angular components and hence the fixed value \eqref{eq:ALF-beta}; no separate AF alignment equation is imposed. The same common exterior shear and normalization as in \eqref{eq:app-units} recovers $L>0$. The normalized lengths are $(z_{j+1}-z_j)/L$.

The search uses logarithmic weights and gap ratios. Outward-rounded interval contraction boxes, evaluated with Arb \cite{Arb}, certify exact roots and an invertible Jacobian. Positivity of the occupations, separation of the nodes and every denominator used in the equations are enclosed on each box. These checks prove local isolation after the two stated coordinate normalizations. They do not exclude a second root elsewhere in the same sector.

\subsection{Scope of the catalogue}
The supplement contains all 65 words, the allowed reduced sectors, 52 balance certificates, the primitive rod and fibre data, and the separate homogeneous sign certificates. The sign comparison uses \eqref{eq:cluster-trial-bounds}; a failed common-cone search is not assigned a sign without a separate trial certificate. In particular, both positive boundary-minimization cases in the eight-nut census have independently certified positive lower bounds.

The original endpoint labels $\mathcal C_4$, $\mathcal C_5^a$, $\mathcal C_5^b$, and $\mathcal C_6^{a},\ldots,\mathcal C_6^{e}$ are retained in the electronic data. The four-nut component has an exact radical metric with recorded smooth completion; an independent form is given in \cite[Section~11.5]{ADN}. The symmetric five-nut example has exact algebraic balance data specified by a degree-sixteen polynomial and an isolating interval. Such algebraic specifications describe exact solutions without requiring radical expressions. The later finite census supplies interval balance certificates for the original eight endpoint types.

The whole-axis and end analysis is supplied analytically by \cite[Propositions~7.5, 8.7--8.9 and 11.4]{ADN}, through the normalization dictionary in Section~\ref{subsec:ADN-dictionary}. Applied to the certified positive, separated roots of \eqref{eq:ALF-balance-system}, it gives smooth complete ALF metrics. Thus bulk grid sampling is not used to establish their global positivity or completion. The interval certificates prove existence and local isolation; they do not prove global uniqueness in each marked sector. The census excludes rank-zero multi--Taub--NUT and the classical one- and two-nut families, and therefore is not an exhaustive classification of ALF instantons.

\subsection{Certified word--sector points}\label{app:ALF-table}
Table~\ref{tab:ALF-points} lists all 52 certified points. The identifier includes the number of genuine nuts. The rank $k$ gives the phantom count $h=N-2k$ and BZ slot count $2(N-k)$. The fibre is written in the coefficient basis $(A,-B)$, up to overall sign. Masses are rounded values of $M/L$ at the certified roots, not exact algebraic expressions. Repeated identifiers denote different marked sectors, not duplicated metrics. Words are displayed in the orientation used for the rank-$k$ solve; for \texttt{ALF7-016}, this reverses every sign of the canonical catalogue word. All values of $\beta$ and $f_\infty$ refer to the displayed orientation.

{\footnotesize\setlength{\tabcolsep}{4pt}
\begin{longtable}{@{}llrrcr@{}}
\caption{The ALF balance catalogue through eight genuine turning points.}\label{tab:ALF-points}\\
\toprule
Identifier & Self-intersection word & $k$ & $\beta$ & $f_\infty$ & $M/L$\\\midrule
\endfirsthead
\multicolumn{6}{c}{\tablename\ \thetable\ (continued)}\\
\toprule Identifier & Self-intersection word & $k$ & $\beta$ & $f_\infty$ & $M/L$\\\midrule
\endhead
\midrule\multicolumn{6}{r}{Continued on next page}\\\endfoot
\bottomrule\endlastfoot
\texttt{ALF4-001} & $(-1,1,-1)$ & 2 & $-\frac{1}{3}$ & $(1,-3)$ & 0.75300368\\
\texttt{ALF5-001} & $(-1,0,-1,1)$ & 2 & $-\frac{2}{3}$ & $(2,-3)$ & 0.76553516\\
\texttt{ALF5-002} & $(-1,1,1,-1)$ & 2 & $-\frac{2}{3}$ & $(2,-3)$ & 0.76737186\\
\texttt{ALF6-001} & $(-1,0,-1,-1,1)$ & 3 & $-\frac{1}{3}$ & $(1,-3)$ & 0.78669946\\
\texttt{ALF6-002} & $(-1,0,-1,0,-1)$ & 2 & $-\frac{2}{3}$ & $(2,-3)$ & 0.78287974\\
\texttt{ALF6-004} & $(-1,0,0,1,-1)$ & 3 & $-\frac{1}{3}$ & $(1,-3)$ & 0.78555972\\
\texttt{ALF6-005} & $(-1,1,-1,-1,1)$ & 3 & $-\frac{1}{2}$ & $(1,-2)$ & 1.00702968\\
\texttt{ALF6-006} & $(-1,1,0,1,-1)$ & 3 & $-\frac{1}{2}$ & $(1,-2)$ & 1.00487937\\
\texttt{ALF7-001} & $(-1,0,-1,-1,0,-1)$ & 3 & $-\frac{1}{3}$ & $(1,-3)$ & 0.80670110\\
\texttt{ALF7-002} & $(-1,0,-1,-1,1,-1)$ & 3 & $-\frac{1}{2}$ & $(1,-2)$ & 1.02175123\\
\texttt{ALF7-003} & $(-1,0,-1,0,-1,1)$ & 3 & $-\frac{1}{2}$ & $(1,-2)$ & 1.01330007\\
\texttt{ALF7-004} & $(-1,0,-1,0,0,1)$ & 3 & $-\frac{2}{3}$ & $(2,-3)$ & 0.80494940\\
\texttt{ALF7-006} & $(-1,0,-1,1,1,-1)$ & 3 & $-\frac{1}{2}$ & $(1,-2)$ & 1.01516190\\
\texttt{ALF7-007} & $(-1,0,0,0,-1,1)$ & 3 & $-\frac{2}{3}$ & $(2,-3)$ & 0.80944995\\
\texttt{ALF7-009} & $(-1,0,0,1,1,-1)$ & 3 & $-\frac{2}{3}$ & $(2,-3)$ & 0.81025821\\
\texttt{ALF7-010} & $(-1,1,-1,-1,1,-1)$ & 3 & $-\frac{3}{5}$ & $(3,-5)$ & 1.25277890\\
\texttt{ALF7-011} & $(-1,1,0,-1,-1,1)$ & 3 & $-\frac{2}{3}$ & $(2,-3)$ & 0.81829517\\
\texttt{ALF7-013} & $(-1,1,0,0,1,-1)$ & 3 & $-\frac{2}{3}$ & $(2,-3)$ & 0.81817727\\
\texttt{ALF7-016} & $(1,-1,-1,-1,-1,1)$ & 3 & $\frac{1}{3}$ & $(1,3)$ & 0.81836746\\
\texttt{ALF8-001} & $(-1,0,-1,-1,-1,-1,1)$ & 3 & $-\frac{2}{3}$ & $(2,-3)$ & 0.84027110\\
\texttt{ALF8-003} & $(-1,0,-1,-1,0,0,1)$ & 4 & $-\frac{1}{3}$ & $(1,-3)$ & 0.83091072\\
\texttt{ALF8-004} & $(-1,0,-1,-1,0,1,-1)$ & 3 & $-\frac{2}{3}$ & $(2,-3)$ & 0.84009294\\
\texttt{ALF8-005} & $(-1,0,-1,-1,1,0,1)$ & 4 & $-\frac{1}{2}$ & $(1,-2)$ & 1.03038698\\
\texttt{ALF8-008} & $(-1,0,-1,0,-1,0,-1)$ & 3 & $-\frac{1}{2}$ & $(1,-2)$ & 1.02312317\\
\texttt{ALF8-010} & $(-1,0,-1,0,0,-1,1)$ & 4 & $-\frac{1}{3}$ & $(1,-3)$ & 0.83831812\\
\texttt{ALF8-011} & $(-1,0,-1,0,0,0,-1)$ & 3 & $-\frac{2}{3}$ & $(2,-3)$ & 0.83041074\\
\texttt{ALF8-013} & $(-1,0,-1,0,1,1,-1)$ & 4 & $-\frac{1}{3}$ & $(1,-3)$ & 0.83857917\\
\texttt{ALF8-017} & $(-1,0,-1,1,1,-1,1)$ & 4 & $-\frac{2}{5}$ & $(2,-5)$ & 1.25854700\\
\texttt{ALF8-018} & $(-1,0,0,0,-1,-1,1)$ & 4 & $-\frac{1}{3}$ & $(1,-3)$ & 0.83474379\\
\texttt{ALF8-019} & $(-1,0,0,0,-1,1,-1)$ & 3 & $-\frac{3}{5}$ & $(3,-5)$ & 1.25052615\\
\texttt{ALF8-019} & $(-1,0,0,0,-1,1,-1)$ & 4 & $-\frac{1}{5}$ & $(1,-5)$ & 1.25087950\\
\texttt{ALF8-020} & $(-1,0,0,0,0,1,-1)$ & 4 & $-\frac{1}{3}$ & $(1,-3)$ & 0.83412893\\
\texttt{ALF8-021} & $(-1,0,0,1,-1,-1,1)$ & 4 & $-\frac{1}{2}$ & $(1,-2)$ & 1.03090301\\
\texttt{ALF8-022} & $(-1,0,0,1,0,0,-1)$ & 4 & $-\frac{1}{3}$ & $(1,-3)$ & 0.82759798\\
\texttt{ALF8-023} & $(-1,0,0,1,0,1,-1)$ & 4 & $-\frac{1}{2}$ & $(1,-2)$ & 1.02920163\\
\texttt{ALF8-024} & $(-1,0,0,1,1,-1,1)$ & 4 & $-\frac{1}{2}$ & $(1,-2)$ & 1.03751978\\
\texttt{ALF8-025} & $(-1,1,-1,-1,-1,-1,1)$ & 3 & $-\frac{4}{5}$ & $(4,-5)$ & 1.25132304\\
\texttt{ALF8-025} & $(-1,1,-1,-1,-1,-1,1)$ & 4 & $-\frac{2}{5}$ & $(2,-5)$ & 1.25280610\\
\texttt{ALF8-026} & $(-1,1,-1,-1,-1,1,-1)$ & 3 & $-\frac{2}{3}$ & $(2,-3)$ & 0.94312493\\
\texttt{ALF8-028} & $(-1,1,-1,-1,0,1,-1)$ & 3 & $-\frac{4}{5}$ & $(4,-5)$ & 1.25135225\\
\texttt{ALF8-028} & $(-1,1,-1,-1,0,1,-1)$ & 4 & $-\frac{2}{5}$ & $(2,-5)$ & 1.25291835\\
\texttt{ALF8-030} & $(-1,1,-1,0,0,-1,1)$ & 4 & $-\frac{1}{2}$ & $(1,-2)$ & 1.04786289\\
\texttt{ALF8-031} & $(-1,1,-1,0,1,1,-1)$ & 4 & $-\frac{1}{2}$ & $(1,-2)$ & 1.04738469\\
\texttt{ALF8-032} & $(-1,1,0,-1,-1,-1,1)$ & 4 & $-\frac{1}{3}$ & $(1,-3)$ & 0.84649583\\
\texttt{ALF8-033} & $(-1,1,0,-1,0,1,-1)$ & 4 & $-\frac{1}{3}$ & $(1,-3)$ & 0.84619067\\
\texttt{ALF8-034} & $(-1,1,0,0,-1,-1,1)$ & 4 & $-\frac{1}{2}$ & $(1,-2)$ & 1.02642786\\
\texttt{ALF8-035} & $(-1,1,0,0,0,1,-1)$ & 4 & $-\frac{1}{2}$ & $(1,-2)$ & 1.02507012\\
\texttt{ALF8-036} & $(-1,1,0,1,-1,-1,1)$ & 4 & $-\frac{3}{5}$ & $(3,-5)$ & 1.25778951\\
\texttt{ALF8-037} & $(-1,1,0,1,0,1,-1)$ & 4 & $-\frac{3}{5}$ & $(3,-5)$ & 1.25597809\\
\texttt{ALF8-038} & $(-1,1,1,-1,1,1,-1)$ & 4 & $-\frac{2}{5}$ & $(2,-5)$ & 1.25967010\\
\texttt{ALF8-039} & $(-1,1,1,0,1,1,-1)$ & 4 & $-\frac{1}{2}$ & $(1,-2)$ & 1.02779163\\
\texttt{ALF8-040} & $(-1,1,1,1,1,1,-1)$ & 4 & $-\frac{2}{3}$ & $(2,-3)$ & 0.84680026\\
\end{longtable}}

The three doubled words are \texttt{ALF8-019}, \texttt{ALF8-025} and \texttt{ALF8-028}. Their two values of $\beta$ give different primitive fibres. The exact ALF equations, rather than a choice of polynomial root label alone, distinguish their solutions.

\subsection{Candidates without an ALF balance certificate}
Table~\ref{tab:ALF-other} lists the other sixteen words. Each has certified positive invariant for its normalized homogeneous model. These are not certified nonexistence statements for the full ALF problem: the general implication from the homogeneous sign to finite ALF balance remains open.

{\small
\begin{longtable}{@{}ll@{}}
\caption{ALF candidates with positive homogeneous invariant and no certified balance root.}\label{tab:ALF-other}\\
\toprule Identifier & Self-intersection word\\\midrule
\endfirsthead
\multicolumn{2}{c}{\tablename\ \thetable\ (continued)}\\
\toprule Identifier & Self-intersection word\\\midrule\endhead
\bottomrule\endfoot
\texttt{ALF6-003} & $(-1,0,-1,1,-1)$\\
\texttt{ALF7-005} & $(-1,0,-1,1,0,1)$\\
\texttt{ALF7-008} & $(-1,0,0,1,-1,1)$\\
\texttt{ALF7-012} & $(-1,1,0,0,-1,1)$\\
\texttt{ALF7-014} & $(-1,1,0,1,1,-1)$\\
\texttt{ALF7-015} & $(-1,1,1,-1,-1,1)$\\
\texttt{ALF8-002} & $(-1,0,-1,-1,0,-1,1)$\\
\texttt{ALF8-006} & $(-1,0,-1,-1,1,1,-1)$\\
\texttt{ALF8-007} & $(-1,0,-1,0,-1,-1,1)$\\
\texttt{ALF8-009} & $(-1,0,-1,0,-1,1,-1)$\\
\texttt{ALF8-012} & $(-1,0,-1,0,0,1,-1)$\\
\texttt{ALF8-014} & $(-1,0,-1,1,-1,-1,1)$\\
\texttt{ALF8-015} & $(-1,0,-1,1,0,0,-1)$\\
\texttt{ALF8-016} & $(-1,0,-1,1,0,1,-1)$\\
\texttt{ALF8-027} & $(-1,1,-1,-1,0,-1,1)$\\
\texttt{ALF8-029} & $(-1,1,-1,-1,1,1,-1)$\\
\end{longtable}}

\subsection{The twelve-nut heptagon refinement}
A targeted refinement tests the last row beyond the exhaustive census. The eight-nut word
\[
 C_8=(-1,0,-1,1,-1,-1,1)
\]
has the non-fan heptagon arrangement and $0.09454<\mu(C_8)<0.09541$. Extend its first zero string and then insert a zero pair near the right end:
\[
 \begin{split}
 C_{10}&=(-1,0,0,0,-1,1,-1,-1,1),\\
 C_{12}&=(-1,0,0,0,-1,1,-1,0,0,-1,1).
 \end{split}
\]
The certified homogeneous values change from $\mu(C_{10})=0.07155305867\ldots$ to $\mu(C_{12})=-0.01968371494\ldots$. The twelve-nut word has full ALF balance certificates in both allowed sectors, at $\beta=-4/7$ and $-2/7$. All three words have the same direction arrangement and the same two possible fibre lines. This is a targeted example, not an ALF classification through twelve nuts.

\begin{figure}[htbp]\centering
\includegraphics[width=.72\textwidth]{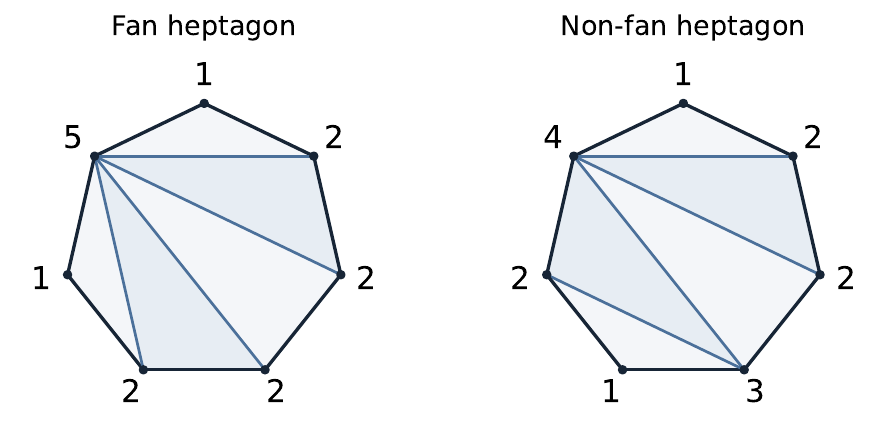}
\caption{The two heptagon triangulations in the ALF census. The right-hand type has no certified balance root through eight nuts but admits the twelve-nut refinement $C_{12}$. Vertex labels are the incidence coefficients from \eqref{eq:polygon-cycle}.}\label{fig:ALF-heptagons}
\end{figure}

Together with the AF zigzag example, this refinement supports studying realizability at the level of marked walks on a polygon. It leaves open which refinements succeed, whether an arrangement always supports a regular walk, and whether the same criterion governs the AF and non-GH ALF sectors.

For the word $C_{12}$ above, the two full-system certificates have the following data. The primitive fibres use coefficient coordinates in $(A,-B)$.
\begin{table}[htbp]\centering\small
\begin{tabular}{@{}rrrrll@{}}\toprule
$k$&$h$&$n=N+h$&$\beta$&$f_\infty$&Contraction bound\\\midrule
5&2&14&$-4/7$&$(4,-7)$&$<2.338\cdot10^{-11}$\\
6&0&12&$-2/7$&$(2,-7)$&$<6.346\cdot10^{-14}$\\\bottomrule
\end{tabular}
\caption{Two isolated positive ALF balance roots with the same twelve-nut word and different asymptotic fibres.}
\end{table}
The computation used 120 decimal digits. Independent evaluation of the original homogeneous BZ determinants at the refined centres gives rod residuals below $1.7\cdot10^{-73}$ in both sectors; exact existence follows from the interval boxes, not from the residuals. The auxiliary negative common-cone model is certified separately. Its rank is five and its cone minimum lies within $4.0\cdot10^{-35}$ of $-0.0196837149414386576970427241280012$.

\section{Supplementary analytic arguments}\label{app:analytic-proofs}
This appendix proves the coupled transport lemma and supplies the local comparison, collision, and insertion arguments used in Section~\ref{sec:analytic}. It also records finite enclosures and derivative formulas used to interpret the numerical catalogue.

\subsection{Transport through singular principal blocks}\label{app:paired-transport}
The following finite-matrix statement supplies \eqref{eq:triangular-transport} without a nonvanishing assumption on the principal minors of $T$.

\begin{lemma}[Coupled triangular transport]\label{lem:paired-transport}
Let $T$ be a real symmetric $N\times N$ matrix and $v\in\R^N$. For their leading blocks, set
\[
 D_j=\det T_j,\qquad A_j=v_{[j]}^T\adj(T_j)v_{[j]},\qquad
 \mathsf r_j=(-1)^j(2A_j,D_j),\qquad \mathsf r_0=(0,1).
\]
Suppose consecutive rods are nonparallel. Every infinitesimal variation $(\delta T,\delta v)$ that preserves all the projective directions $[\mathsf r_j]$, $1\le j\le N$, admits a unique lower triangular $K$ such that
\[
 \delta T=KT+TK^T,\qquad \delta v=Kv,\qquad
 \delta\mathsf r_j=2\tr(K_j)\mathsf r_j.
\]
The matrix $K$ depends real-analytically on the data and the variation locally on this domain, including where $D_j=0$.
\end{lemma}
\begin{proof}
Choose a real scalar $c$, fixed during the variation, and set $Q=T+c vv^T$. The rank-one determinant identity gives
\begin{equation}\label{eq:transport-pencil}
 \det Q_j=D_j+c A_j,\qquad
 \adj(Q_j)v_{[j]}=\adj(T_j)v_{[j]}.
\end{equation}
The second identity follows first for invertible $T_j$ by the rank-one inverse formula and then everywhere by polynomial continuation. Each pair $(A_j,D_j)$ is nonzero, so $c$ can avoid the finitely many zeros of the first expression. All leading minors of $Q$ are then nonzero in a neighborhood of the given data. Differentiating its $LDL^T$ factorization gives a unique lower triangular $K$ with $\delta Q=KQ+QK^T$.

Put $e=\delta v-Kv$. The quadratic adjugate $A_j$ is unchanged by replacing $T_j$ with $Q_j$. Its variation and that of the shifted determinant are therefore
\begin{align*}
 \delta A_j&=2\tr(K_j)A_j
       +2e_{[j]}^T\adj(Q_j)v_{[j]},\\
 \delta(D_j+c A_j)&=2\tr(K_j)(D_j+c A_j).
\end{align*}
The shifted pair $(2A_j,D_j+c A_j)$ is obtained from $(2A_j,D_j)$ by one fixed shear. Its direction is preserved, and its second component is nonzero. Hence
\[
 e_{[j]}^T\adj(Q_j)v_{[j]}=0,\qquad 1\le j\le N.
\]
These equations form a lower triangular system for $e$. Its $j$th diagonal coefficient is
\[
 g_j=\bigl(\adj(T_j)v_{[j]}\bigr)_j,\qquad
 2g_j^2=\det(\mathsf r_{j-1},\mathsf r_j)>0.
\]
The last equality is the polynomial form of \eqref{eq:nut-square}; it also holds when $T_{j-1}$ is singular. Thus $e=0$. Subtracting $c\,\delta(vv^T)$ from the equation for $\delta Q$ proves the asserted equation for $\delta T$. Congruence gives the rod variation. Any second coupled transport would give a second transport of $Q$, proving uniqueness. A fixed admissible $c$ works on a neighborhood, where the construction is analytic. Uniqueness makes these local constructions agree.
\end{proof}

The pencil is only an auxiliary coordinate choice; it is not required to satisfy the involution or displacement identities. Those identities are applied to the recovered transport of the original $(T,v)$. In particular, no limiting argument in $\beta$ is needed at a zero angular component.

\subsection{Canonical identification}\label{app:canonical-proof}
The following proves Proposition~\ref{prop:canonical-comparison}.
\begin{proof}
We first show that $\widehat{\mathcal B}$ in
\eqref{eq:augmented-bulk} is positive definite everywhere on the closed
orbit half-plane. Away from nuts all $R_i$ are positive, so this is immediate.
At the nut $z_i$, a vector annihilated by its quadratic form would have to
be a multiple of the $i$th coordinate vector, and $Te_i$ would also have to
be a multiple of $e_i$. This is impossible. Indeed the weighted polynomial
projector $\mathcal P=(I+T)/2$ satisfies $0<\mathcal P_{ii}<1$:
its $i$th row is nonzero, and deleting that row leaves a full-rank
Vandermonde matrix. Thus $Te_i$ is not a coordinate eigenvector.

On a rod interior, \eqref{eq:regular-gram} is analytic in $\rho^2$ and $z$.
Its axis value is positive semidefinite of rank one. Rank zero would give
$G=O(\rho^2)$ and contradict $\det G=\rho^2$.
In a basis containing the vanishing rod, the Gram matrix therefore has a
positive tangential entry, a normal entry of order $\rho^2$, and mixed entry
$O(\rho^2)$. Consequently $G/\rho$ is at bounded hyperbolic distance from
the corresponding rod model. The adjugate rod formula gives the same
cone-normalized vector by continuity, also when $\det T_j=0$.

Near a nut put $\zeta=z-z_i$, $r=(\rho^2+\zeta^2)^{1/2}$,
$x=r-\zeta$ and $y=r+\zeta$. Thus $xy=\rho^2$ and $x+y=2r$.
All entries of \eqref{eq:regular-gram} are analytic in $(x,y)$.
Let $S=(\mathsf r_{i-1}\ \mathsf r_i)$ and $\delta=\det S>0$.
The two distinct vanishing directions and analyticity give
\[
 S^TGS=
 \begin{pmatrix}
 y A(x,y)&xy E(x,y)\\
 xy E(x,y)&x D(x,y)
 \end{pmatrix},\qquad
 AD-xyE^2=\delta^2.
\]
Here $A(0,0),D(0,0)>0$. Also
$e^{2\nu}=F(x,y)/(x+y)$ with $F$ analytic and positive.
Cone normalization on the two open axes gives
$A(x,0)=F(x,0)$ and $D(0,y)=F(0,y)$.
Hence $A(0,0)=D(0,0)=F(0,0)=\delta$, and
\begin{equation}\label{eq:all-chart-nut}
 \delta^{-1}S^TGS=
 \begin{pmatrix}
 y[1+O(r)]&O(xy)\\
 O(xy)&x[1+O(r)]
 \end{pmatrix},\qquad
 e^{2\nu}=\frac{\delta}{2r}[1+O(r)].
\end{equation}
The relative Gram error is $O(r)$ uniformly up to either axis.
This proves bounded-distance comparison with the nut model without a
principal-minor restriction.

At infinity, $\mathcal B=2rI+O(1)$.
Expansion of \eqref{eq:positive-bulk}, using $u^Tv=0$ and the exact
vanishing of the rotation potential on both exterior rods, gives
\[
 f=1-\frac{2M}{r}+O(r^{-2}),\qquad
 \omega=O(\rho^2/r^3),\qquad
 e^{2k}=1+O(\rho^2/r^4).
\]
These estimates are uniform in angle; they follow by expanding in $r^{-1}$
and factoring the exterior-axis zeros. They imply hyperbolic distance
$O(r^{-1})$ to the prescribed AF end.
The bounded-distance uniqueness argument of
\cite[Theorem~4.24]{LiSun} now identifies the harmonic maps.
The augmentation equations determine $\nu$ up to a constant, and the
normalization at infinity fixes that constant
\cite[Proposition~4.26]{LiSun}. Finally, Proposition~\ref{prop:direction-coordinates}
supplies the direction-matched data at every positive length vector.
\end{proof}

For Corollary~\ref{cor:smooth-completion}, the completion argument is as follows.
\begin{proof}
The rod and nut estimates above put the metric in the scope of the local
regularity results \cite[Propositions~4.15 and~4.23]{LiSun}. Cone balance
removes the rod singularities, and the lattice-basis condition removes
orbifold quotients at nuts. The normalized AF expansion gives a complete end
and finite curvature energy.
\end{proof}

\subsection{The homogeneous moment model}\label{app:homogeneous-model}
Let a cluster have $r$ distinct ordered nodes $z_1<\cdots<z_r$, positive weights $\omega_i$, and lower occupation degree $k$. The latter is chosen to realize its prescribed signed direction lift. Define the signed moments and their Hankel determinants by
\begin{equation}\label{eq:cluster-moments}
 H_a^{[j]}=\det\left(\sum_{i=1}^{r}\epsilon_i^{[j]}\omega_i z_i^{s+t}\right)_{s,t=0}^{a-1},\qquad
 H_0^{[j]}=1,\qquad \mathsf t_j=(H_{k+1}^{[j]},H_k^{[j]}),
\end{equation}
where $\epsilon_i^{[j]}=-1$ for $i\le j$ and $+1$ otherwise. A single real linear map sends the two independent exterior vectors to the prescribed primitive exterior generators. Its action on every internal vector defines the multipliers $a_j>0$ by $\mathsf r_j=a_j\mathsf v_j$. This fixes $b_{C,j}=\log a_j$, including the exterior normalization. The degree $k$ in this auxiliary model need not equal the reduced degree of a finite ALF realization of the same word.

For completeness, the homogeneous first law and convexity have a direct moment interpretation. Let $Z=\operatorname{diag}(z_i)$ and let $P_a$ be the orthogonal projector onto the vectors $(\sqrt{\omega_i}p(z_i))_i$ with $\deg p<a$. Multiplication by $z$ gives
\[
 P_kZ(I-P_{k+1})=0.
\]
Write $T_-=I-2P_k$, $T_+=I-2P_{k+1}$ and $P_{k+1}-P_k=qq^T$ with $\|q\|=1$. The matrix determinant lemma gives $H_a^{[j]}/H_a^{[0]}=\det(I-2P_a)_{[j]}$. After a constant change of the two rod coordinates, the direction data are therefore the pairs of leading determinants of $T_+$ and $T_-$. Since $T_+=T_--2qq^T$, these pairs are a fixed linear transform of the rod pairs for $(T_-,q)$ in Lemma~\ref{lem:paired-transport}. The exterior frame is held fixed here, as in the definition of the direction-matched model. The lemma gives $\delta T_-=KT_-+T_-K^T$ and $\delta q=Kq$, and hence also $\delta T_+=KT_++T_+K^T$. This proves simultaneous transport even when a leading minor vanishes. Congruence gives $\delta b_j=2\operatorname{tr}K_{[j]}$, where $K_{[j]}$ is the leading $j\times j$ block. With $S=K+K^T$, differentiating the two involutions gives $ST_\pm+T_\pm S=0$. Their difference, together with $\|q\|=1$, gives $Sq=0$. Writing $\xi_j=z_{j+1}-z_j$ for the finite gaps, the polynomial multiplication identity and its derivative give
\[
 \operatorname{tr}(ZS)=0,\qquad \sum_j \xi_j\,\delta b_j=0,\qquad
 \delta\xi^T(D_\xi b_C)\delta\xi=\sum_{i>j}(z_i-z_j)K_{ij}^2\ge0.
\]
These identities establish the homogeneous potential and its convexity on the positive direction-matched family. The local direction-coordinate hypothesis in Section~\ref{subsec:comparison} specifies the family on which the transport is used. In rank zero the potential can be linear; no strict convexity assertion is needed for the sign tests.

\subsection{Proof of the collision criterion}\label{app:collision-proof}
Parallel blocks require a bound preventing other cone errors from canceling the divergent negative error. The involution and AF residue identity give
\[
 \|v\|^2\le S:=z_N-z_1,\qquad
 |\det T_j|\le1,\qquad \|\adj(T_j)\|\le1.
\]
Equation~\eqref{eq:T-rods} implies
\begin{equation}\label{eq:bounded-span-cones}
 |\mathsf r_j^t|\le2S,\qquad |\mathsf r_j^\phi|\le1.
\end{equation}
Since $\mathsf r_j=e^{b_j}(-q_j,p_j-\beta q_j)$, every $b_j$ is bounded above on a bounded-span family. If $q_j=0$, primitivity gives $|p_j|=1$ and hence $b_j\le0$; otherwise $b_j\le\log(2S/|q_j|)$. A parallel block has at least one internal direction different from its exterior direction, so its diverging cone gives $b_j\to-\infty$. Any fixed strictly positive opening of that block has cost $-\infty$. Complete contraction is also covered by \eqref{eq:bounded-span-cones}, since $S\to0$.

\begin{proof}
At a finite stratum minimum, open each nonparallel block in a shape of negative cost and assign the blocks positive total-length fractions. Equation~\eqref{eq:opening-cost} gives a negative total derivative. If a parallel block occurs, \eqref{eq:bounded-span-cones} instead gives derivative $-\infty$. The full AF exterior pair is parallel, so every nonparallel collision block is proper. Continuity of a proper convex function along a segment from a finite boundary point to an interior point converts either derivative into strict decrease. No collision face can contain a minimum.

For the upper side, suppose a sequence has total span tending to infinity. Rescale to span one and pass to a subsequence. At least one rod retains a positive fraction of the span. The large-span cone limit \cite[Proposition~5.15]{LiSun}, together with $\max_j|p_j-\beta q_j|<1$, supplies a retained rod with limiting angle below $2\pi$. If a retained angular component vanishes, its limiting angle is zero. Otherwise all retained limiting directions agree with the exterior direction. In either case, some length tends to infinity with $b_j>0$.

Minimize $\Phi_W$ on $[\delta,R]^{N-1}$, where $0<\delta\le1<R$. At a minimizer, $x_j>\delta$ implies $b_j\le0$, because shortening is feasible. The preceding large-span conclusion therefore bounds all such minimizers above by a constant $M(W,\beta)$ independent of $\delta,R$. Fix $R>M+1$ and let $\delta\downarrow0$. A subsequence converges to $q\in[0,M]^{N-1}$. An affine supporting plane bounds the potential below on this compact set, while comparison with any fixed interior point bounds it above. Lower semicontinuity and minimality give
\[
 \overline\Phi_W(q)\le\liminf\Phi_W(x_\delta)\le\Phi_W(y)
 \quad\text{for every }y\in(0,R)^{N-1}.
\]
Thus $q$ minimizes the closed-box extension. Boundary descent forces $q$ to be interior, hence $b(q)=0$. Strict convexity gives uniqueness; the positive Hessian and the implicit-function theorem give analytic dependence. Corollary~\ref{cor:smooth-completion} gives the final metric statement.
\end{proof}

\subsection{Negative opening tests}\label{app:negative-tests}
Three elementary tests suffice for the families below. First, write the internal rods in the real exterior basis,
\[
 \mathsf v_j=\alpha_j\mathsf v_0+\gamma_j\mathsf v_n,
 \qquad D(C)=\max_{0<j<n}(|\alpha_j|+|\gamma_j|).
\]
Cauchy--Binet expresses $H_m^{[j]}/H_m^{[0]}$ as a signed average of positive weights, so its absolute value is at most one. The exterior raw components, after separate component normalization, are $(1,1)$ and $((-1)^{k+1},(-1)^k)$. Consequently
\begin{equation}\label{eq:diamond-bound}
 b_{C,j}\le-\log(|\alpha_j|+|\gamma_j|),\qquad
 \mu(C)\le-\log D(C).
\end{equation}
The second inequality follows by concentrating positive normalized lengths near an index maximizing $D$. Geometrically, $D>1$ means that an internal primitive rod lies outside the diamond spanned by the two signed exterior rods; it guarantees $\mu<0$.

Second, deletion of an entry $d_i=\sigma\in\{-1,1\}$ changes each existing neighbor by $d\mapsto d-\sigma$. This operation has a positive local inverse. For insertion of a $-1$ rod, split a moment mass $W$ into $sW$ and $(1-s)W$ at a common node. Rank-one affinity of both determinants gives the intermediate raw rod
\[
 \mathsf t_{\rm new}=(1-s)\mathsf t_-+s\mathsf t_+.
\]
Since $\mathsf v_{\rm new}=\mathsf v_-+\mathsf v_+$, direction matching uniquely determines
\begin{equation}\label{eq:single-insertion}
 s=\frac{a_-}{a_-+a_+},\qquad
 a_{\rm new}=\frac{a_-a_+}{a_-+a_+}.
\end{equation}
The new direction equation has nonzero derivative in $s$ and the old equations are independent of $s$ at collision. Local invertibility of the old direction map therefore lifts the construction to positive separation. The $+1$ case follows by complementing occupation subsets:
\begin{equation}\label{eq:particle-hole}
 H_m^{[j]}(\omega)=(-1)^j\mathcal A\,
 H_{n-m}^{[j]}(\widetilde\omega),\qquad
 \widetilde\omega_i=\frac1{\omega_iD'(z_i)^2},\quad
 \mathcal A=\Delta(z)^2\prod_i\omega_i,
\end{equation}
where $D(z)=\prod_i(z-z_i)$ and $\Delta(z)=\prod_{i<j}(z_j-z_i)$. The induced frame change reverses the self-intersection signs and preserves cone multipliers. Vanishing inserted length gives
\[
 \mu(C)\le\mu(C_{\rm deleted}).
\]
A nonempty Gibbons--Hawking chain $(-2,\ldots,-2)$ has all $b_j=0$ in the homogeneous exterior normalization. Its first inverse deletion has opening derivative $-\log2$, including at an exterior nut; subsequent insertions preserve negativity. Thus any word reducing by at least one such deletion to a nonempty Gibbons--Hawking chain has $\mu<0$.

Third, an even-length all-zero word has a negative Weyl trial. For $n=2k+1$ nodes take $\omega_i=1/|D'(z_i)|$. Then \eqref{eq:particle-hole} gives $H_{k+1}^{[j]}=(-1)^jH_k^{[j]}$, so the rod directions alternate between the two exterior directions. Strict positivity of the occupation sums gives $0<|H_k^{[j]}|/H_k^{[0]}<1$. Every finite cone multiplier is therefore below one. For two finite rods the exact formula is
\begin{equation}\label{eq:zero-entropy}
 \Psi_{(0,0)}(s,1-s)=s\log s+(1-s)\log(1-s),\qquad 0<s<1.
\end{equation}
Odd-length all-zero words have parallel exterior rods and belong to the separate descent argument.

\subsection{The confluent insertion calculation}\label{app:insertion-proof}
This proves Lemma~\ref{lem:zero-insertion}; $D$ denotes the unrefined word in this calculation.
\begin{proof}
Use positive occupation variables $\zeta_i=\omega_i|D'(z_i)|$. Split a node $(c,\zeta_*)$ into three at $c-\varepsilon A,c,c+\varepsilon B$, with $A,B>0$ and all three occupation variables initially equal to $\zeta_*$. If the old moment mass at $c$ is $W$, the leading new masses are $\varepsilon^{-2}\lambda_i$, where
\[
 (\lambda_1,\lambda_2,\lambda_3)
 =W\left(\frac1{A(A+B)},\frac1{AB},\frac1{B(A+B)}\right).
\]
An outside mass $\omega_i$ tends to $\omega_i/(z_i-c)^2$. For $j=0,1,2,3$ passed cluster nodes, put $\xi=(-A,0,B)$ and
\[
 M_j=\sum_i\epsilon_i^{[j]}\lambda_i,\quad
 S_j=\sum_i\epsilon_i^{[j]}\lambda_i\xi_i,\quad
 Q_j=\sum_i\epsilon_i^{[j]}\lambda_i\xi_i^2,\quad
 V_j=Q_j-S_j^2/M_j.
\]
Direct calculation gives
\begin{equation}\label{eq:confluent-masses}
 (V_0,V_1,V_2,V_3)=(W,-W,W,-W),\qquad
 \frac{(M_0,M_1,M_2,M_3)}{M_0}=(1,s,-(1-s),-1),
 \quad s=\frac A{A+B}.
\end{equation}
Use the monic basis $1,(z-c),(z-c)z,\ldots,(z-c)z^{m-1}$ in the new Gram determinant. Its first entry is $\varepsilon^{-2}M_j$; the cross block is $\varepsilon^{-1}S_j(1,c,\ldots,c^{m-1})+O(1)$. The Schur complement therefore replaces the cluster by effective mass $V_j$, giving
\begin{equation}\label{eq:confluent-Hankel}
 \varepsilon^2H_{m+1}^{{\rm new},[j]}
 \longrightarrow M_jH_m^{\rm old}(V_j).
\end{equation}
The common factor is the same in the two consecutive degrees. After exterior normalization, the surviving rods converge to the old rods, with the signs in \eqref{eq:zero-rod-operation}, and the two new rods tend to $s\mathsf r_+$ and $-(1-s)\mathsf r_-$. This proves the limiting multipliers.

To obtain exact directions at positive separation, use the effective mass $V_0$ and two relative cluster weights as coordinates. The surviving equations at $\varepsilon=0$ are the old equations; the new ones are $V_1=-V_0$, $V_2=V_0$. With the middle leading weight set to one and the other two denoted $x,y$, their Jacobian at $x=B/(A+B)$, $y=A/(A+B)$ is
\[
 \frac{\partial(V_1+V_0,V_2-V_0)}{\partial(x,y)}
 =\begin{pmatrix}-B(2A+B)&B^2\\-A^2&A(A+2B)\end{pmatrix},
 \qquad\det=-2AB(A+B)^2\ne0.
\]
Fixing $V_0=W>0$ preserves nonsingularity. Independence of the two old adjacent directions makes these effective-mass equations equivalent to direction matching. After removal of the common pole, the full system is analytic with block triangular invertible Jacobian. The implicit-function theorem supplies positive exact solutions, with all three corrected occupation variables tending to $\zeta_*$.

The new errors have finite limits and their lengths vanish, so $\Psi_{D^\sharp}(h_\varepsilon)\to\Psi_D(h)$. Taking infima proves the inequality. For the symmetric interior opening, the first law in \eqref{eq:cluster-potential} gives
\[
 \Psi'_{D^\sharp}(0+)
 =\tfrac12\log(a_+/2)+\tfrac12\log(a_-/2)
   -\tfrac12\log a_--\tfrac12\log a_+=-\log2.
\]
Integration gives \eqref{eq:zero-opening-derivative}. The remainder is local in the original data; no uniform finite-opening estimate is required.
\end{proof}

\subsection{Verification of the explicit families}\label{app:family-proof}
\begin{proof}
For the five base words at $r=s=0$, the rod recurrence verifies the stated windows and sectors. Every nonparallel proper block either has $D(C)>1$, reduces by elementary $\pm1$ deletions to a nonempty Gibbons--Hawking chain, or is an even all-zero word. These finite integer checks and their reduction paths are included in the supplementary verification. Appendix~\ref{app:negative-tests} gives their signs, and Theorem~\ref{thm:zero-extension} proves the result for all $r,s$.

For the second family there is also a direct description at every order. In lattice coefficient coordinates its finite unoriented directions are $(0,1),(-1,1),(1,-2)$, while the exterior is $(1,0)$. The angular components, up to sign, are $-\beta,-1-\beta,1+2\beta$, giving exactly $(-1,0)$. The sole determinant of magnitude two pairs an exterior direction with $(1,-2)$, at separation at least three.

Its proper blocks, up to reversal and orientation, are all-zero blocks, a single end entry $\pm1$ with zeros, $(-1,1,0^{2a})$, $(-1,1,0^{2a+1})$, $(1,0^{2r+1},-1)$, and $(-1,1,0^{2r+1},-1)$. The middle two-sided block is parallel; the one-end blocks and the remaining odd-string blocks have diamond constants two or three. Writing $B_a=(-1,1,0^{2a})$, the only repeating reduction is
\[
 B_a\longrightarrow(-2,-1,0^{2a-1})\longrightarrow B_{a-1},
 \qquad B_0=(-1,1)\longrightarrow(-2).
\]
This proves all its block signs directly.
\end{proof}

\subsection{Finite interior certificates}\label{app:finite-certificates}
A complementary finite certificate gives effective enclosures of balanced lengths without separate collision calculations. Normalize $L=1$ and write $x=\ell/L$ and $\Phi=\Phi_W$. All following derivatives use these fixed-period coordinates; in particular $H=D_xb$. Choose finitely many positive trial vectors $q_0,q_1,\ldots,q_s$ and evaluate only their cone gradients $g_k=b(q_k;\beta)$. None of the trial metrics needs to be balanced.

Two supporting-plane inequalities imply
\begin{equation}\label{eq:supporting-plane}
 \Phi(x)-\Phi(q_0)
 \ge g_k\cdot(x-q_0)-(g_k-g_0)\cdot(q_k-q_0).
\end{equation}
Define
\begin{equation}\label{eq:certificate-polyhedron}
 \mathcal P=\left\{x\in\R^{N-1}:
 g_k\cdot(x-q_0)\le(g_k-g_0)\cdot(q_k-q_0)\quad\text{for all }k\right\}.
\end{equation}

\begin{theorem}[Finite interior certificate]\label{thm:finite-certificate}
If $\mathcal P$ is bounded and contained strictly inside the positive orthant, then there exists exactly one balanced length vector. It belongs to $\mathcal P$, so every component satisfies
\[
 x_j\ge\min_{y\in\mathcal P}y_j>0.
\]
If the same conditions hold uniformly over a twist interval, the bounds and the existence conclusion hold throughout that interval.
\end{theorem}
\begin{proof}
Equation~\eqref{eq:supporting-plane} places the nonempty sublevel set $\{\Phi\le\Phi(q_0)\}$ inside $\mathcal P$. It is consequently compactly contained in the domain. The potential attains its infimum at an interior point, where its gradient vanishes. Strict convexity gives uniqueness.
\end{proof}

The argument controls arbitrary ratios of collapsing lengths and arbitrary sizes of the other lengths. It can also give a cone-error estimate away from balance. For example, let $E$ be a positive diagonal matrix specifying the side lengths of a box, and let $\alpha>0$ and $B_0\ge0$ be constants. Suppose the supporting planes yield
\[
 \Phi(x)-\Phi(q_0)\ge\alpha\|E^{-1}(x-q_0)\|_\infty-B_0,
 \qquad\alpha>0,
\]
and the box $q_0+E[-B_0/\alpha,B_0/\alpha]^{N-1}$ lies strictly inside the orthant. Convexity at $x$ then converts distance outside that box into a positive lower bound for $\|b(x)\|_\infty$. This supplies the desired persistent cone mismatch near every collision face. The appendix gives two verified examples.

This criterion is useful when explicit quantitative bounds are required. The collision argument proves existence for the families above without first constructing such an enclosure; the supporting planes can subsequently yield effective bounds on selected parameter intervals.

\subsection{Twist derivatives}\label{app:twist-derivatives}
At fixed $\beta$, metrics with prescribed cone-error vector $c$ are critical points of $\Phi(x)-c\cdot x$, so their local response is $\delta x=H^{-1}\delta c$. The flow $\dot x=-H^{-1}b$ satisfies
\[
 \dot b=-b,\qquad \frac{d\Phi}{dt}=-b^TH^{-1}b.
\]
It balances nearby cones and converges locally around a minimum. Global convergence still requires the trajectory to remain in the domain. This is finite-dimensional cone stability, rather than stability under arbitrary metric perturbations.

Changing $\beta$ by an increment $\tau$ at $L=1$ shears a rod vector by $(\mathsf r^t,\mathsf r^\phi)\mapsto(\mathsf r^t,\mathsf r^\phi+\tau\mathsf r^t)$. In \eqref{eq:T-rods}, the shear is represented exactly by
\[
 (T,v)\mapsto(T+2\tau vv^T,v).
\]
The physical deformation restores $T^2=I$ by triangular transport. Along a balanced branch, with primes denoting $d/d\beta$,
\begin{equation}\label{eq:twist-transport}
 T'=KT+TK^T+2vv^T,\qquad v'=Kv,\qquad u'=-K^Tu,
 \qquad K_{ii}=0.
\end{equation}
The same trace argument as above gives the extended first law
\begin{equation}\label{eq:full-first-law}
 2\,dM=x\cdot db+4J\,d\beta,\qquad
 d\Phi=b\cdot dx-4J\,d\beta.
\end{equation}
For arbitrary $L$, define the dimensionless invariants $\overline M=M/L$ and $\overline J=J/L^2$. A balanced branch at fixed $L$ therefore satisfies $d\overline M/d\beta=2\overline J$.

The relative length response is
\begin{equation}\label{eq:relative-response}
 R=\diag(x)^{-1}H^{-1}\partial_\beta b
       =-\frac{d}{d\beta}\log x.
\end{equation}
The BZ representation refines this generic implicit-derivative formula. In the positive polynomial row gauge, $T_{j,j+1}>0$ by Cauchy--Binet, and the displacement identity gives
\begin{equation}\label{eq:positive-gap-quotient}
 x_j=\frac{v_j u_{j+1}-u_jv_{j+1}}{T_{j,j+1}},\qquad
 v_j u_{j+1}-u_jv_{j+1}>0.
\end{equation}
Thus $R_j$ is the difference of logarithmic derivatives of two positive quantities. For initial boundary data $(T_*,u_*,v_*)$ satisfying \eqref{eq:T-identities} and $u_*^Tv_*=0$, the finite bare shear is explicitly invertible:
\[
 (T_*+2\tau v_*v_*^T)^{-1}=T_*-2\tau u_*u_*^T.
\]
These identities isolate concrete objects to estimate. A quantitative response estimate additionally requires control of the compensating triangular transport and the relative node motion.

If a branch has candidate interval $(\beta_-,\beta_+)$, define the endpoint weight and its associated parameter by
\[
 w(\beta)=\frac{(\beta-\beta_-)(\beta_+-\beta)}{\beta_+-\beta_-},\qquad
 s(\beta)=\log\frac{\beta-\beta_-}{\beta_+-\beta}.
\]
A bound $\|wR\|_\infty\le C$ along its existing part would imply
\begin{equation}\label{eq:integrated-response}
 x_j(\beta_0)e^{-C|s(\beta)-s(\beta_0)|}
 \le x_j(\beta)\le
 x_j(\beta_0)e^{C|s(\beta)-s(\beta_0)|}.
\end{equation}
These are uniform upper and lower length bounds on compact subintervals. These estimates would provide quantitative continuation bounds. For words satisfying Theorem~\ref{thm:collision-existence}, interior continuation already follows from existence and uniqueness; a bound on the relative response would additionally control its rate. No uniform estimate of this kind at the twist endpoints is asserted.

\bibliographystyle{amsplain}
\providecommand{\bysame}{\leavevmode\hbox to3em{\hrulefill}\thinspace}
\providecommand{\MR}{\relax\ifhmode\unskip\space\fi MR }
\providecommand{\MRhref}[2]{%
  \href{http://www.ams.org/mathscinet-getitem?mr=#1}{#2}
}
\providecommand{\href}[2]{#2}

\end{document}